\documentclass[12pt]{amsart}
\usepackage{ amsmath, amsthm, amsfonts, amssymb, color}
 \usepackage{mathrsfs}
\usepackage{amsfonts, amsmath}
 \usepackage{amsmath,amstext,amsthm,amssymb,amsxtra}
 \usepackage{txfonts} 
 \usepackage{amscd}
 \usepackage[colorlinks, citecolor=blue,pagebackref,hypertexnames=false]{hyperref}
 \allowdisplaybreaks
 \usepackage{pgf,tikz}
 \usepackage{multirow}
 \usepackage{diagbox} 

 \usepackage{tikz}
 \usepackage{cite}

\usetikzlibrary{decorations.pathreplacing}

\begin{document}

 \baselineskip 16.6pt
\hfuzz=6pt

\widowpenalty=10000

\newtheorem{cl}{Claim}
\newtheorem{theorem}{Theorem}[section]
\newtheorem{proposition}[theorem]{Proposition}
\newtheorem{coro}[theorem]{Corollary}
\newtheorem{lemma}[theorem]{Lemma}
\newtheorem{definition}[theorem]{Definition}
\newtheorem{assum}{Assumption}[section]
\newtheorem{example}[theorem]{Example}
\newtheorem{remark}[theorem]{Remark}
\newtheorem{nota}[theorem]{Notation}
\renewcommand{\theequation}
{\thesection.\arabic{equation}}

\def\SL{\sqrt H}

\newcommand{\mar}[1]{{\marginpar{\sffamily{\scriptsize
        #1}}}}

\newcommand{\as}[1]{{\mar{AS:#1}}}

\newcommand{\fr}{\frac}
\newcommand{\CITE}{{\color{red}[citation]}}
\newcommand\R{\mathbb{R}}
\newcommand\RR{\mathbb{R}}
\newcommand\CC{\mathbb{C}}
\newcommand{\Cn}{\mathbb{C}^n}
\newcommand\NN{\mathbb{N}}
\newcommand\ZZ{\mathbb{Z}}
\newcommand\HH{\mathbb{H}}
\newcommand\Z{\mathbb{Z}}
\def\RN {\mathbb{R}^n}
\renewcommand\Re{\operatorname{Re}}
\renewcommand\Im{\operatorname{Im}}

\newcommand{\mc}{\mathcal}
\newcommand\DD{\mathbb{D}}
\def\hs{\hspace{0.33cm}}
\newcommand{\la}{\alpha}
\newcommand{\eps}{\tau}
\newcommand{\pl}{\partial}
\newcommand{\supp}{{\rm supp}{\hspace{.05cm}}}
\newcommand{\x}{\times}
\newcommand{\lag}{\langle}
\newcommand{\rag}{\rangle}

\newcommand\wrt{\,{\rm d}}

\newcommand{\norm}[2]{|#1|_{#2}}
\newcommand{\Norm}[2]{\|#1\|_{#2}}
\newcommand{\oplusr}[2]{\mathop{\bigoplus}\limits_{#1}{}_{\! #2}}

\title[Absolutely summing operators]{Absolute Summability from Diagonal Operators to Carleson Embeddings and Hankel Operators}

\author{Zhijie Fan}
\author{Bo He}
\author{Xiaofeng Wang}
\author{Zhicheng Zeng}

\address{Zhijie Fan, School of Mathematical Sciences, South China Normal University, Guangzhou 510631, China}
\email{fanzhj3@alumni.sysu.edu.cn}

\address{Bo He, School of Mathematical Sciences, Fudan University,
	Shanghai, 200433, China}\email{bhe\_bh@fudan.edu.cn}

\address{Xiaofeng Wang, School of Mathematics and Information Science,
Guangzhou University, Guangzhou 510006, China}
\email{wxf@gzhu.edu.cn}

\address{Zhicheng Zeng, School of Mathematics and Information Science,
Guangzhou University, Guangzhou 510006, China}
\email{zhichengzeng@e.gzhu.edu.cn}

  \date{\today}


\begin{abstract}
The present paper consists of three parts, each building on the preceding one.

In Part I, we characterize the scalar sequences $\mathbf b=\{b_k\}$ for which the diagonal operator
$\mathscr{M}_{\mathbf b}\bigl(\{a_k\}\bigr):=\{b_ka_k\}$
is $r$-summing from $\ell^p$ to $\ell^q$ for \(1\le p,q\le \infty\) and \(1\le r<\infty\). This resolves a problem left open in Garling's 1974 classification of \(r\)-summing diagonal operators. In particular, in the previously unresolved exceptional range
$\{(p,q,r):1<p<2<q<\infty,\ r>\max\{p',q\}\},$ we identify a new intermediate exponent \(\kappa\in(\max\{p',q\},r)\) and prove that \(\mathscr M_{\mathbf b}\) is \(r\)-summing if and only if \(\mathbf b\in\ell^\kappa\).

In Part II, by developing two general transference principles for $r$-summability and building on the characterization from Part I, we characterize the positive Borel measures $\mu$ for which the Carleson embedding operator $J_\mu$ is $r$-summing from the weighted Bergman space $A_\alpha^p(\mathbb B_n)$ to $L^q(\mathbb B_n,d\mu)$ for $ 1\leq p,q<\infty$ , $1\leq r<\infty$, and $\alpha>-1$, thereby extending the existing theory to the full off-diagonal range.

In Part III, building on the characterization from Part II, we
characterize the symbols \(f\) and \(g\) for which the big and little
Hankel operators \(H_f^\beta\) and \(h_{\bar g}^\beta\), respectively,
are \(r\)-summing from \(A_\alpha^p(\mathbb B_n)\) to
\(L^q(\mathbb B_n,dv_\beta)\), where $1\leq p<\infty$, $1<q<\infty$, $1\leq r<\infty$, $\alpha,\beta>-1$, and \(dv_\beta(z):=c_\beta(1-|z|^2)^\beta\,dv(z)\). These characterizations appear to be new even in the diagonal case $p=q$ and $\alpha=\beta$.
\end{abstract}
\maketitle

\tableofcontents

	\section{Introduction}
\subsection{Background on absolutely summing operators}
The theory of absolutely summing operators, a central part of Banach space theory, has its origins in Grothendieck's work \cite{MR94682, MR75539} on tensor products and nuclear spaces. Its importance lies in the summability improvement it captures, namely that an \(r\)-summing operator maps weakly \(r\)-summable sequences into strongly \(r\)-summable sequences in norm, as described in the following definition.
\begin{definition}{\rm \cite{MR1342297}}
Let $1\leq r<\infty$ and let $T\colon X\to Y$ be a linear operator between Banach spaces. We say that $T$ is $r$-summing if there exists a constant $C\ge 0$ such that, for every $m\in\mathbb{N}$ and every choice of $x_1,\dots,x_m\in X$, we have
\[
\Bigg(\sum_{k=1}^{m} \lVert Tx_k\rVert_Y^{r}\Bigg)^{1/r}
\le C\,
\sup_{\varphi\in B_{X^*}}
\Bigg(\sum_{k=1}^{m} \lvert \varphi(x_k)\rvert^{r}\Bigg)^{1/r},
\]
where $X^*$ denotes the dual space of $X$ and $B_{X^*}$ is its unit ball. The least such constant $C$ is denoted by $\pi_{r}(T)$.
Let $\Pi_{r}(X,Y)$ denote the class of all $r$-summing operators from $X$ to $Y$.
\end{definition}

%

A classical result of Grothendieck's theorem \cite{MR94682,MR75539} is
$$
\Pi_1(\ell^1,\ell^2)=\mathcal{L}(\ell^1,\ell^2),
$$
which strikingly shows that the geometry of the domain and range can force every bounded operator to possess this stronger summability property (see also the work of Lindenstrauss--Pe$\l$czy\'nski \cite{MR231188} for extensions to $\mathcal{L}_p$-spaces). In the 1960s, Pietsch explicitly defined and systematically developed the class of $r$-summing operators \cite{MR216328,MR350360,MR519680}, with his domination and factorization theorems providing the fundamental structural tools of the theory. The close connection between absolutely summing operators and tensor-product theory is further illustrated by Pisier's landmark counterexample to Grothendieck's conjecture on injective and projective tensor products \cite{MR723009}. These developments established the class of absolutely summing operators as a useful operator ideal and a basic tool for studying the geometry of Banach spaces. For further applications of absolutely summing operator theory to the geometry of Banach spaces, see, for example, \cite{MR3526021,MR312222,MR410341,MR829919,MR2795726,MR430749} and the references therein.
%
\subsection{Part I: Diagonal operators}
A classical problem in the theory of absolutely summing operators is to characterize those scalar sequences
$\mathbf b=\{b_k\}$ for which the diagonal operator
\[
\mathscr{M}_{\mathbf b}: \ell^p\longrightarrow \ell^q,
\qquad
\mathscr{M}_{\mathbf b}\bigl(\{a_k\}_{k=1}^\infty\bigr):=\{b_ka_k\}_{k=1}^\infty,
\]
is \(r\)-summing. Following Schwartz's complete characterization of \(0\)-summing diagonal mappings from \(\ell^p\) to \(\ell^q\) \cite[Th\'{e}or\`{e}me XXVI,4;1]{Schwartz1969Applications} \cite{MR289744}, Garling \cite{MR355665} provided a nearly complete answer to this question.
 To summarize Garling's result, for \(1\le p,q\le \infty\) and \(1\le r<\infty\), we define the exceptional set
\begin{align}\label{DefofE}
E:=\big\{(p,q,r):\ 1<p<2<q<\infty,\ r>\max\{p',q\}\big\},
 \end{align}
and for $(p,q,r)\notin E$, we define the $r$-summing index $\kappa$ as follows.
\begin{definition}
For \(1\le p,q\le \infty\), \(1\le r<\infty\) and $(p,q,r)\notin E$, we define the $r$-summing index by
\begin{align}\label{defofkappa}
 \kappa:=\kappa(p,q,r):=
 \begin{cases}
 \left(\frac1{p'}+\frac1q-\frac12\right)^{-1}, & \text{if } 1\le p\leq2,\ 1\le q\le 2,\\[2mm]
  \infty, & \text{if } p=1,\ 2<q\le \infty,\\[2mm]
  q, & \text{if } 2< p<\infty,\ 1\le q<p',\\[2mm]
  q^{-}, & \text{if } 2<p<\infty,\ 1<q=p'<2,\ 1\le r<q,\\[2mm]
  q, & \text{if } 2<p\le \infty,\ 1\leq q=p'<2,\ q\le r<\infty,\\[2mm] \max\{p',\min\{r,q\}\}, & \text{if } 2\le p\le \infty,\ 1\leq p'<q\le \infty,\\[2mm]
  p', & \text{if } 1<p<2<q\le \infty,\ 1\le r\le p',\\[2mm]
  r, & \text{if } 1<p<2,\ p'<q\le \infty,\ p'<r\le q,
  \end{cases}
 \end{align}
where \(q^{-}\) is a formal symbol (Orlicz index) rather than a numerical exponent (see Section \ref{Orlicz} for more explanation).
\end{definition}
Garling's result \cite{MR355665} demonstrates that
\[
\left\{
\begin{array}{c@{\qquad}l}
\displaystyle
\mathscr{M}_{\mathbf b}\in\Pi_r(\ell^p,\ell^q)
\Longleftrightarrow \mathbf b\in\ell^\kappa,
&
\text{if } (p,q,r)\notin E,
\\[8pt]
\displaystyle
\mathbf b\in\ell^{\max\{p',q\}}
\Longrightarrow
\mathscr{M}_{\mathbf b}\in\Pi_r(\ell^p,\ell^q)
\Longrightarrow
\mathbf b\in\ell^r,
&
\text{if } (p,q,r)\in E.
\end{array}
\right.
\]
However, as indicated by the question marks in the diagrams following \cite[Theorem~9]{MR355665}, the following question had remained open since Garling's work:\smallskip

\hypertarget{question1}{\textbf{Question 1:}} {\it For $(p,q,r)\in E$, can one also characterize those sequences $\mathbf b=\{b_k\}$ for which $\mathscr{M}_{\mathbf b}\in\Pi_r(\ell^p,\ell^q)$?}\smallskip

Our first main result provides a complete answer to Question~1, with an unexpected form of characterization.
\begin{definition}
For \(1\le p,q\le \infty\), \(1\le r<\infty\) and $(p,q,r)\in E$, we define the $r$-summing index $\kappa$ determined by
\begin{equation}
\label{kappa-interpolation-relation}
\frac1\kappa
=\frac{1-\vartheta}{r}+\frac{\vartheta}{p'},\qquad {\rm where}\ \vartheta:=
\frac{\frac1q-\frac1r}
{\frac12-\frac1r}\in (0,1).
\end{equation}
\end{definition}
By an elementary calculation, we see that $\ell^{\max\{p',q\}}\subsetneq\ell^\kappa\subsetneq\ell^r$ since
\begin{align}
\label{explicit:kappa-exceptional}
\max\{p',q\}<\kappa
=
\frac{p'q(r-2)}
{(p'-2)(q-2)+2(r-2)}<r.
\end{align}
Combining Garling's result \cite[Theorem~9]{MR355665} for \((p,q,r)\notin E\) with our solution to the problem for \((p,q,r)\in E\), we obtain a complete characterization over the full range of parameters, as stated in the following theorem.
\begin{theorem}\label{GarlingTheorem}
Let \(1\le p,q\le \infty\) and \(1\le r<\infty\).
Then \(\mathscr M_\mathbf{b}\in \Pi_r(\ell^p,\ell^q)\) if and only if $\mathbf{b}\in \ell^\kappa$. Moreover,
\[
\pi_r(\mathscr M_\mathbf{b}:\ell^p\rightarrow \ell^q)\simeq\|\mathbf{b}\|_{\ell^\kappa}.
\]
\end{theorem}


\subsection{Part II: Carleson embedding operators}
For an integer \(n\geq 1\), let \(\mathbb{C}^n\) denote the \(n\)-dimensional complex Euclidean space and \(\mathbb{B}_n\) its unit ball.  For $z=(z_1,\dots,z_n)$ and $w=(w_1,\dots,w_n)$ in $\mathbb{C}^n$, we define
$
\langle z,w\rangle := \sum_{j=1}^{n} z_j\,\overline{w_j}$
and $
|z| := \sqrt{\langle z,z\rangle}.
$
 For $\alpha>-1$, the weighted Lebesgue measure $dv_\alpha$ is defined by
\[
dv_\alpha(z):=c_\alpha (1-|z|^2)^\alpha\,dv(z),
\]
where $dv$ denotes the normalized Lebesgue measure on
\(\mathbb B_n\),
\begin{equation}\label{vc}
c_\alpha:=\frac{\Gamma(n+1+\alpha)}{n!\,\Gamma(\alpha+1)}
\end{equation}
is chosen so that $v_\alpha(\mathbb{B}_n)=1$, and $\Gamma$ denotes the Gamma function.
\begin{nota}
Denote by $d\lambda$ the M\"{o}bius invariant measure on $\mathbb{B}_n$. Specifically,
$$d\lambda(z):=\frac{dv(z)}{(1-|z|^2)^{n+1}}.$$
\end{nota}

Given $\alpha>-1$ and $0<p<\infty$, for $f \in L^p(dv_\alpha):=L^p(\mathbb{B}_n,dv_\alpha)$,  we write
\[
\|f\|_{p,\alpha}:=\left(\int_{\mathbb{B}_n}|f(z)|^p\,dv_\alpha(z)\right)^{1/p}.
\]
For a domain $\Omega\subset\mathbb C^n$, we denote by $H(\Omega)$ the space of all holomorphic functions on $\Omega$. For $\alpha>-1$ and $0<p<\infty$, we define the weighted Bergman space by
$
A_\alpha^p:=H(\mathbb B_n)\cap L^p(dv_\alpha).
$

Carleson measures were originally introduced by Carleson in his study of
interpolation problems and the corona theorem for the algebra \(H^\infty\) of
bounded analytic functions on the unit disk \cite{MR0117349,MR0141789}.
Since then, they have become an important tool in the analysis of function spaces
and related operator theory. In the setting of Bergman spaces, early contributions on such measures include the works of Oleinik \cite{MR0369705}, Hastings \cite{MR0374886}, Cima and Wogen \cite{MR0650200}, and Luecking \cite{MR0687635,MR0778090}.
\begin{definition}
Let \(\mu\) be a positive Borel measure on \(\mathbb{B}_n\), and let
\(0< p,q<\infty\) and \(\alpha>-1\). Whenever the natural inclusion
$
A_\alpha^p \subset L^q(d\mu)
$
is well-defined, we denote by \(J_\mu\) the associated Carleson embedding operator, that is,
\[
J_\mu:A_\alpha^p \to L^q(d\mu),\qquad J_\mu(f):=f.
\]
\end{definition}

Recently, the study of absolutely summing operators on various spaces of analytic functions has attracted substantial attention (see, e.g., \cite{MR4534902,MR4676254,HW2025,MR3886175,LL2025,MR4690506,MR4905107,MR5040750}). A significant breakthrough in this direction was obtained by Lef\'evre and Rodr\'iguez-Piazza \cite{MR3886175}, who obtained the first complete characterization of absolutely summing Carleson embeddings on the Hardy spaces $H^p(\mathbb{D})$, thereby laying the groundwork for much of the subsequent development. In particular, He-Jreis-Lef\`{e}vre-Lou \cite{MR4690506} proved that if $1<p<\infty$, $1\leq r<\infty$ and $\alpha>-1$, then the Carleson embedding operator
$
J_\mu:A_\alpha^p(\mathbb{D}) \to L^p(\mathbb{D},d\mu)$
is $r$-summing if and only if
$$\hat{\mu}_p(z):=\bigg(\frac{\mu(B(z,\delta))}{v_\alpha(B(z,\delta))}\bigg)^{1/p}\in L^\kappa(\mathbb{D},d\lambda),$$
where $B(z,r)$ denotes the Bergman metric ball centered at $z$ with radius $r$ (see \eqref{defofball}).
The argument in \cite{MR4690506} can also be adapted to the unit ball
\(\mathbb B_n\) in arbitrary dimension for \(1<p<\infty\), while Hu and Wang \cite{HWarxiv} have recently extended the corresponding characterization to the endpoint \(p=1\). To understand
how absolute summability is affected by the interaction between different integrability scales, it is natural to extend the above characterization to the off-diagonal setting and ask the following question:\smallskip

{\bf Question 2:} {\it For $ 1\leq p,q<\infty$, $1\leq r<\infty$ and $\alpha>-1$, can one characterize the measures $\mu$ for which $J_\mu\in\Pi_r(A_\alpha^p(\mathbb B_n),L^q(\mathbb B_n,d\mu))$?}\smallskip

In the course of characterizing the \(r\)-summing property of Toeplitz
operators, Hu and Wang \cite{HWarxiv} obtained, as a byproduct, a solution to Question~2 in the low-exponent range \(1\leq p,q\leq 2\). Despite this valuable progress, Question~2 had remained largely open, particularly in the off-diagonal range $p\neq q$ with \(\max\{p,q\}>2\). Indeed, compared with the low-exponent range, the remaining parameter regimes are substantially more delicate and involve two main challenges. First, in the exceptional parameter range, a complete \(r\)-summing characterization was not known even for the simpler class of diagonal operators, which is closely related to the study of Carleson embedding operators. Second, for a large range of parameters, no suitable block-diagonal principle was available to transfer such a characterization from diagonal operators to Carleson embeddings.

To formulate our second main result, we introduce the following
averaging function.

\begin{definition}\label{def:averaging-function}
Let \(\mu\) be a positive Borel measure on \(\mathbb B_n\), let
\(0< p,q<\infty\), and let \(\alpha>-1\). For \(0<\delta\le 1\),
the \((p,q,\alpha)\)-averaging function of \(\mu\) is defined by
\[
\hat{\mu}_{p,q,\alpha,\delta}(z)
:=
\frac{\mu(B(z,\delta))^{1/q}}
     {v_\alpha(B(z,\delta))^{1/p}},
\qquad z\in\mathbb B_n.
\]
\end{definition}

With this notation and building on the result of Theorem~\ref{GarlingTheorem}, our second main result provides a complete answer to Question~2 (see Definition~\ref{defoflattice} for the definition of a $\delta$-lattice).
\begin{theorem}\label{MainTheorem2}
Let $ 1\leq p,q<\infty$, $1\leq r<\infty$, and $\alpha>-1$. Then the following statements are equivalent:
\begin{enumerate}
  \item\label{carleson1} $J_\mu\in\Pi_r(A^p_\alpha,L^q(d\mu));$
    \item\label{carleson3} $\hat{\mu}_{p,q,\alpha,\delta}\in L^\kappa(\mathbb{B}_n,d\lambda)$ \text{for some (equivalently, every)} $0<\delta\leq 1;$
  \item\label{carleson2} $\{\hat{\mu}_{p,q,\alpha,\delta}(a_k)\}\in \ell^\kappa$ \text{for some (equivalently, every)} $\delta$-lattice $\{a_k\}$ in $\mathbb{B}_n$, $0<\delta\leq 1.$
\end{enumerate}
Moreover, for every $0<\delta\leq 1$ and every $\delta$-lattice $\{a_k\}$ in $\mathbb{B}_n$,
\begin{align*}
\pi_r(J_\mu:A^p_\alpha\rightarrow L^q(d\mu))\simeq \|\hat{\mu}_{p,q,\alpha,\delta}\|_{L^\kappa(\mathbb{B}_n,d\lambda)}\simeq \|\{\hat{\mu}_{p,q,\alpha,\delta}(a_k)\}\|_{\ell^\kappa}.
\end{align*}
The implicit constants are independent of \(\mu\) and of the lattice
points \(a_k\).
\end{theorem}



\subsection{Part III: Hankel operators}
For \(1\le p\le\infty\), let $\overline{A^p_\alpha}:=\{\overline{f} : f \in A^p_\alpha\}$ be the space of conjugate analytic functions in $L^p(dv_\alpha)$.  Let $P_\alpha$ and $\overline{P_\alpha}$ denote the orthogonal projection from $L^2(dv_\alpha)$ onto $A^2_\alpha$ and $\overline{A^2_\alpha}$, respectively. By means of their integral representations (see \eqref{Parepresent} and \cite[p.~707]{MR1093426}), $P_\alpha$ and $\overline{P_\alpha}$ can be extended to operators on \(L^1(dv_\alpha)\).
\begin{definition}
Let $\alpha>-1$. For every $f\in L^1(dv_\alpha)$, the (big) Hankel operator $H^\alpha_f$ and the little Hankel operator $h_f^\alpha$, initially defined on bounded holomorphic functions on $\mathbb{B}_n$, are given by
$$H^\alpha_f(g):=(I-P_\alpha)(fg),$$
$$h_f^\alpha (g) := \overline{P_\alpha}(fg),$$
respectively,
where $I$ denotes the identity operator.
\end{definition}

Hankel operators constitute one of the most important classes of operators on analytic function spaces and play a significant role in several areas of mathematics, notably in functional analysis, complex analysis, operator theory, and certain parts of control theory (for further background and applications, see the monograph of V. Peller \cite{P2}).

A central theme in the theory of Hankel operators has long been to characterize
increasingly refined operator properties in terms of the behavior of their symbols.  For big Hankel operators, Luecking
\cite{MR1194989} showed in 1992 that boundedness and compactness from
\(A^p(\mathbb D)\) to \(L^p(\mathbb D)\), \(1<p<\infty\), as well as
Schatten class membership on \(A^2(\mathbb D)\), can be characterized in
terms of the symbol. For little Hankel operators, Zhu \cite{MR1093426} established corresponding characterizations of boundedness, compactness, and membership in
\(S^{\mathfrak{p}}\), \(1\leq \mathfrak{p}<\infty\), on \(A^2(\mathbb D)\).  See also \cite{PRarxiv} for the
extension of Zhu's boundedness characterization to $A^p(\mathbb{D})$ for
\(1<p<\infty\). By contrast, despite the substantial progress made over the years in these three
directions (see, e.g.,
\cite{MR850538,MR1113389,MR970119,MR4402674,
MR4668087,MR3158507,MR1194989,MR3551773,
MR3010276,MR1860488,MR1951248,MR4850427,MR1013987,
MR1087805,MR4631347,
MR3394382,MR4864874,MR1093426,MR971880,
MR1044786,MR4413302,MR2311536,MR4140739,MR1073289,MR3936102,MR412721,MR2148810,MR5088306,MR1260841}
and the references therein),
a general characterization of \(r\)-summing
Hankel operators beyond the Hilbertian setting remained unavailable.

Motivated by the central role of absolutely summing operators in Banach space geometry and of Hankel operators in operator theory on analytic function spaces, it is natural to ask the following basic question:\smallskip

{\bf Question 3:} {\it For $1\leq p<\infty, 1<q<\infty, 1\leq r<\infty$, and $\alpha,\beta>-1$, can one characterize the symbols $f$ and $g$ for which $H_f^\beta\in\Pi_r(A_\alpha^p,L^q(dv_\beta))$ and $h_{\bar g}^\beta\in \Pi_r(A_\alpha^p,L^q(dv_\beta))$, respectively?}\smallskip

Building on the result of Theorem~\ref{MainTheorem2}, we answer Question~3 separately for big and little Hankel operators. To characterize the $r$-summing property for big Hankel operators, we introduce the weighted IDA spaces as the appropriate characterization spaces, which are defined as follows.
\begin{definition}
For any $\delta>0$, $q\geq 1$, and $f\in L^q_{\rm loc}(\mathbb{B}_n)$, we define the $IDA$ function $G_{q,\delta}(f)$ as
\begin{align*}
G_{q,\delta}(f)(z):=\inf_{h\in H(B(z,\delta))}\left\{\frac{1}{v(B(z,\delta))}\int_{B(z,\delta)}|f(w)-h(w)|^q\,dv(w)\right\}^{1/q}.
\end{align*}
\end{definition}

\begin{definition}\label{weightedIDA}
Let \(1\le q<\infty\) and
$
s\in[1,\infty]\cup\{s_0^{-}:1<s_0<2\}.
$ Let \(t\in\mathbb R\). The weighted IDA spaces $IDA_{t,q}^s(\mathbb{B}_n)$ (Integral Distance to Analytic Functions)  are defined as the sets of $f\in L^q_{\rm loc}(\mathbb{B}_n)$ such that
$$\|f\|_{IDA_{t,q}^s(\mathbb{B}_n)}:=\|(1-|\cdot|^2)^{t}G_{q,\delta}(f)\|_{L^s(\mathbb{B}_n,d\lambda)}<+\infty\quad {\rm for}\ {\rm some}\  0<\delta\leq 1.$$
We shall see in Theorem \ref{mainhankel} that the definition is independent of $0<\delta\leq 1$. When $q=2$ and $t=0$, we briefly write $IDA^s:=IDA_{0,2}^s$.
\end{definition}
The idea of introducing weighted IDA spaces is inspired by the viewpoint of using local distance to analytic functions in the influential works of Luecking {\rm \cite{MR1194989}} and Hu--Virtanen {\rm \cite{MR4668087,MR4402674}}. In particular, Hu--Virtanen {\rm \cite{MR4668087,MR4402674}} introduced the unweighted $ IDA^{\mathfrak{p}}(\mathbb{C})$ space in the setting of Fock space $F^2$ and proved that for any $0<\mathfrak{p}<\infty$,
$$H_f\in S^{\mathfrak{p}}(F^2\rightarrow L^2(\mathbb{C},e^{-|z|^2}dA)) \Longleftrightarrow f\in IDA^{\mathfrak{p}}(\mathbb{C}).$$
This resolved a problem that had remained open for more than twenty years. Using the notion of $IDA$ spaces introduced in {\rm \cite{MR4668087,MR4402674}}, one can
reformulate the classical characterization of Schatten-class Hankel operators
on Bergman spaces (see \cite{MR1194989,MR1332063,MR3770363}) as follows: for every  \(0<\mathfrak{p}<\infty\) and $\alpha>-1$,
\begin{align}\label{SchattenHankel}
H_f^\alpha\in S^{\mathfrak{p}}\bigl(A_\alpha^2(\mathbb{B}_n),L^2(\mathbb{B}_n,dv_\alpha)\bigr)
\quad\Longleftrightarrow\quad
f\in IDA^{\mathfrak{p}}(\mathbb{B}_n).
\end{align}

In the present paper, we introduce weighted and Orlicz extensions of the IDA spaces, and use these extensions to characterize the $r$-summing property of Hankel operators $H^\beta_f:A^p_\alpha\to L^q(dv_\beta)$. Our third main result is the following.
\begin{theorem}\label{Maintheorem}
	Let $1\leq p<\infty, 1<q<\infty, 1\leq r<\infty$, $\alpha,\beta>-1$, and $f\in L^1(dv_\beta)$. Set $$\gamma:=\frac{n+1+\beta}{q}-\frac{n+1+\alpha}{p}.$$  Then
$H^\beta_f\in \Pi_{r}(A^p_\alpha,L^q(dv_\beta))$ if and only if $f\in IDA_{\gamma,q}^{\kappa}(\mathbb{B}_n).$ Moreover,
\begin{align*}
\pi_r(H_f^\beta:A^p_\alpha\rightarrow L^q(dv_\beta))\simeq \|f\|_{IDA_{\gamma,q}^{\kappa}(\mathbb{B}_n)}.
\end{align*}
\end{theorem}
\begin{remark}\label{keyremarkHankel}
Building on Theorem~\ref{Maintheorem}, we carry out a further analysis
to obtain the following results:
\begin{enumerate}
  \item In Theorem~\ref{mainhankel}, we provide an equivalent formulation
  of Theorem~\ref{Maintheorem}, together with a discrete characterization
  in terms of \(\delta\)-lattices;
  \item In Corollary~\ref{mainhankelb}, we characterize the simultaneous
  \(r\)-summability of the associated Hankel operators
  \[
  H^\beta_f:A^p_\alpha\to L^q(dv_\beta)
  \quad\text{and}\quad
  H^\beta_{\bar f}:A^p_\alpha\to L^q(dv_\beta)
  \]
  in terms of membership in the weighted $\mathrm{IMO}$ space;
  \item In Corollary~\ref{mainhankelc}, we characterize the \(r\)-summability of the Hankel operator with
  anti-analytic symbol \(\bar f\),
  \[
  H^\beta_{\bar f}:A^p_\alpha\to L^q(dv_\beta),
  \]
  in terms of membership in the weighted Besov space.

  Furthermore, in Corollary \ref{mainhankelcconstant}, we determine the
sharp parameter ranges for which
$
H^\beta_{\bar f}:A^p_\alpha\to L^q(dv_\beta)
$
is \(r\)-summing if and only if \(f\) is constant.
\end{enumerate}
\end{remark}

To characterize the $r$-summing property of little Hankel operators, we define the integral operator $V_{b,\beta}$ as follows.

\begin{definition}\label{Defofintegralop}
For $b>0$ and $\beta>-1$, we let $V_{b,\beta}$  be the integral operator defined by
\begin{align}\label{integraloperatorkey}
V_{b,\beta} (f)(z) :=\frac{c_{b+\beta}(1 - |z|^2)^{b} }{c_\beta}\int_{\mathbb{B}_n} \frac{f(w)}{(1 - \langle z,w\rangle )^{n+1+\beta+b}} \, dv_\beta(w),\quad f\in L^1(dv_\beta),
\end{align}
where $c_{b+\beta}$ and $c_\beta$ are the volume constants defined in \eqref{vc}.
\end{definition}

The idea of using integral transforms like \eqref{integraloperatorkey} to characterize little Hankel operators can be traced back to \cite{MR1093426} (see also the earlier work \cite{MR1008445} for a related use of this idea in the study of Hankel forms). Using this integral transform, one can formulate the classical characterization of Schatten-class little Hankel operators
on Bergman spaces (see \cite{MR1093426,MR1333834}) as follows: for every  \(1\leq \mathfrak{p}<\infty\), $\alpha>-1$,  and \(f\in L^2(\mathbb B_n,dv_\alpha)\),
\begin{align}\label{SchattenliHankel}
h_{\bar{f}}^\alpha\in S^{\mathfrak{p}}\bigl(A_\alpha^2(\mathbb{B}_n),L^2(\mathbb{B}_n,dv_\alpha)\bigr)
\quad\Longleftrightarrow\quad
 V_{n+1+\alpha,\alpha}(f)\in L^{\mathfrak{p}}(\mathbb{B}_n,d\lambda).
\end{align}


In the present paper, we introduce weighted and Orlicz extensions of the function space $\{f\in L^2(\mathbb B_n,dv_\alpha):V_{n+1+\alpha,\alpha}(f)\in L^\mathfrak{p}(\mathbb{B}_n,d\lambda)\}$, under a slightly weaker a priori integrability assumption, and use these extensions to characterize the $r$-summing property of little Hankel operators $h^\beta_{\bar{f}}:A^p_\alpha\to L^q(dv_\beta)$. For \(0<\delta\leq1\), let \(\mathfrak L_\delta\) be the collection of all
\(\delta\)-lattices in \(\mathbb B_n\). Our fourth main result is the following.
\begin{theorem}\label{Maintheoremlittle}
	Let $1\leq p<\infty, 1<q<\infty, 1\leq r<\infty$, and $\alpha,\beta>-1$. Let $f\in L^1(dv_\beta)$, and let $$\gamma=\frac{n+1+\beta}{q}-\frac{n+1+\alpha}{p}.$$ Then the following statements are equivalent:

\begin{enumerate}
  \item\label{littlehankel1} $h^\beta_{\bar{f}}\in\Pi_r(A^p_\alpha,L^q(dv_\beta));$
  \item\label{littlehankel3} $(1-|\cdot|^2)^\gamma |V_{b,\beta}(f)|\in L^\kappa(\mathbb{B}_n,d\lambda) $ \text{for some (equivalently, every)} $b>n+(\alpha+1)/p;$
    \item\label{littlehankel2} $\{(1-|a_k|^2)^\gamma |V_{b,\beta}(f)(a_k)|\}\in \ell^\kappa$ for some (equivalently, every)
\(b>n+(\alpha+1)/p\), some (equivalently, every)
\(0<\delta\leq1\) and for every \(\delta\)-lattice \(\{a_k\}\) in \(\mathbb B_n.\)
\end{enumerate}
Moreover, for every $b>n+(\alpha+1)/p$  and every
\(0<\delta\leq1\),
\begin{align*}
\pi_r(h^\beta_{\bar{f}}:A^p_\alpha\rightarrow L^q(dv_\beta))\simeq \|(1-|\cdot|^2)^\gamma |V_{b,\beta}(f)|\|_{L^\kappa(\mathbb{B}_n,d\lambda)}\simeq \sup_{\{a_k\}\in\mathfrak L_\delta}\|\{(1-|a_k|^2)^\gamma |V_{b,\beta}(f)(a_k)|\}\|_{\ell^\kappa}.
\end{align*}
\end{theorem}

%

\begin{remark}
Recall from {\rm \cite[Chapter~4]{MR1342297}} that for any Hilbert spaces \(H_1\), \(H_2\) and every
\(1\leq r<\infty\),
\[
\Pi_r(H_1,H_2)=S^2(H_1,H_2),
\]
with equivalence of the corresponding norms. Therefore, when \(p=q=2\) and
\(\alpha=\beta\), Theorems~\ref{Maintheorem} and \ref{Maintheoremlittle} reduce to the \(\mathfrak{p}=2\) case of the
Schatten-class characterizations in \eqref{SchattenHankel} and \eqref{SchattenliHankel}, respectively. Hence,
Theorems~\ref{Maintheorem} and \ref{Maintheoremlittle} may be regarded as Banach-space extensions
of the Hilbert--Schmidt class characterizations of big Hankel operators and little Hankel operators, respectively.
\end{remark}

\begin{remark}
The present paper incorporates and supersedes {\rm arXiv: 2511.22165v1},
with its title changed to the present one and substantially new arguments added. In particular, the case \(p=q\) of Theorems \ref{Maintheorem} and
\ref{Maintheoremlittle} was already contained in that earlier version.
To the best of our knowledge, the conclusions of Theorems \ref{Maintheorem} and \ref{Maintheoremlittle}
are both new even in the case \(p=q\) and $\alpha=\beta$. In particular, apart from the
classical theory of Schatten classes on Hilbert spaces, namely the case $p=q=2$, the present work, together with the independent work {\rm \cite{HLarxiv}} on Fock spaces, appears to be among
the first systematic studies of \(r\)-summing Hankel operators
on analytic function spaces.
\end{remark}
\subsection{Difficulties and Strategies}\label{technicalpoints}
The main difficulty in Part I is to identify the correct
characterization space \(\ell^\kappa\) in the exceptional range
\((p,q,r)\in E\). Our strategy is as follows:
\begin{itemize}
  \item For the sufficiency part, we first lift the diagonal operator to the
bilinear map
\[
\mathcal B\bigl(\mathbf b,\{x^{(m)}\}_m\bigr)
:=
\bigl\{\mathscr M_{\mathbf b}x^{(m)}\bigr\}_m
\]
and establish its boundedness at the two endpoint pairs
\[
\mathcal B:
\ell^r\times\ell^r_w(\ell^p)
\longrightarrow
\ell^r(\ell^r)
\qquad\text{and}\qquad
\mathcal B:
\ell^{p'}\times\ell^r_w(\ell^p)
\longrightarrow
\ell^r(\ell^2).
\]
See Definition \ref{def:weak-strong-r-summable} for the notation of these mixed-norm spaces.
We then apply the bilinear complex interpolation theorem with parameter
\(\vartheta\) to show that
\[
\mathcal B:
\ell^\kappa\times\ell^r_w(\ell^p)
\longrightarrow
\ell^r(\ell^q)
\]
is bounded, which gives the desired \(r\)-summing estimate.

  \item For the necessity part, we combine a randomization argument
  with Rosenthal's inequality to establish a self-improvement phenomenon (see Lemma \ref{lem:GarlingNec-t-q-estimate}): the \(r\)-summability of
  \(\mathscr M_{\mathbf b}:\ell^p\to\ell^q\) forces
  \(\mathscr M_{\mathbf b}\) to admit a bounded extension from the
  strictly larger space \(\ell^t\) to \(\ell^q\), where \(t\) is the
sharp exponent for this domain enlargement. The boundedness characterization for diagonal operators then identifies the required characterization space as \(\ell^\kappa\).
\end{itemize}

Our main strategies of Parts II and III lie in the development of two new general principles for $r$-summability, established in Propositions~\ref{lem:block-diagonal-transfer2} and~\ref{proq2pr}. Their respective roles are as follows:
\begin{itemize}
  \item For the sufficiency part, Proposition~\ref{lem:block-diagonal-transfer2}
transfers the quantitative \(r\)-summing norm estimate for an associated
diagonal operator, whose diagonal entries are the \(1\)-summing norms of
the individual blocks, to the corresponding estimate for the
block-diagonal operator. The underlying block-diagonal argument is inspired by the
 works \cite{MR3886175,MR3078637,MR3538650,MR3742997,MR4745822,MR4690506}. In the aforementioned literature, the available block-diagonal principles take different
forms according to the relative values of \((p,q,r)\), while no
corresponding principle is available in those works for a large range of parameters. In contrast,
Proposition~\ref{lem:block-diagonal-transfer2} presents a unified block-diagonal principle that holds for all \(1\leq p\leq\infty\), \(1\leq q<\infty\), and \(1\leq r<\infty\), with a rather simple proof.

  \item For the necessity part, inspired by the ingenious randomization argument developed in \cite{MR3886175}, Proposition \ref{proq2pr} transfers the quantitative
\(r\)-summing norm estimate for an arbitrary operator
\(T:A_\alpha^p\to L^q(d\mu)\) to the corresponding estimate for an associated diagonal operator determined by certain localized testing quantities. This principle applies in particular to, but is not restricted to, Hankel operators, and is therefore also of independent interest for operators on Bergman spaces.
\end{itemize}

An additional technical difficulty, which persists throughout all three parts
of the paper, arises in the borderline regime
\(2<p<\infty\), \(1<q=p'<2\), and \(1\le r<q\).
In this regime, the natural summability scale is a logarithmically corrected
Orlicz scale rather than an ordinary Lebesgue scale, and the corresponding
analysis requires additional care in handling the Luxemburg norm.
\subsection{Organization of the paper}
This  paper is organized as follows. In Section \ref{Sect2}, we collect preliminary properties of Bergman metric balls, $\delta$-lattices and absolutely summing operators, and introduce the normalized modified test functions, which serve as substitutes for normalized reproducing kernels in the endpoint case $p=1$. In Section \ref{ASDOsect}, we introduce the endpoint Orlicz spaces in order to study the borderline case, record their
basic properties needed later, and provide a proof of Theorem \ref{GarlingTheorem}. In Section \ref{Sect3}, we establish two general principles for $r$-summability. In Section \ref{Sect4}, we characterize those measures $\mu$ for which the Carleson embedding operator $J_\mu:A_\alpha^p \to L^q(d\mu)$ is $r$-summing. Sections \ref{Sect5} and \ref{Sect6} are devoted to the characterization of absolutely summing big and little Hankel operators, respectively.

\section{Preliminaries}\label{Sect2}
\setcounter{equation}{0}
\subsection{Notation}
Throughout the paper, \(C\) denotes a positive constant, which may change
from line to line and whose dependence is restricted to the fixed
structural parameters under consideration. We write \(A\lesssim B\) if
\(A\le CB\) for some such constant \(C>0\), \(A\gtrsim B\) if
\(B\lesssim A\), and \(A\simeq B\) if both \(A\lesssim B\) and
\(A\gtrsim B\).

We let $\mathbb{N}$ denote the set of positive integers. For $j\in \mathbb{N}$, let \(e_j\) be the \(j\)-th canonical unit vector. The indicator function of a subset $F\subset \mathbb{B}_n$ is denoted by $\chi_{F}$.  When the domain and target space are clear from the context, \(\|T\|\) denotes the corresponding operator norm of $T$. $\mathcal L(X,Y)$ denotes the Banach space of all bounded linear operators from $X$ to $Y$.

For \(1\le p\le \infty\), we denote by \(p'\) the conjugate exponent
determined by
$
1/p+1/p'=1,
$
with the convention that \(1/\infty=0\).
\subsection{Standard estimates}
We let $d$ be the Bergman metric on $\mathbb{B}_n$ (see \cite{MR2115155} for the definition). For $z\in \mathbb{B}_n$ and $r>0$, we denote by $B(z,r)$ the Bergman metric ball centered at $z$ with radius $r$, that is
\begin{align}\label{defofball}
B(z,r):=\{w\in \mathbb{B}_n: d(z,w)<r\}.
\end{align}
It is well known that (see, e.g., \cite[Lemmas 1.24 and 2.20]{MR2115155}) for any fixed $\delta>0$,
\begin{equation}\label{a2}
	v_\alpha(B(z,\delta))\simeq (1-|z|^2)^{n+1+\alpha},
\end{equation}
and for any $z\in \mathbb{B}_n$ and $w\in B(z,\delta)$, we have
\begin{equation}\label{a3}
	|1-\langle z,w\rangle|\simeq 1-|w|^2\simeq 1-|z|^2,
\end{equation}
where the implicit constants are independent of $z$ and $w$.  By  \eqref{a2} and \eqref{a3},  we have the following local volume comparability: for $w\in \mathbb{B}_n$ and $z\in B(w,\delta)$,
\begin{equation}\label{b1}
v_\alpha(B(z,\delta))\simeq (1-|z|^2)^{n+1+\alpha}\simeq (1-|w|^2)^{n+1+\alpha}\simeq v_\alpha\bigl(B(w,2\delta)\bigr).
\end{equation}

\begin{definition}\label{defoflattice}
	Given $\delta>0$, we say that a sequence $\{a_j\}$ in $\mathbb{B}_n$ is a $\delta$-lattice in the Bergman metric if the following conditions are satisfied:
	\begin{enumerate}
		\item  $\{B(a_j,\delta)\}$ is a covering of $\mathbb{B}_n$;
		\item  $B(a_i,\delta/4)\cap B(a_j,\delta/4)=\emptyset$ whenever $i\neq j$.
	\end{enumerate}
\end{definition}
By \cite[Theorem 2.23]{MR2115155},  $\delta$-lattices exist for any $0<\delta\leq 1$.
A sequence \(\{a_j\}\) in \(\mathbb B_n\) is called
\(\sigma\)-separated, where \(\sigma>0\), if
$
d(a_i,a_j)\ge \sigma
$
whenever \(i\ne j\). It is called separated if it is
\(\sigma\)-separated for some \(\sigma>0\). It is well known that if $\{a_j\}$ is a $\delta$-lattice, then we have the following finite multiplicity property: for each fixed constant $c>0$,
\begin{equation}\label{a1}
	\sum_j \chi_{B(a_j,c\delta)}(z)\leq N,\quad z\in \mathbb{B}_n,
\end{equation}
where $N$ depends on $c$ and $\delta$ but is independent of $z$. More generally, the same argument as in the proof of \cite[Theorem~2.23(3)]{MR2115155} shows that \eqref{a1} holds for every $\delta$-separated sequence, with $N$ depending only
on $c$ and $\delta$.

We have the following local mean value inequality on Bergman balls.
\begin{lemma}\label{Bergmaninequality}
	Let $\alpha>-1$, $0<t,\delta<\infty$. Then there exists a constant
	$C=C(n,\alpha,\delta)>0$ such that
	\[
	|f(z)|^t
	\le
	\frac{C}{v_\alpha(B(z,\delta))}
	\int_{B(z,\delta)} |f(w)|^t\,dv_\alpha(w)
	\]
	for every $z\in\mathbb B_n$ and every holomorphic or conjugate-holomorphic
	function  $f$ on  $B(z,\delta)$.
\end{lemma}
\begin{proof}
Lemma \ref{Bergmaninequality} is a local version of \cite[Lemma 2.24]{MR2115155}. Its proof is obtained by the same argument, with only the obvious restriction to the Bergman ball under consideration. Therefore, we omit the details.
\end{proof}


\subsection{Test functions}\label{testfunction}
In this subsection, we introduce a modified test function that serves as a substitute for the reproducing kernel and is useful to study the $r$-summing operators in the endpoint case $p=1$.

It is well known that $A^2_\alpha$ is a reproducing kernel Hilbert space. The associated reproducing (Bergman) kernel of $A^2_\alpha$ is given explicitly by
\[
K^\alpha_w(z):=K^\alpha(z,w)=\frac{1}{(1-\langle z,w \rangle)^{n+1+\alpha}}.
\]
Then $P_\alpha$ is an integral operator with kernel $K^\alpha(z,w)$. That is,
\begin{align}\label{Parepresent}
P_\alpha(f)(z)=\int_{\mathbb{B}_n} f(w)\,K^\alpha(z,w)\,dv_\alpha(w),\quad f\in L^2(dv_\alpha).
\end{align}
For each fixed $z\in\mathbb{B}_n$, the function $w\mapsto K^\alpha(z,w)$ is bounded on $\mathbb{B}_n$, so the integral above also makes sense for every $f\in L^1(dv_\alpha)$ and $z\in\mathbb{B}_n$.


For $1\leq p<\infty$ and $b>n+(\alpha+1)/p$, we define the modified reproducing kernel $H_{p,z}^b$ by
\begin{align}
H_{p,z}^b(w):=\frac{1}{(1-\langle w,z\rangle)^{b}},\quad z,w\in\mathbb{B}_n.
\end{align}
For any $1<p<\infty$,  we have $H_{p,z}^{n+1+\alpha}=K^\alpha_z$. Thus, $H_{p,z}^b$ may be viewed as a natural modification of the reproducing kernel. Moreover, it follows from \cite[Theorem 1.12]{MR2115155} that for any $1\leq p<\infty$,
\begin{equation}\label{a4s}
	\|H_{p,z}^b\|_{p,\alpha}\simeq (1-|z|^2)^{(-pb+n+1+\alpha)/p},\quad z\in\mathbb{B}_n.
\end{equation}
Next, for any $1\leq p<\infty$, we define the $L^p(dv_\alpha)$-normalized function of $H_{p,z}^b$ by
\begin{align*}
h_{p,z}^b(w):=\frac{H_{p,z}^b(w)}{\|H_{p,z}^b\|_{p,\alpha}},\quad z,w\in\mathbb{B}_n.
\end{align*}
By \eqref{a4s} and \eqref{a3}, for any fixed $\delta>0$,
\begin{equation}\label{a5s}
	|h_{p,z}^b(w)|\simeq\frac{(1-|z|^2)^{(pb-n-1-\alpha)/p}}{|1-\langle z,w \rangle|^{b} }\simeq (1-|z|^2)^{-(n+1+\alpha)/p},\quad w\in B(z,\delta).
\end{equation}

\subsection{Absolutely summing operators}
In this subsection, we collect several general properties about $r$-summing operators from \cite{MR1342297}.
\begin{definition}\label{def:weak-strong-r-summable}
Let \(X\) be a Banach space and let \(1\leq r<\infty\). The space of
strongly \(r\)-summable sequences in \(X\), denoted by \(\ell^r(X)\),
consists of all \(X\)-valued sequences
\(\mathbf x=\{x^{(m)}\}_{m\geq1}\) such that
\[
\|\mathbf x\|_{\ell^r(X)}
:=
\Bigg(
\sum_{m=1}^{\infty}
\|x^{(m)}\|_X^r
\Bigg)^{1/r}
<\infty.
\]
The space of weakly \(r\)-summable sequences in \(X\), denoted by
\(\ell^r_w(X)\), consists of all \(X\)-valued sequences
\(\mathbf x=\{x^{(m)}\}_{m\geq1}\) such that
\[
\|\mathbf x\|_{\ell^r_w(X)}
:=
\sup_{\varphi\in B_{X^*}}
\Bigg(
\sum_{m=1}^{\infty}
|\varphi(x^{(m)})|^r
\Bigg)^{1/r}
<\infty.
\]
\end{definition}
Let $X$ and $Y$ be two Banach spaces, and let $1\le r<\infty$. By \cite[Proposition~2.1]{MR1342297}, \(T:X\to Y\) is \(r\)-summing if and only if it induces a bounded operator
\[
\widetilde T:\ell^r_w(X)\longrightarrow \ell^r(Y),\qquad
\widetilde T\bigl(\{x^{(m)}\}_{m\geq1}\bigr)
:=
\{Tx^{(m)}\}_{m\geq1}.
\]
Moreover, the space $\Pi_{r}(X,Y)$, equipped with the norm $\pi_{r}$, is a Banach space. The space $\Pi_{r}(X,Y)$ is a linear subspace of $\mathcal{L}(X,Y)$, and every $T\in \Pi_{r}(X,Y)$ satisfies
\[
\|T\|\le \pi_{r}(T).
\]

The following ideal property plays a crucial role in our analysis
(see \cite[p.~37]{MR1342297}). Let \(1\leq r<\infty\), and let
\(X_0,X,Y\), and \(Y_0\) be Banach spaces. If
\(T\in\Pi_r(X,Y)\), \(S\in\mathcal L(X_0,X)\), and
\(U\in\mathcal L(Y,Y_0)\), then \(UTS\in\Pi_r(X_0,Y_0)\) and
\begin{equation}\label{ideal-a}
\pi_r(UTS)\le \|U\|\,\pi_r(T)\,\|S\|.
\end{equation}
%


%

\section{Absolutely summing diagonal operators}\label{ASDOsect}
\setcounter{equation}{0}

\subsection{Orlicz spaces}\label{Orlicz}
In this subsection, we collect several useful properties of Orlicz spaces, which will be used to characterize the $r$-summing property of operators in the borderline case.
\begin{definition}{\rm \cite[Definition 3.5.2]{MR3243741}}
A Young function is a continuous increasing convex function $\Phi$ on $[0,\infty)$ that satisfies $\Phi(0)=0$ and $\lim_{t\rightarrow \infty}\Phi(t)=\infty$.
\end{definition}

\begin{definition}
Let $\Phi$ be a Young function. We define $L_\Phi(\mathbb{B}_n,d\lambda)$ to be the Orlicz space of all \(d\lambda\)-measurable functions \(f:\mathbb{B}_n\to\mathbb C\) such that \[ \int_{\mathbb{B}_n} \Phi\left(\frac{|f(z)|}{\rho}\right)\,d\lambda(z)<\infty \] for some \(\rho>0\).
It is equipped with the Luxemburg norm
\[ \|f\|_{L_\Phi(\mathbb{B}_n,d\lambda)} := \inf\left\{ \rho>0: \int_{\mathbb{B}_n} \Phi\left(\frac{|f(z)|}{\rho}\right)\,d\lambda(z)\le 1 \right\}. \]
\end{definition}

\begin{definition}
Let $\Phi$ be a Young function.  We define \(\ell_{\Phi}\) to be the Orlicz sequence space of all complex-valued sequences \(a=\{a_k\}_{k=1}^{\infty}\) such that \[ \sum_{k=1}^{\infty} \Phi\left(\frac{|a_k|}{\rho}\right)<\infty \] for some \(\rho>0\). It is equipped with the Luxemburg norm \[ \|a\|_{\ell_\Phi} := \inf\left\{ \rho>0: \sum_{k=1}^{\infty} \Phi\left(\frac{|a_k|}{\rho}\right)\le 1 \right\}. \]
\end{definition}
It follows directly from the definition that if
\(\Phi(t)=t^p\) for some \(p\geq 1\), then the corresponding Luxemburg
norm coincides with the usual \(L^p\)-norm. Thus,
$
L_\Phi(\mathbb{B}_n,d\lambda)
=
L^p(\mathbb{B}_n,d\lambda)
$
isometrically. However, in the borderline regime $\{2<p<\infty, 1<q=p'<2, 1\le r<q\}$, the natural summability
scale is a logarithmically corrected Orlicz scale rather than the usual
\(L^q\) and \(\ell^q\) scales. Therefore, we introduce the Orlicz index \(q^{-}\) as follows. For \(1<q<\infty\) and \(t>0\), set
\[
M_q(t):=t^q\bigl(1+\log^+(t^{-1})\bigr),
\]
and set \(M_q(0):=0\). Choose \(\varepsilon\in(0,1)\) sufficiently small
so that \(M_q\) is increasing and convex on \([0,\varepsilon]\), and define
\[
\Psi_q(t):=
\begin{cases}
M_q(t), & 0\le t\le \varepsilon,\\[1mm]
M_q(\varepsilon)
+M_q'(\varepsilon)(t-\varepsilon)
+(t-\varepsilon)^2,
& t>\varepsilon.
\end{cases}
\]
It is easy to verify that \(\Psi_q\) is a Young function. Accordingly, we define
\[
L^{q^-}(\mathbb{B}_n,d\lambda)
:=L_{\Psi_q}(\mathbb{B}_n,d\lambda),
\qquad
\ell^{q^-}:=\ell_{\Psi_q},
\]
equipped with the corresponding Luxemburg norms.

%
%
%

\begin{remark}\label{orliczremark}
By {\rm \cite[Propositions~4.a.4 and 4.a.5]{MR500056}}, for any \(1< q<\infty\), $a=\{a_k\}\in \ell^{q^-}$ if and only if $\sum_{k=1}^{\infty}|a_k|^q\bigl(1+\log^+|a_k|^{-1}\bigr)<\infty$, where the summand is understood to be zero whenever \(a_k=0\). This equivalence is only qualitative. It should \textbf{not} be interpreted as a norm equivalence of the form
\begin{align}\label{falsetwoside}
\|a\|_{\ell^{q^-}}\simeq \Bigg(\sum_{k=1}^{\infty}
|a_k|^q
\bigl(1+\log^+|a_k|^{-1}\bigr)\Bigg)^{1/q}.
\end{align}
This can be explained by taking \(a=t e_1\), where \(0<t<1\) and $e_1:=(1,0,0,\ldots)$. By the homogeneity of the Luxemburg norm,
$
\|t e_1\|_{\ell^{q^-}}=t\|e_1\|_{\ell^{q^-}},
$
whereas
\[
\Bigg(\sum_{k=1}^{\infty}
|(te_1)_k|^q
\bigl(1+\log^+|(te_1)_k|^{-1}\bigr)\Bigg)^{1/q}
=
t\bigl(1+\log(1/t)\bigr)^{1/q}.
\]
Thus the ratio between these two quantities is unbounded as \(t\to0^+\). Hence the two-sided estimate \eqref{falsetwoside} is false.
\end{remark}

The following elementary estimate shows that the coordinate functionals on \(\ell^{q^-}\) are continuous.
\begin{lemma}\label{coordinateconver}
Let $1< q<\infty$. For every \(a=\{a_k\}\in\ell^{q^-}\) and every
\(j\in\mathbb N\), we have
\[
|a_j|\le \sup\{t\ge 0:\Psi_q(t)\le 1\}\|a\|_{\ell^{q^-}}.
\]
\end{lemma}
\begin{proof}
Fix \(\lambda>\|a\|_{\ell^{q^-}}\). By the definition of $\|a\|_{\ell^{q^-}}$, we have
\[
\Psi_q\left(\frac{|a_j|}{\lambda}\right)\le
\sum_{k=1}^{\infty}\Psi_q\left(\frac{|a_k|}{\lambda}\right)
\le 1.
\]
Hence,
\[
|a_j|\le \sup\{t\ge 0:\Psi_q(t)\le 1\}\lambda.
\]
Letting \(\lambda\downarrow \|a\|_{\ell^{q^-}}\), we complete the proof of Lemma \ref{coordinateconver}.
\end{proof}
\subsection{Proof of Theorem \ref{GarlingTheorem}}
We begin with recalling a useful inequality from Rosenthal.
\begin{lemma}{\rm \cite[Theorem 3]{MR271721}}\label{RosenthalLemma}
Let $2< r<\infty$. There exists $R_r>0$ such that for every finite family
\(\{\xi_j\}_{j=1}^N\) of independent mean-zero scalar random
variables satisfying $\mathbb{E}|\xi_j|^r<\infty$, $j=1,2,\ldots,N$,
\begin{equation}
\label{scalar-Rosenthal}
\Bigg(\mathbb{E}\bigg|\sum_{j=1}^N \xi_j\bigg|^r\Bigg)^{1/r}
\le
R_r
\Bigg(
\bigg(
\sum_{j=1}^N\mathbb E|\xi_j|^2
\bigg)^{1/2}
+
\bigg(
\sum_{j=1}^N\mathbb E|\xi_j|^r
\bigg)^{1/r}
\Bigg).
\end{equation}
\end{lemma}
The following construction is standard in probability theory. Nevertheless, for completeness, we include a self-contained proof.
\begin{lemma}
\label{lem:Bernoulli-Rademacher-product}
Let \(F\) be a finite set and let
$
\{\rho_j\}_{j\in F}\subset[0,1].
$
Then there exists a finite probability space supporting mutually
independent random variables
$
\{\varsigma_j,\varepsilon_j\}_{j\in F}
$
such that for every \(j\in F\),
\[
\mathbb P(\varsigma_j=1)=\rho_j,
\quad
\mathbb P(\varsigma_j=0)=1-\rho_j,\quad {\rm and}\quad
\mathbb P(\varepsilon_j=1)
=
\mathbb P(\varepsilon_j=-1)
=
\frac12.
\]
\end{lemma}

\begin{proof}
Set
$
\Omega_F:=\{0,1\}^{F}\times\{-1,1\}^{F}.
$
For \(a\in\mathbb R\), let \(\delta_a\) denote the Dirac probability
measure concentrated at \(a\). Equip \(\Omega_F\) with the product
probability measure
\[ \mathbb P := \Bigg( \bigotimes_{j\in F} \bigl[ (1-\rho_j)\delta_0 + \rho_j\delta_1 \bigr] \Bigg) \otimes \Bigg( \bigotimes_{j\in F} \bigl[ \tfrac12\delta_{-1} + \tfrac12\delta_1 \bigr] \Bigg). \]
For
$
\omega=(\{\omega_j\}_{j\in F},\{\eta_j\}_{j\in F})\in\Omega_F,
$
define
$
\varsigma_j(\omega):=\omega_j,$ and $
\varepsilon_j(\omega):=\eta_j.
$
The required properties follow immediately from
the product construction. This finishes the proof of Lemma \ref{lem:Bernoulli-Rademacher-product}.
\end{proof}

\begin{lemma}
\label{lem:Bernoulli-Rademacher-estimate}
Let \(1<p<2\), \(1\le q<\infty\), and
\(\max\{p',q\}<r<\infty\). Define \(s\) by
\[
\frac1s=\frac12-\frac1{p'}.
\]
Let \(\mathbf b=\{b_j\}_{j\ge1}\) be a scalar sequence such that
$
\mathscr M_{\mathbf b}\in\Pi_r(\ell^p,\ell^q).
$
Then for every finite set \(F\subset\mathbb N\), every scalar family
\(\{a_j\}_{j\in F}\), and every family
\(\{\rho_j\}_{j\in F}\subset(0,1]\),
\begin{align}\label{eq:Bernoulli-Rademacher-estimate}
\Bigg(
\sum_{j\in F}\rho_j|b_ja_j|^q
\Bigg)^{1/q}
\le
R_r\pi_r(\mathscr M_{\mathbf b})
\Bigg(
\|\{\rho_j^{1/2}a_j\}_{j\in F}\|_{\ell^s}
+
\|\{\rho_j^{1/r}a_j\}_{j\in F}\|_{\ell^\infty}
\Bigg).
\end{align}
\end{lemma}
\begin{proof}
Let
$
\{\varsigma_j,\varepsilon_j\}_{j\in F}
$ be mutually
independent random variables
provided in Lemma~\ref{lem:Bernoulli-Rademacher-product}.
Define the finitely supported random vector
$
Z
:=
\sum_{j\in F}
a_j\varsigma_j\varepsilon_je_j.
$
By \cite[(1.1)]{MR3886175} (see also \cite[p.~56]{MR1342297}),
\begin{align}
\Bigg(
\mathbb E
\|\mathscr M_{\mathbf b}Z\|_{\ell^q}^r
\Bigg)^{1/r}
\le
\pi_r(\mathscr M_{\mathbf b})
\sup_{\|y\|_{\ell^{p'}}\le1}
\Bigg(
\mathbb E
\bigg|
\sum_{j\in F}
a_jy_j\varsigma_j\varepsilon_j
\bigg|^r
\Bigg)^{1/r}.
\label{random-summing-inequality}
\end{align}

Since \(r/q>1\), we apply Jensen's inequality to see that
\begin{align}
\Bigg(
\mathbb E
\|\mathscr M_{\mathbf b}Z\|_{\ell^q}^r
\Bigg)^{q/r}
=
\Bigg(
\mathbb E (\sum_{j\in F}|b_ja_j|^q\varsigma_j)^{r/q}
\Bigg)^{q/r}
\ge
\mathbb E\Big(\sum_{j\in F}|b_ja_j|^q\varsigma_j\Big)
=
\sum_{j\in F}\rho_j|b_ja_j|^q.
\label{random-left-lower-root}
\end{align}

For fixed \(y\in\ell^{p'}\), the scalar random variables
$
\xi_j
:=
a_jy_j\varsigma_j\varepsilon_j, j\in F,
$
are independent and have mean zero. By Lemma \ref{RosenthalLemma},
\begin{align}
\Bigg(
\mathbb E
\bigg|
\sum_{j\in F}
a_jy_j\varsigma_j\varepsilon_j
\bigg|^r
\Bigg)^{1/r}
\le
R_r
\Bigg(
\bigg(
\sum_{j\in F}
\rho_j|a_jy_j|^2
\bigg)^{1/2}
+
\bigg(
\sum_{j\in F}
\rho_j|a_jy_j|^r
\bigg)^{1/r}
\Bigg).
\label{Rosenthal-applied}
\end{align}

Using H\"older's inequality, we obtain
\begin{equation}
\label{Rosenthal-first-term}
\Bigg(
\sum_{j\in F}
\rho_j|a_jy_j|^2
\Bigg)^{1/2}
\le
\|\{\rho_j^{1/2}a_j\}_{j\in F}\|_{\ell^s}
\|y\|_{\ell^{p'}}.
\end{equation}
Since \(p'<r\), we also have
\begin{align}
\Bigg(
\sum_{j\in F}
\rho_j|a_jy_j|^r
\Bigg)^{1/r}
&\le
\|\{\rho_j^{1/r}a_j\}_{j\in F}\|_{\ell^\infty}
\|y\|_{\ell^r}
\le
\|\{\rho_j^{1/r}a_j\}_{j\in F}\|_{\ell^\infty}
\|y\|_{\ell^{p'}}.
\label{Rosenthal-second-term}
\end{align}
Combining \eqref{random-summing-inequality}--\eqref{Rosenthal-second-term}, we deduce \eqref{eq:Bernoulli-Rademacher-estimate}.
This finishes the proof of Lemma \ref{lem:Bernoulli-Rademacher-estimate}.
\end{proof}

\begin{lemma}
\label{lem:GarlingNec-t-q-estimate}
Let
 $
1<p<2<q<\infty,$ and $
r>\max\{p',q\}.
$
Define \(s\), \(\vartheta\), and \(t\) by
\[
\frac1s=\frac12-\frac1{p'},
\qquad
\vartheta
:=
\frac{\frac1q-\frac1r}
     {\frac12-\frac1r},
\qquad
t:=\frac{s}{\vartheta},
\]
respectively.
Let \(\mathbf b=\{b_j\}_{j\geq1}\) be a scalar sequence such that
$
\mathscr M_{\mathbf b}\in\Pi_r(\ell^p,\ell^q).
$
Then \(\mathscr M_{\mathbf b}\) extends uniquely to a bounded operator
from \(\ell^t\) to \(\ell^q\). Moreover,
\[
\|\mathscr M_{\mathbf b}\|_{\mathcal L(\ell^t,\ell^q)}
\leq
2R_r\,
\pi_r\left(
\mathscr M_{\mathbf b}:\ell^p\to\ell^q
\right).
\]
\end{lemma}
\begin{proof}
By a standard density argument, it suffices to show that for every finitely supported scalar sequence \(w\neq0\), we have
\begin{equation}
\label{eq:GarlingNec-t-q-bound}
\|\mathscr M_{\mathbf b}w\|_{\ell^q}
\leq
2R_r\pi_r(\mathscr M_{\mathbf b})
\|w\|_{\ell^t}.
\end{equation}
To this end, we set
$
F:=\operatorname{supp}w,
$
and for each \(j\in F\), define
\begin{equation}
\label{eq:GarlingNec-rho-d}
\rho_j
:=
\left(
\frac{|w_j|}{\|w\|_{\ell^t}}
\right)^{
1/\left(\frac1q-\frac1r\right)
},
\qquad
a_j
:=
\rho_j^{-1/r}\frac{w_j}{|w_j|}.
\end{equation}
Since \(r>q\), we see that for each $j\in F$,
$
0<\rho_j\leq1.
$
Thus, the families \(\{a_j\}_{j\in F}\) and
\(\{\rho_j\}_{j\in F}\) satisfy the hypotheses of
Lemma~\ref{lem:Bernoulli-Rademacher-estimate}.
It follows from the definition of \(\rho_j\) that for every \(j\in F\),
\begin{align}\label{eq:GarlingNec-first-normalization}
\rho_j^{1/q}a_j
&=
\rho_j^{\frac1q-\frac1r}
\frac{w_j}{|w_j|}
=
\frac{|w_j|}{\|w\|_{\ell^t}}
\frac{w_j}{|w_j|}
=
\frac{w_j}{\|w\|_{\ell^t}},\qquad {\rm and}\qquad \rho_j^{1/r}a_j=
\frac{w_j}{|w_j|}.
\end{align}

We next calculate the other term appearing in
Lemma~\ref{lem:Bernoulli-Rademacher-estimate}.
Since
$
|a_j|=\rho_j^{-1/r},
$
we have
\begin{align}\label{secondequalttt}
\|
\{\rho_j^{1/2}a_j\}_{j\in F}
\|_{\ell^s}^{\,s}
=
\sum_{j\in F}
\rho_j^{s\left(\frac12-\frac1r\right)}
=
\sum_{j\in F}
\left(
\frac{|w_j|}{\|w\|_{\ell^t}}
\right)^{
s/\vartheta
}=
\sum_{j\in F}
\left(
\frac{|w_j|}{\|w\|_{\ell^t}}
\right)^t
=
1.
\end{align}

Applying Lemma~\ref{lem:Bernoulli-Rademacher-estimate}, equalities
\eqref{eq:GarlingNec-first-normalization} and \eqref{secondequalttt}, we obtain
\begin{align*}
\frac{\|\mathscr M_{\mathbf b}w\|_{\ell^q}}
     {\|w\|_{\ell^t}}
&=
\Bigg(
\sum_{j\in F}
\left|
b_j\frac{w_j}{\|w\|_{\ell^t}}
\right|^q
\Bigg)^{1/q}
\\
&=
\Bigg(
\sum_{j\in F}
\rho_j|b_ja_j|^q
\Bigg)^{1/q}
\\
&\leq
R_r\pi_r(\mathscr M_{\mathbf b})
\Big(
\|
\{\rho_j^{1/2}a_j\}_{j\in F}
\|_{\ell^s}
+
\|
\{\rho_j^{1/r}a_j\}_{j\in F}
\|_{\ell^\infty}
\Big)\\
&\leq
2R_r\pi_r(\mathscr M_{\mathbf b}).
\end{align*}
Multiplying both sides by \(\|w\|_{\ell^t}\), we complete the proof of Lemma \ref{lem:GarlingNec-t-q-estimate}.
\end{proof}

The proof of Theorem~\ref{GarlingTheorem} is divided into two distinct parts.
In the exceptional range \((p,q,r)\in E\), it provides a complete
solution to Question\hyperlink{question1}{~1}, together with the corresponding quantitative
estimates. In the non-exceptional range \((p,q,r)\notin E\),  it gives a quantitative refinement of Garling's qualitative characterization \cite[Theorems~4 and~9]{MR355665}. To the best of our knowledge, although quantitative characterizations have been stated in the literature for some special parameter ranges (see, e.g., \cite{MR355665},\cite[Proposition 4.1]{MR3886175},\cite[Theorem 3.4]{MR1342297}), the corresponding estimates in the full non-exceptional range of parameters do not seem to be explicitly available. The main subtlety in this quantitative refinement arises in the borderline case $2<p<\infty$, $1<q=p'<2$ and $1\le r<q$. In this case, Garling's qualitative result \cite[Theorem~9(ii)]{MR355665} says that the diagonal operator \(\mathscr M_{\mathbf b}\) is \(r\)-summing if and only if $\sum_{k=1}^{\infty}
|b_k|^q
\bigl(1+\log^+|b_k|^{-1}\bigr)<+\infty$. However, as observed in Remark \ref{orliczremark}, the \(q\)-th root of this
modular expression is not equivalent to the Luxemburg norm $\|\mathbf{b}\|_{\ell^{q^-}}$, and thus cannot be used as a norm-level quantitative substitute for $\pi_r(\mathscr M_\mathbf{b})$. Therefore, for completeness, in the non-exceptional range we explain below how the desired quantitative assertions can be obtained by combining Garling's qualitative characterization with the closed graph theorem and the inverse mapping theorem.

\begin{proof}[Proof of Theorem \ref{GarlingTheorem}]
Note that in each case (especially in the borderline case in which $\kappa=q^-$), $\ell^\kappa$ is always a Banach space. Now we divide the proof into two cases.

{\bf Case 1.} Suppose that $(p,q,r)\in E$.

{\it (Sufficiency.)} Consider the bilinear map
\[
\mathcal B:
(\ell^r+\ell^{p'})
\times
\ell_w^r(\ell^p)
\longrightarrow
\ell^r(\ell^r)+\ell^r(\ell^2)
\]
defined by
\[
\mathcal B
\bigl(
\mathbf b,\{x^{(m)}\}_{m\geq1}
\bigr)
:=
\bigl\{
\mathscr M_{\mathbf b}x^{(m)}
\bigr\}_{m\geq1}.
\]

We first establish the two endpoint estimates. Let
\(\mathbf b\in\ell^r\) and
\(\mathbf x=\{x^{(m)}\}_{m\geq1}\in\ell_w^r(\ell^p)\).
Since the \(j\)-th coordinate functional \(e_j^*\) belongs to
\(B_{(\ell^p)^*}\), Tonelli's theorem gives
\begin{align}
\left\|
\mathcal B(\mathbf b,\mathbf x)
\right\|_{\ell^r(\ell^r)}^r
=
\sum_{m=1}^{\infty}
\sum_{j=1}^{\infty}
|b_jx_j^{(m)}|^r
=
\sum_{j=1}^{\infty}
|b_j|^r
\sum_{m=1}^{\infty}
|e_j^*(x^{(m)})|^r
\leq
\|\mathbf b\|_{\ell^r}^r
\|\mathbf x\|_{\ell_w^r(\ell^p)}^r.
\label{eq:bilinear-endpoint-r}
\end{align}
Thus,
$
\mathcal B:
\ell^r\times\ell_w^r(\ell^p)
\longrightarrow
\ell^r(\ell^r)
$
is bounded with norm at most \(1\).

Let
\(\mathbf b\in\ell^{p'}\), and let
\(\{r_j\}_{j\geq1}\) be the standard Rademacher system on
\([0,1]\), which satisfies
\begin{equation}\label{orthogonality}
\int_0^1 r_i(\tau)r_j(\tau)\,d\tau=\delta_{ij},
\qquad i,j\geq1,
\end{equation}
where $\delta_{ij}$ denotes the Kronecker delta.  For each
\(\tau\in[0,1]\), define
\[
\varphi_\tau(x)
:=
\sum_{j=1}^{\infty}
r_j(\tau)b_jx_j,
\qquad x\in\ell^p.
\]
By H\"older's inequality, the series converges absolutely and
\begin{align}\label{dualityestima}
\|\varphi_\tau\|_{(\ell^p)^*}
=
\|\{r_j(\tau)b_j\}_{j\geq1}\|_{\ell^{p'}}
=
\|\mathbf b\|_{\ell^{p'}}.
\end{align}
Since \(r>2\), we apply the Rademacher orthogonality relation \eqref{orthogonality} and \eqref{dualityestima} to see that
\begin{align}
\left\|
\mathcal B(\mathbf b,\mathbf x)
\right\|_{\ell^r(\ell^2)}^r
&=
\sum_{m=1}^{\infty}
\Bigg(
\sum_{j=1}^{\infty}
|b_jx_j^{(m)}|^2
\Bigg)^{r/2}
\notag\\
&=
\sum_{m=1}^{\infty}
\Bigg(
\int_0^1
\bigg|
\sum_{j=1}^{\infty}
r_j(\tau)b_jx_j^{(m)}
\bigg|^2
\,d\tau
\Bigg)^{r/2}
\notag\\
&\leq
\sum_{m=1}^{\infty}
\int_0^1
|\varphi_\tau(x^{(m)})|^r
\,d\tau
\notag\\
&=
\int_0^1
\sum_{m=1}^{\infty}
|\varphi_\tau(x^{(m)})|^r
\,d\tau
\notag\\
&\leq
\|\mathbf b\|_{\ell^{p'}}^r
\|\mathbf x\|_{\ell_w^r(\ell^p)}^r.
\label{eq:bilinear-endpoint-two}
\end{align}
Therefore,
$
\mathcal B:
\ell^{p'}\times\ell_w^r(\ell^p)
\longrightarrow
\ell^r(\ell^2)
$
is bounded with norm at most \(1\).

Applying the bilinear complex interpolation theorem
\cite[Theorem~4.4.1]{MR0482275}, we obtain
\[
\mathcal B:
[\ell^r,\ell^{p'}]_\vartheta
\times
[\ell_w^r(\ell^p),\ell_w^r(\ell^p)]_\vartheta
\longrightarrow
[\ell^r(\ell^r),\ell^r(\ell^2)]_\vartheta
\]
with norm at most \(1\). By the standard complex interpolation formulas
\cite[Theorems~5.1.1 and~5.1.2]{MR0482275}, together with
\eqref{kappa-interpolation-relation} and the definition of
\(\vartheta\), we have
\[
[\ell^r,\ell^{p'}]_\vartheta
=
\ell^\kappa\qquad {\rm and}\qquad [\ell^r(\ell^r),\ell^r(\ell^2)]_\vartheta
=
\ell^r\bigl([\ell^r,\ell^2]_\vartheta\bigr)=
\ell^r(\ell^q).
\]
Consequently,
$
\mathcal B:
\ell^\kappa\times\ell_w^r(\ell^p)
\longrightarrow
\ell^r(\ell^q)
$
is bounded with norm at most \(1\). Hence, for every
\(\mathbf b\in\ell^\kappa\) and every
\(\mathbf x=\{x^{(m)}\}_{m\geq1}\in\ell_w^r(\ell^p)\),
\begin{align}
\Bigg(
\sum_{m=1}^{\infty}
\|\mathscr M_{\mathbf b}x^{(m)}\|_{\ell^q}^r
\Bigg)^{1/r}
\leq
\|\mathbf b\|_{\ell^\kappa}
\|\mathbf x\|_{\ell_w^r(\ell^p)}
=
\|\mathbf b\|_{\ell^\kappa}
\sup_{\varphi\in B_{(\ell^p)^*}}
\Bigg(
\sum_{m=1}^{\infty}
|\varphi(x^{(m)})|^r
\Bigg)^{1/r}.
\label{eq:bilinear-interpolated-estimate}
\end{align}
This implies that
\(\mathscr M_{\mathbf b}\in\Pi_r(\ell^p,\ell^q)\), and
\begin{equation}
\label{general-upper-kappa}
\pi_r\left(
\mathscr M_{\mathbf b}:\ell^p\to\ell^q
\right)
\leq
\|\mathbf b\|_{\ell^\kappa}.
\end{equation}

{\it (Necessity.)} Assume that
$
\mathscr M_{\mathbf b}\in\Pi_r(\ell^p,\ell^q).
$
Define \(s\) by
$
1/s
=
1/2-1/p',
$
and put
$
t:=s/\vartheta.
$
By \eqref{kappa-interpolation-relation},
\begin{align}\label{eq:GarlingNec-exponent-relation}
\frac1q-\frac1t
=
\frac1q
-
\vartheta\left(
\frac12-\frac1{p'}
\right)=
\frac{\vartheta}{p'}
+
\frac{1-\vartheta}{r}
=
\frac1\kappa.
\end{align}
By \eqref{eq:GarlingNec-exponent-relation}, $t>q$ and
$
1/\kappa=1/q-1/t.
$
As a direct consequence of H\"older's inequality,
\[
\|\mathscr M_{\mathbf b}\|_{\mathcal L(\ell^t,\ell^q)}
\simeq
\|\mathbf b\|_{\ell^\kappa}.
\]
Therefore, by Lemma \ref{lem:GarlingNec-t-q-estimate},
\[
\|\mathbf b\|_{\ell^\kappa}
\lesssim
\|\mathscr M_{\mathbf b}\|_{\mathcal L(\ell^t,\ell^q)}
\lesssim
\pi_r(\mathscr M_{\mathbf b}).
\]

{\bf Case 2.} Suppose that $(p,q,r)\notin E$.

To show the upper estimate, we first apply
\cite[Theorems 4 and 9]{MR355665} to see that for every scalar sequence \(\mathbf{b}\), if $\mathbf{b}\in \ell^\kappa$, then
$\mathscr M_\mathbf{b}\in \Pi_r(\ell^p,\ell^q)$.
Therefore, the linear map
\[
\Gamma:\ell^\kappa\to \Pi_r(\ell^p,\ell^q),
\qquad
\Gamma(\mathbf{b}):=\mathscr M_\mathbf{b},
\]
is well-defined. We {\it claim} that \(\Gamma\) has closed graph. Suppose that \(\mathbf{b}^{(n)}\to \mathbf{b}\) in \(\ell^\kappa\), and that
$
\mathscr M_{\mathbf{b}^{(n)}}\to T\
\text{in } \Pi_r(\ell^p,\ell^q).
$
Then \(\mathscr M_{\mathbf{b}^{(n)}}\to T\) in the operator norm. For every
\(x=\{x_j\}\in\ell^p\) and \(j\in\mathbb N\), we have
\begin{align}\label{coordingateagume}
(Tx)_j
=
\lim_{n\to\infty}(\mathscr M_{\mathbf b^{(n)}}x)_j
=
\lim_{n\to\infty} b_j^{(n)}x_j=b_jx_j,
\end{align}
where we used the fact that \(\mathbf{b}^{(n)}\to \mathbf{b}\) in \(\ell^\kappa\) implies \(b_j^{(n)}\to b_j\) for every fixed
\(j\in\mathbb N\) (this is clear if \(\ell^\kappa\) is one of the usual
\(\ell^s\) spaces, while in the case \(\kappa=q^-\) it follows from
Lemma \ref{coordinateconver}). This implies that
\(T=\mathscr M_\mathbf{b}\), and thus the graph of \(\Gamma\) is closed. By the closed
graph theorem,
\begin{align}\label{upperGarling}
\pi_r(\mathscr M_\mathbf{b})\le C\|\mathbf{b}\|_{\ell^\kappa},
\qquad \mathbf{b}\in \ell^\kappa.
\end{align}

To show the lower estimate, we define
$
\mathcal D:=
\{\mathscr M_\mathbf{b}:\mathscr M_\mathbf{b}\in \Pi_r(\ell^p,\ell^q)\}.
$
We regard \(\mathcal D\) as a subspace of \(\Pi_r(\ell^p,\ell^q)\), equipped
with the inherited norm \(\pi_r\). We {\it claim} that \(\mathcal D\) is closed. Indeed, if \(\mathscr M_{\mathbf{b}^{(n)}}\to T\) in \(\Pi_r(\ell^p,\ell^q)\), then a similar coordinate argument, with $b_j:=(Te_j)_j$, shows that \eqref{coordingateagume} holds for every
\(x=\{x_j\}\in\ell^p\) and \(j\in\mathbb N\).
This implies that \(T=\mathscr M_\mathbf{b}\in\mathcal D\).
Therefore, \(\mathcal D\) is closed in \(\Pi_r(\ell^p,\ell^q)\). Consequently,
\(\mathcal D\) is a Banach space. By \eqref{upperGarling} and \cite[Theorems 4 and 9]{MR355665}, the map $\Gamma(\mathbf{b}):=\mathscr M_\mathbf{b}$
is a bounded bijective linear operator from $\ell^\kappa$ to $\mathcal{D}$. By the inverse mapping theorem,
\(\Gamma^{-1}:\mathcal D\to \ell^\kappa\) is bounded. Therefore,
\[
\|\mathbf{b}\|_{\ell^\kappa}\le C\pi_r(\mathscr M_\mathbf{b}),
\qquad \mathbf{b}\in \ell^\kappa.
\]
This finishes the proof of Theorem \ref{GarlingTheorem}.
\end{proof}

%

\section{Two general principles for absolutely summing operators}\label{Sect3}
\setcounter{equation}{0}

In this section, we establish two useful transference principles that reduce $r$-summability for the operators under consideration to the corresponding problem for scalar diagonal operators, thereby providing the passage from Theorem \ref{GarlingTheorem} to the $r$-summing criteria for Carleson embeddings and Hankel operators.
\begin{definition}\label{def:lpsums}
Let $
\{X_k\}_{k\ge 1}$ be a sequence of Banach spaces and let $1\le p\leq\infty$. The $\ell^{p}$-direct sum of $\{X_k\}_{k\ge 1}$ is the vector space
\[
\oplusr{k\ge 1}{\ell^{p}} X_k
:= \left\{\{x_k\}_{k\ge 1}  :\ x_k\in X_k,\ \  \big\|\{\|x_k\|_{X_k}\}_{k\ge 1}\big\|_{\ell^p}<\infty \right\},
\]
with coordinatewise addition and scalar multiplication, equipped with the norm
\[
\big\|\{x_k\}_{k\ge 1}\big\|_{\oplusr{k\ge 1}{\ell^{p}} X_k}
:=\big\|\big\{\|x_k\|_{X_k}\big\}_{k\ge 1}\big\|_{\ell^p}.
\]
\end{definition}

\subsection{A general principle in the sufficiency direction}
\begin{proposition}\label{lem:block-diagonal-transfer2}
	Let \(1\le p\leq \infty,1\le q<\infty\), and $1\leq r<\infty$. Let \(\{X_k\}_{k\ge1}\) and \(\{Y_k\}_{k\ge1}\) be sequences of Banach spaces. For each \(k\ge1\), let \( T_k : X_k \longrightarrow Y_k \) be a $1$-summing operator.  Set $\mathcal{T}:=\{\pi_1(T_k)\}_{k\geq 1}$.
	If the diagonal operator $\mathscr{M}_\mathcal{T}:\ell^p\rightarrow \ell^q$
	is \(r\)-summing, then the block-diagonal operator
	\[
	T : \{x_k\}_{k\geq 1} \in \oplusr{k\ge 1}{\ell^p} X_k \longrightarrow \{ T_k(x_k) \}_{k\geq 1} \in \oplusr{k\ge 1}{\ell^q} Y_k
	\]
	is \(r\)-summing and satisfies
	\[
	\pi_r(T)\le \pi_r(\mathscr{M}_\mathcal{T}).
	\]
\end{proposition}
\begin{proof}
Set
\[
X:=\oplusr{k\geq1}{\ell^p}X_k,
\qquad
Y:=\oplusr{k\geq1}{\ell^q}Y_k.
\]

We first verify that \(T\) is well defined and bounded. Indeed, for every $x=\{x_k\}_{k\geq1}\in X$,
\begin{align*}
\Bigg(
\sum_{k=1}^{\infty}\|T_k(x_k)\|_{Y_k}^{q}
\Bigg)^{1/q}
&\leq
\Bigg(
\sum_{k=1}^{\infty}
\pi_1(T_k)^{q}\|x_k\|_{X_k}^{q}
\Bigg)^{1/q}
=
\left\|
\mathscr M_{\mathcal T}
\bigl(\{\|x_k\|_{X_k}\}_{k\geq1}\bigr)
\right\|_{\ell^q}
\leq
\|\mathscr M_{\mathcal T}\|\|x\|_X.
\end{align*}
Thus, \(T:X\to Y\) is well defined and
bounded.

By Pietsch's domination theorem
\cite[Theorem~2.12]{MR1342297}, for every \(k\geq1\), there exists a regular Borel
probability measure \(\mu_k\) on the weak-$*$ compact unit ball
\(B_{X_k^*}\) such that for each $x\in X_k$,
\begin{equation}\label{eq:Pietsch-block-transfer3}
\|T_k(x)\|_{Y_k}
\leq
\pi_1(T_k)
\int_{B_{X_k^*}}
|x^*(x)|\,d\mu_k(x^*).
\end{equation}

Fix \(N\in\mathbb N\) and consider the finite product probability
space
\[
(\Omega_N,\mu^{(N)})
:=
\prod_{k=1}^{N}(B_{X_k^*},\mu_k).
\]
For every
$
\omega=(x_1^*,\ldots,x_N^*)\in\Omega_N,
$
 define
$
A_{\omega,N}:X\longrightarrow\ell^p
$
by
\[
A_{\omega,N}\bigl(\{x_k\}_{k\geq1}\bigr)
:=
\bigl(
x_1^*(x_1),\ldots,x_N^*(x_N),0,0,\ldots
\bigr).
\]
It is immediate that for every $1\leq p\leq \infty$, we have $\|A_{\omega,N}\|_{\mathcal L(X,\ell^p)}\leq 1$.

Define the \(N\)-th block truncation of \(T\) by
\[
T^{(N)}(x)
:=
(T_1(x_1),\ldots,T_N(x_N),0,0,\ldots),
\qquad
x=\{x_k\}_{k\geq1}\in X.
\]
We apply \eqref{eq:Pietsch-block-transfer3} to see that for every \(1\leq k\leq N\),
\[
\|T_k(x_k)\|_{Y_k}
\leq
\int_{\Omega_N}
\pi_1(T_k)|x_k^*(x_k)|\,d\mu^{(N)}(\omega).
\]
This, together with Minkowski's integral inequality, yields
\begin{align}
\|T^{(N)}(x)\|_Y
=
\Bigg(
\sum_{k=1}^{N}\|T_k(x_k)\|_{Y_k}^{q}
\Bigg)^{1/q}
&\leq
\Bigg(
\sum_{k=1}^{N}
\Bigg(
\int_{\Omega_N}
\pi_1(T_k)|x_k^*(x_k)|\,d\mu^{(N)}(\omega)
\Bigg)^{q}
\Bigg)^{1/q}
\notag\\
&\leq
\int_{\Omega_N}
\Bigg(
\sum_{k=1}^{N}
\pi_1(T_k)^{q}|x_k^*(x_k)|^{q}
\Bigg)^{1/q}
\,d\mu^{(N)}(\omega)
\notag\\
&=
\int_{\Omega_N}
\left\|
\mathscr M_{\mathcal T}A_{\omega,N}(x)
\right\|_{\ell^q}
\,d\mu^{(N)}(\omega).
\label{eq:finite-block-domination-transfer3}
\end{align}
Let \(m\in\mathbb N\) and let
$
x^{(1)},\ldots,x^{(m)}\in X.
$
Applying the inequality \eqref{eq:finite-block-domination-transfer3} and  Minkowski's integral inequality, we obtain
\begin{align}
\Bigg(
\sum_{j=1}^{m}
\|T^{(N)}(x^{(j)})\|_Y^{r}
\Bigg)^{1/r}
&\leq
\int_{\Omega_N}
\Bigg(
\sum_{j=1}^{m}
\left\|
\mathscr M_{\mathcal T}
A_{\omega,N}(x^{(j)})
\right\|_{\ell^q}^{r}
\Bigg)^{1/r}
\,d\mu^{(N)}(\omega).
\label{eq:Minkowski-r-transfer3}
\end{align}

For every fixed \(\omega\in\Omega_N\), since \(\mathscr M_{\mathcal T}\) is \(r\)-summing from $\ell^p$ to $\ell^q$, we have
\begin{align}
&\Bigg(
\sum_{j=1}^{m}
\left\|
\mathscr M_{\mathcal T}
A_{\omega,N}(x^{(j)})
\right\|_{\ell^q}^{r}
\Bigg)^{1/r}\leq
\pi_r(\mathscr M_{\mathcal T})
\sup_{\xi\in B_{(\ell^p)^*}}
\Bigg(
\sum_{j=1}^{m}
\left|
\xi(A_{\omega,N}(x^{(j)}))
\right|^{r}
\Bigg)^{1/r}.
\label{eq:apply-diagonal-r-summing-transfer3}
\end{align}
Since for every \(\xi\in B_{(\ell^p)^*}\), the functional
$
\xi\circ A_{\omega,N}\in B_{X^*}
$, we have
\begin{align}
\sup_{\xi\in B_{(\ell^p)^*}}
\Bigg(
\sum_{j=1}^{m}
\left|
\xi(A_{\omega,N}(x^{(j)}))
\right|^{r}
\Bigg)^{1/r}
&\leq
\sup_{\varphi\in B_{X^*}}
\Bigg(
\sum_{j=1}^{m}
|\varphi(x^{(j)})|^{r}
\Bigg)^{1/r}.
\label{eq:weak-r-transfer3}
\end{align}

Combining \eqref{eq:Minkowski-r-transfer3},
\eqref{eq:apply-diagonal-r-summing-transfer3}, and
\eqref{eq:weak-r-transfer3}, we obtain
\begin{align}
\Bigg(
\sum_{j=1}^{m}
\|T^{(N)}(x^{(j)})\|_Y^{r}
\Bigg)^{1/r}
&\leq
\pi_r(\mathscr M_{\mathcal T})
\sup_{\varphi\in B_{X^*}}
\Bigg(
\sum_{j=1}^{m}
|\varphi(x^{(j)})|^{r}
\Bigg)^{1/r}.
\label{eq:uniform-truncated-transfer3}
\end{align}

Finally, letting \(N\to\infty\) in
\eqref{eq:uniform-truncated-transfer3}, we conclude that
\[
\Bigg(
\sum_{j=1}^{m}
\|T(x^{(j)})\|_Y^{r}
\Bigg)^{1/r}
\leq
\pi_r(\mathscr M_{\mathcal T})
\sup_{\varphi\in B_{X^*}}
\Bigg(
\sum_{j=1}^{m}
|\varphi(x^{(j)})|^{r}
\Bigg)^{1/r}.
\]
Since \(m\) and \(x^{(1)},\ldots,x^{(m)}\) are arbitrary, we see that \(T\) is
\(r\)-summing with
$
\pi_r(T)
\leq
\pi_r(\mathscr M_{\mathcal T}).
$
This finishes the proof of Proposition \ref{lem:block-diagonal-transfer2}.
\end{proof}

\subsection{A general principle in the necessity direction}

\begin{lemma}\label{lem:synthesis-hpa}
Let \( 1\leq p <\infty\), \(\alpha>-1\), $b>n+(\alpha+1)/p$, and let $\{a_k\}$ be a separated sequence.
For $\{c_k\}\in \ell^p$, define
\[
S(\{c_k\}):=\sum_{k=1}^\infty c_k h_{p,a_k}^b.
\]
Then \(S\) is a bounded operator from $\ell^p$ to $A_\alpha^p$.
\end{lemma}

\begin{proof}
The argument is standard (see \cite[Theorem 2.30 and p.~72]{MR2115155} for the proof).
\end{proof}

%

\begin{proposition}\label{proq2pr}
Let $ 1\leq p,q<\infty$, $1\leq r<\infty$, $\alpha>-1$, $b>n+(\alpha+1)/p$, and $\delta>0$. Let $\{a_k\}$ be a separated sequence and $\mu$ be a positive Borel measure on $\mathbb{B}_n$. Assume that $ T:A^p_\alpha\to L^q(d\mu) $ is an $r$-summing operator. Set \[ b_k:=\left(\int_{B(a_k,\delta)} |T(h_{p,a_k}^b)(z)|^q\,d\mu(z)\right)^{1/q},\qquad \mathbf{b}:=\{b_k\}_{k\geq 1}. \]
Then the diagonal operator $\mathscr{M}_\mathbf{b}:\ell^p\rightarrow \ell^q$ is $r$-summing and
\[ \pi_r(\mathscr{M}_\mathbf{b})\lesssim \pi_r(T). \] \end{proposition}
\begin{proof}
Since removing the indices for which \(b_k=0\) does not affect the conclusion, we may assume that \(b_k>0\) for all \(k\). Consider the Banach space
\[
Y:=\oplusr{k\ge 1}{\ell^q} L^q(B(a_k,\delta),d\mu).
\]

For each \(k\), define
$
u_k:=b_k^{-1}T(h_{p,a_k}^b)\chi_{B(a_k,\delta)}.
$
Then \(\|u_k\|_{L^q(B(a_k,\delta),d\mu)}=1\). Let
\[
Z:=\Big\{\{c_k u_k\}_{k=1}^\infty:\{c_k\}\in \ell^q\Big\}\subset Y,
\]
and define
\[
\Phi:\ell^q\to Z,\qquad \Phi(\{c_k\}):=\{c_k u_k\}_{k=1}^\infty.
\]
Since
\[
\|\Phi(\{c_k\})\|_Y^q
=\sum_{k=1}^\infty |c_k|^q \|u_k\|_{L^q(B(a_k,\delta),d\mu)}^q
=\sum_{k=1}^\infty |c_k|^q,
\]
the map \(\Phi\) is an isometric isomorphism from \(\ell^q\) onto \(Z\). In particular, \(Z\) is a closed subspace of \(Y\).

For each \(t\in[0,1]\), define
$
D_t(\{c_k\}):=\{r_k(t)c_k\}
$
and
$\Theta_t:=S\circ D_t$. By Lemma \ref{lem:synthesis-hpa}, $S$ is bounded from $\ell^p$ to $A_\alpha^p$. Since
 \(D_t\) is an isometry on $\ell^p$, $\Theta_t$ is bounded from $\ell^p$ to $A_\alpha^p$ uniformly in $t$.
Next, define
$
R_t(f):=\{r_k(t)\chi_{B(a_k,\delta)}f\}_{k=1}^\infty,
$
where \(\{r_k\}_{k\geq1}\) is the standard Rademacher system on
\([0,1]\).
By the finite overlap property \eqref{a1},
\[
\|R_t(f)\|_Y^q
=\sum_{k=1}^\infty \int_{B(a_k,\delta)}|f(z)|^q\,d\mu(z)
\le N\int_{\mathbb B_n}|f(z)|^q\,d\mu(z).
\]

For \(c\in\ell^p\), define
\[
\widetilde{\Delta}(c):=\int_0^1 R_tT\Theta_t(c)\,dt,
\]
where the integral is taken in the Bochner sense in \(Y\). This is well defined, since for every \(t\in[0,1]\),
\[
\|R_tT\Theta_t(c)\|_Y
\le \|R_t\|\,\|T\|\,\|\Theta_t\|\,\|c\|_{\ell^p}
\lesssim \pi_r(T)\|c\|_{\ell^p},
\]
and \(t\mapsto R_tT\Theta_t(c)\) is strongly measurable. By the triangle inequality for \(\pi_r\) and the ideal property,
\begin{align*}
\pi_r(\widetilde{\Delta})
\le \int_0^1 \pi_r(R_t T \Theta_t)\,dt
\le \int_0^1 \|R_t\|\,\pi_r(T)\,\|\Theta_t\|\,dt
\lesssim \pi_r(T).
\end{align*}

Now let \(c=\{c_k\}\in c_{00}\), where \(c_{00}\) denotes the space of finitely supported scalar sequences.
By definition, the \(k\)-th coordinate of \(R_tT\Theta_t(c)\) is
\[
r_k(t) \chi_{B(a_k,\delta)}\sum_{j=1}^\infty c_j r_j(t) T(h_{p,a_j}^b).
\]
Using the orthogonality relation (see \eqref{orthogonality}) for Rademacher sequences, we obtain
\[
\widetilde{\Delta}(c)
=
\left\{c_k\,T(h_{p,a_k}^b) \chi_{B(a_k,\delta)}\right\}_{k=1}^\infty=\{b_kc_ku_k\}_{k=1}^\infty.
\]
Since \(c\in c_{00}\), the sequence \(\{b_kc_k\}_{k=1}^\infty\) has finite support, and hence belongs to \(\ell^q\). By the definition of $Z$,  we see that \(\widetilde{\Delta}(c)\in Z\) for every \(c\in c_{00}\).
Since $c_{00}$ is dense in $\ell^p$, for every \(c\in \ell^p\), there exists a sequence \(\{c^{(m)}\}\subset c_{00}\) such that \(c^{(m)}\to c\) in \(\ell^p\). Since $\widetilde{\Delta}$ is bounded from $\ell^p$ to $Y$,
$
\widetilde{\Delta}(c^{(m)})\to \widetilde{\Delta}(c)\ \text{in }Y.
$
Since \(\widetilde{\Delta}(c^{(m)})\in Z\) for every \(m\) and \(Z\) is closed in \(Y\), it follows that \(\widetilde{\Delta}(c)\in Z\). Hence,
$
\widetilde{\Delta}(\ell^p)\subset Z.
$
We may therefore define
\[
M:=\Phi^{-1}\circ \widetilde{\Delta}:\ell^p\to\ell^q.
\]
Since \(\Phi^{-1}\) is an isometry and \(\widetilde{\Delta}\) is \(r\)-summing, the ideal property yields
\[
\pi_r(M)\le \|\Phi^{-1}\|\,\pi_r(\widetilde{\Delta})
=\pi_r(\widetilde{\Delta})
\lesssim \pi_r(T).
\]
Moreover, for every \(c=\{c_k\}\in c_{00}\),
\[
M(c)
=
\Phi^{-1}\big(\{c_k\,T(h_{p,a_k}^b) \chi_{B(a_k,\delta)}\}_{k=1}^\infty\big)
=
\{b_k c_k\}_{k=1}^\infty.
\]
Thus, the formal diagonal map \(c\mapsto \{b_k c_k\}\), initially defined on \(c_{00}\), extends to a bounded operator
$
\mathscr M_\mathbf{b}:\ell^p\to\ell^q,
$
and this extension coincides with \(M\). Hence,
\[
\pi_r(\mathscr M_\mathbf{b})=\pi_r(M)\lesssim \pi_r(T).
\]
This completes the proof of Proposition \ref{proq2pr}.
\end{proof}


\section{Absolutely summing Carleson embeddings}\label{Sect4}
\setcounter{equation}{0}

In this section, we characterize those measures $\mu$ for which the Carleson embedding operator $J_\mu:A_\alpha^p \to L^q(d\mu)$ is $r$-summing.
To begin with, we summarize the boundedness characterization of the Carleson embedding operator $J_\mu:A_\alpha^p \to L^q(d\mu)$ as follows (see, for example, \cite[Theorems~A and~B]{MR3426615} and \cite[Lemma~2.5]{MR3520714}).
\begin{lemma}\label{thm:Carleson-embedding}
Let \(0<p,q<\infty\), \(\alpha>-1\), and \(0<\delta\leq 1\). Let
\(\{a_k\}\) be a \(\delta\)-lattice in \(\mathbb B_n\).
\begin{enumerate}
\item If \(p\le q\), then
$
J_\mu:A^p_\alpha\to L^q(d\mu)$
 is bounded  if and only if
$\hat{\mu}_{p,q,\alpha,\delta}\in L^\infty(\mathbb B_n,d\lambda)$, if and only if $\{\hat{\mu}_{p,q,\alpha,\delta}(a_k)\}\in\ell^\infty$.
Moreover,
\[
\|J_\mu\|_{A^p_\alpha\to L^q(d\mu)}
\simeq
\|\hat{\mu}_{p,q,\alpha,\delta}\|_{L^\infty(\mathbb B_n,d\lambda)}\simeq \|\{\hat{\mu}_{p,q,\alpha,\delta}(a_k)\}\|_{\ell^\infty}.
\]

\item If \(q<p\), then $
J_\mu:A^p_\alpha\to L^q(d\mu)$ is bounded if and only if
$\hat{\mu}_{p,q,\alpha,\delta}\in L^{pq/(p-q)}(\mathbb B_n,d\lambda)$, if and only if $\{\hat{\mu}_{p,q,\alpha,\delta}(a_k)\}\in \ell^{pq/(p-q)}$.
Moreover,
\[
\|J_\mu\|_{A^p_\alpha\to L^q(d\mu)}
\simeq
\|\hat{\mu}_{p,q,\alpha,\delta}\|_{L^{pq/(p-q)}(\mathbb B_n,d\lambda)}\simeq \|\{\hat{\mu}_{p,q,\alpha,\delta}(a_k)\}\|_{\ell^{pq/(p-q)}}.
\]
\end{enumerate}
\end{lemma}
It follows from Lemma \ref{thm:Carleson-embedding} that, for fixed
\(\alpha>-1\), the class of positive Borel measures \(\mu\) for which the Carleson
embedding
$
J_\mu:A^p_\alpha\to L^q(d\mu)
$
is bounded depends on \(p\) and \(q\) only through the ratio
$
\varpi:=\frac{q}{p}.
$
Therefore, it is called a $(\varpi,\alpha)$--Bergman--Carleson measure in \cite{MR3426615}.
In contrast, we shall see in this section that such a ratio-invariance
phenomenon no longer persists for the \(r\)-summing property of the
Carleson embedding operator
$
J_\mu:A^p_\alpha\to L^q(d\mu).
$
Instead, the summability condition involves the parameters \(p\), \(q\), and
the summing exponent in a more refined way.

\begin{lemma}\label{cc1p4}
Let $ 1\leq p,q<\infty$, $\alpha>-1$, and $0<\delta\leq 1$, and let $\{a_k\}$ be a $\delta$-lattice in $\mathbb{B}_n$.
Let $\mu$ be a positive Borel measure on $\mathbb{B}_n$.
For each $k$, set $d\nu_k:=v_\alpha(B(a_k,2\delta))^{-1}\,dv_\alpha$ and define
\[
T_{p,q}^{(k)}(f):=\frac{f|_{B(a_k,\delta)}}{v_\alpha(B(a_k,2\delta))^{1/p}}.
\]
Then $T_{p,q}^{(k)}$ is bounded from $A^1(B(a_k,2\delta),d\nu_k)$ to $L^q(B(a_k,\delta),d\mu)$ if and only if
$
\mu\bigl(B(a_k,\delta)\bigr)<\infty.
$
Moreover,
\[
\|T_{p,q}^{(k)}\|\ \simeq\ \hat{\mu}_{p,q,\alpha,\delta}(a_k).
\]
\end{lemma}
\begin{proof}
{\it (Sufficiency.)} Fix $k$ and let $f\in A^1\!\left(B(a_k,2\delta),d\nu_k\right)$.
By Lemma~\ref{Bergmaninequality}, for each $z\in B(a_k,\delta)$,
\[
|f(z)| \lesssim \frac{1}{v_\alpha(B(z,\delta))}\int_{B(z,\delta)} |f(w)|\,dv_\alpha(w).
\]
Since $B(z,\delta)\subset B(a_k,2\delta)$ and by \eqref{b1}, we obtain
\begin{equation}\label{eq:sup-est}
\sup_{z\in B(a_k,\delta)}|f(z)|
\lesssim\ \frac{1}{v_\alpha(B(a_k,2\delta))}\int_{B(a_k,2\delta)} |f(w)|\,dv_\alpha(w)
= \|f\|_{L^1(B(a_k,2\delta),d\nu_k)}.
\end{equation}
Therefore,
\begin{align*}
\|T_{p,q}^{(k)} (f)\|_{L^q(B(a_k,\delta),d\mu)}
&= v_\alpha(B(a_k,2\delta))^{-1/p}\Bigg(\int_{B(a_k,\delta)} |f(z)|^q d\mu(z)\Bigg)^{1/q}\\
&\le v_\alpha(B(a_k,2\delta))^{-1/p}\mu(B(a_k,\delta))^{1/q}\,
\sup_{z\in B(a_k,\delta)}|f(z)| \\
&\lesssim v_\alpha(B(a_k,2\delta))^{-1/p}\,\mu(B(a_k,\delta))^{1/q}\,
\|f\|_{L^1(B(a_k,2\delta),d\nu_k)}\\
&\lesssim \hat{\mu}_{p,q,\alpha,\delta}(a_k)\,
\|f\|_{L^1(B(a_k,2\delta),d\nu_k)}.
\end{align*}
This proves the sufficiency and yields the claimed upper bound.

{\it (Necessity.)} Take the test function $f\equiv 1$ on $B(a_k,2\delta)$. It follows that
$
\|f\|_{L^1(B(a_k,2\delta),d\nu_k)}=1,
$
and $T_{p,q}^{(k)}(f)=v_\alpha(B(a_k,2\delta))^{-1/p}$ on $B(a_k,\delta)$. Combining this with \eqref{a2}, we deduce that
\[
\|T_{p,q}^{(k)}\| \ge \|T_{p,q}^{(k)}(f)\|_{L^q(B(a_k,\delta),d\mu)}
= v_\alpha(B(a_k,2\delta))^{-1/p}\mu(B(a_k,\delta))^{1/q}
\simeq \hat{\mu}_{p,q,\alpha,\delta}(a_k).
\]
This completes the proof of necessity.
\end{proof}

\begin{lemma}\label{cc1p5}
	Let $1\leq p,q<\infty$, $\alpha>-1$, and $0<\delta\leq 1$. Let $\{a_k\}$ be a $\delta$-lattice and $\mu$ be a positive Borel measure on $\mathbb{B}_n$. For each fixed \(k\), define
$
S_k(f):=f|_{B(a_k,\delta)}.
$
Then \(S_k\) defines a \(1\)-summing operator from
\(A^p(B(a_k,3\delta),dv_\alpha)\) to
\(L^q(B(a_k,\delta),d\mu)\) if and only if
$
\mu(B(a_k,\delta))<\infty.
$
Moreover,
\[
\pi_1(S_k)\simeq
\hat{\mu}_{p,q,\alpha,\delta}(a_k).
\]
\end{lemma}
\begin{proof}
{\it (Sufficiency.)} Suppose first that
$
\mu(B(a_k,\delta))<\infty.
$
Let $d\nu_k=v_\alpha(B(a_k,2\delta))^{-1} dv_\alpha$ be the normalized volume measure on \( B(a_k,2\delta )\). Let
$$P_1(f):=v_\alpha(B(a_k,2\delta))^{1/p} f|_{{B(a_k,2\delta)}},\quad \text{for}\  f  \in A^p(B(a_k,3\delta ),dv_\alpha),$$
$$P_2(f) := T_{p,q}^{(k)}(f),\quad \text{for}\ f\in A^1(B(a_k,2\delta ), d\nu_k).$$
Consider the commutative diagram
\[
\begin{CD}
	A^p(B(a_k,3\delta ),dv_\alpha) @>S_k>> L^q(B(a_k,\delta), d\mu) \\
	@V P_1 VV @AA P_2 A \\
	H^\infty(B(a_k,2\delta )) @>j_k>> A^1(B(a_k,2\delta ), d\nu_k)
\end{CD}
\]
where $j_k$ is the canonical inclusion, and its operator norm is bounded by 1 since
\(d\nu_k\) is a probability measure. By Lemma \ref{cc1p4}, the operator norm of $ P_2 $ is dominated by
$
\hat{\mu}_{p,q,\alpha,\delta}(a_k).
$
Applying Lemma \ref{Bergmaninequality} on the ball $B(z,\delta)$, using the inclusion $B(z,\delta)\subset B(a_k,3\delta)$ for $z\in B(a_k,2\delta)$, and invoking the local volume comparability \eqref{b1}, we obtain for every $f\in A^p(B(a_k,3\delta),dv_\alpha)$,
\begin{align*}
\sup_{z \in B(a_k,2\delta )} |f(z)| &\lesssim \sup_{z \in B(a_k,2\delta )}\left(  \frac{1}{v_\alpha(B(z,\delta))}\int_{B(a_k,3\delta )} |f(w)|^p \, dv_\alpha(w) \right)^{1/p}\lesssim v_\alpha(B(a_k,2\delta ))^{-1/p}\|f\|_{L^p(B(a_k,3\delta ),dv_\alpha)}.
\end{align*}
Hence, \( P_1 \) is bounded and \( \|P_1\|\lesssim 1 \).

To see that $j_k$ is $1$-summing, let $\{g_i\}_{i=1}^m\subset H^\infty(B(a_k,2\delta))$.
Since $\nu_k(B(a_k,2\delta))=1$,
\begin{align*}
\sum_{i=1}^m \bigl\|j_k g_i\bigr\|_{A^1(B(a_k,2\delta),d\nu_k)}
=\int_{B(a_k,2\delta)} \sum_{i=1}^m |g_i(z)|\,d\nu_k(z)
\le \sup_{z\in B(a_k,2\delta)} \sum_{i=1}^m |g_i(z)|\le \sup_{\varphi\in B_{H^\infty(B(a_k,2\delta))^*}} \sum_{i=1}^m |\varphi(g_i)|,
\end{align*}
where the last inequality follows from the fact that for each $z\in B(a_k,2\delta)$, the point evaluation functional
\[
\delta_z:H^\infty(B(a_k,2\delta))\to \mathbb C,
\qquad
\delta_z(g)=g(z),
\]
belongs to the unit ball of \(H^\infty(B(a_k,2\delta))^*\).
Hence, \(j_k\) is \(1\)-summing with \(\pi_1(j_k)\le 1\).
%

Since $S_k$ admits the factorization \( S_k = P_2 \circ j_k \circ P_1 \), by the ideal property of $1$-summing operators,
\[
\pi_1(S_k) \lesssim \|P_1\| \pi_1(j_k) \|P_2\| \lesssim \hat{\mu}_{p,q,\alpha,\delta}(a_k).
\]

{\it (Necessity.)} Suppose that \(S_k\) is \(1\)-summing.
Then \(S_k\) is bounded. Since the constant function \(1\) belongs to
\(A^p(B(a_k,3\delta),dv_\alpha)\), we have
\[
\mu(B(a_k,\delta))^{1/q}
=
\|S_k(1)\|_{L^q(B(a_k,\delta),d\mu)}
\leq
\|S_k\|\,v_\alpha(B(a_k,3\delta))^{1/p}
<\infty.
\]
Moreover, since \(\|S_k\|\leq\pi_1(S_k)\), we apply the local volume
comparability \eqref{a2} to see that
\[
\pi_1(S_k)
\geq
\|S_k\|
\geq
\frac{\mu(B(a_k,\delta))^{1/q}}
     {v_\alpha(B(a_k,3\delta))^{1/p}}
\simeq
\hat{\mu}_{p,q,\alpha,\delta}(a_k).
\]

This completes the proof of Lemma \ref{cc1p5}.
\end{proof}
The following lemma establishes a discretization principle for the averaging function $\hat{\mu}_{p,q,\alpha,\delta}$.
\begin{lemma}\label{lem:averaging-discretization}
Let $1\leq p<\infty$, \(1\le q<\infty\), \(\alpha>-1\), $0<\delta\leq 1$.
Let \(\Phi\) be a Young function. Then
\(\hat{\mu}_{p,q,\alpha,\delta}\in L_\Phi(\mathbb B_n,d\lambda)\) if and only if \(\{\hat{\mu}_{p,q,\alpha,\delta}(a_k)\}\in\ell_\Phi\) \text{for some (equivalently, every)} \(\delta\)-lattice
\(\{a_k\}\) in \(\mathbb B_n\).
Moreover,
\[
\|\hat{\mu}_{p,q,\alpha,\delta}\|_{L_\Phi(\mathbb B_n,d\lambda)}
\simeq
\|\{\hat{\mu}_{p,q,\alpha,\delta}(a_k)\}\|_{\ell_\Phi}.
\]
An analogous statement holds with
\(L_\Phi(\mathbb B_n,d\lambda)\) and \(\ell_\Phi\) replaced by
\(L^\infty(\mathbb B_n,d\lambda)\) and \(\ell^\infty\), respectively.
\end{lemma}
\begin{proof}
For the statement regarding the Luxemburg norm, we choose a partition
\(\{D_k\}\) of \(\mathbb B_n\) such that
$
D_k\subset B(a_k,\delta).
$
For \(z\in D_k\), we have
$
B(z,\delta)\subset B(a_k,2\delta).
$
Since \(\{B(a_j,\delta)\}\) covers \(\mathbb B_n\),
$
B(a_k,2\delta)\subset \bigcup_{j\in I_k}B(a_j,\delta),
$
where
$
I_k:=\{j: B(a_j,\delta)\cap B(a_k,2\delta)\ne\emptyset\}.
$
By the separation of the lattice and the finite multiplicity property,
the number of elements in $I_k$ is bounded by a constant independent of $k$. Moreover, since $j\in I_k$ implies $d(a_j,a_k)<3\delta$, the number of $k$ such that
$j\in I_k$ is also bounded by a constant independent of $j$. If
\(j\in I_k\) and \(z\in D_k\), then
$
1-|z|^2\simeq 1-|a_k|^2\simeq 1-|a_j|^2.
$
Thus,
\[
\mu(B(z,\delta))^{1/q}
\le
\mu(B(a_k,2\delta))^{1/q}
\le
\sum_{j\in I_k}\mu(B(a_j,\delta))^{1/q}.
\]
This, together with the volume comparability $v_\alpha(B(z,\delta))\simeq v_\alpha(B(a_j,\delta))$ for all $j\in I_k$, $z\in D_k$, yields
\begin{align}\label{pointwmu}
\hat{\mu}_{p,q,\alpha,\delta}(z)
\lesssim
\sum_{j\in I_k}\hat{\mu}_{p,q,\alpha,\delta}(a_j),
\qquad z\in D_k.
\end{align}

Let
$
\eta>
\|\{\hat{\mu}_{p,q,\alpha,\delta}(a_k)\}\|_{\ell_\Phi}
$ (if the right-hand side is $\infty$, then there is nothing to prove).
By the definition
of the Luxemburg norm,
\[
\sum_j
\Phi\left(
\frac{\hat{\mu}_{p,q,\alpha,\delta}(a_j)}{\eta}
\right)
\le 1.
\]
Using \eqref{pointwmu}, the convexity of \(\Phi\), and
\(\Phi(0)=0\), we apply Jensen's inequality to see that
\[
\Phi\left(
\frac{\hat{\mu}_{p,q,\alpha,\delta}(z)}{C_1\eta}
\right)
\lesssim
\sum_{j\in I_k}
\Phi\left(
\frac{\hat{\mu}_{p,q,\alpha,\delta}(a_j)}{\eta}
\right),
\qquad z\in D_k
\]
for a sufficiently large constant \(C_1\).
Since \(\lambda(D_k)\le \lambda(B(a_k,\delta))\simeq1\), summing over \(k\)
gives
\[
\begin{aligned}
\int_{\mathbb B_n}
\Phi\left(
\frac{\hat{\mu}_{p,q,\alpha,\delta}(z)}{C_1\eta}
\right)\,d\lambda(z)
&=
\sum_k
\int_{D_k}
\Phi\left(
\frac{\hat{\mu}_{p,q,\alpha,\delta}(z)}{C_1\eta}
\right)\,d\lambda(z)  \\
&\lesssim
\sum_k\sum_{j\in I_k}
\Phi\left(
\frac{\hat{\mu}_{p,q,\alpha,\delta}(a_j)}{\eta}
\right)  \\
&\lesssim
\sum_j
\Phi\left(
\frac{\hat{\mu}_{p,q,\alpha,\delta}(a_j)}{\eta}
\right)
\le C_2
\end{aligned}
\]
for some constant \(C_2\ge1\). Since \(\Phi\) is convex and \(\Phi(0)=0\),
\[
\int_{\mathbb B_n}
\Phi\left(
\frac{\hat{\mu}_{p,q,\alpha,\delta}(z)}{C_1C_2\eta}
\right)\,d\lambda(z)
\le 1.
\]
By the definition of the Luxemburg norm,
$
\|\hat{\mu}_{p,q,\alpha,\delta}\|_{L_\Phi(\mathbb B_n,d\lambda)}
\le C_1C_2\eta.
$
Letting
$
\eta\downarrow
\|\{\hat{\mu}_{p,q,\alpha,\delta}(a_k)\}\|_{\ell_\Phi},
$
we obtain
\[
\|\hat{\mu}_{p,q,\alpha,\delta}\|_{L_\Phi(\mathbb B_n,d\lambda)}
\lesssim
\|\{\hat{\mu}_{p,q,\alpha,\delta}(a_k)\}\|_{\ell_\Phi}.
\]

Conversely, for each \(k\), the Bergman ball \(B(a_k,\delta)\) can be covered
by a uniformly bounded number of balls
$
B(\xi_{k,m},\delta/4)$, $m=1,\ldots,M,
$
where \(\xi_{k,m}\in B(a_k,\delta)\) and \(M\) is independent of \(k\). Hence,
\[
\begin{aligned}
\mu(B(a_k,\delta))^{1/q}
\le
\mu\Bigg(\bigcup_{m=1}^M B(\xi_{k,m},\delta/4)\Bigg)^{1/q}
\le
\Bigg(
\sum_{m=1}^M
\mu(B(\xi_{k,m},\delta/4))
\Bigg)^{1/q}
\le
\sum_{m=1}^M
\mu(B(\xi_{k,m},\delta/4))^{1/q}.
\end{aligned}
\]
Consequently,
\begin{align}\label{precees}
\hat{\mu}_{p,q,\alpha,\delta}(a_k)
\le
\sum_{m=1}^M
\frac{\mu(B(\xi_{k,m},\delta/4))^{1/q}}
{v_\alpha(B(a_k,\delta))^{1/p}}.
\end{align}
Moreover, if \(z\in B(\xi_{k,m},\delta/4)\), then
$
B(\xi_{k,m},\delta/4)\subset B(z,\delta),
$
and
$
v_\alpha(B(z,\delta))\simeq v_\alpha(B(a_k,\delta)).
$

Let
$
\eta>
\|\hat{\mu}_{p,q,\alpha,\delta}\|_{L_\Phi(\mathbb B_n,d\lambda)}
$ (if the right-hand side is $\infty$, then there is nothing to prove).
By the definition of the Luxemburg norm,
\[
\int_{\mathbb B_n}
\Phi\left(
\frac{\hat{\mu}_{p,q,\alpha,\delta}(z)}{\eta}
\right)\,d\lambda(z)
\le 1.
\]
Using convexity and \eqref{precees}, we deduce that
\[
\begin{aligned}
\Phi\left(
\frac{\hat{\mu}_{p,q,\alpha,\delta}(a_k)}{C_3\eta}
\right)
\leq \frac{1}{M}
\sum_{m=1}^{M}
\Phi\left(
\frac{\mu(B(\xi_{k,m},\delta/4))^{1/q}}
{C_3\eta M^{-1}\, v_\alpha(B(a_k,\delta))^{1/p}}
\right)  \lesssim
\sum_{m=1}^{M}
\int_{B(\xi_{k,m},\delta/4)}
\Phi\left(
\frac{\hat{\mu}_{p,q,\alpha,\delta}(z)}{\eta}
\right)\,d\lambda(z),
\end{aligned}
\]
for a sufficiently large
constant \(C_3\).
Note that the family \(\{B(\xi_{k,m},\delta/4)\}_{k,m}\) has finite overlap. Thus
\[
\begin{aligned}
\sum_k
\Phi\left(
\frac{\hat{\mu}_{p,q,\alpha,\delta}(a_k)}{C_3\eta}
\right)
\lesssim
\sum_{k}\sum_{m=1}^{M}
\int_{B(\xi_{k,m},\delta/4)}
\Phi\left(
\frac{\hat{\mu}_{p,q,\alpha,\delta}(z)}{\eta}
\right)\,d\lambda(z) \lesssim
\int_{\mathbb B_n}
\Phi\left(
\frac{\hat{\mu}_{p,q,\alpha,\delta}(z)}{\eta}
\right)\,d\lambda(z)
\le C_4
\end{aligned}
\]
for some constant \(C_4\ge1\). By convexity and \(\Phi(0)=0\),
\[
\sum_k
\Phi\left(
\frac{\hat{\mu}_{p,q,\alpha,\delta}(a_k)}{C_3C_4\eta}
\right)
\le 1.
\]
Thus,
$
\|\{\hat{\mu}_{p,q,\alpha,\delta}(a_k)\}\|_{\ell_\Phi}
\le C_3C_4\eta.
$
Letting
$
\eta\downarrow
\|\hat{\mu}_{p,q,\alpha,\delta}\|_{L_\Phi(\mathbb B_n,d\lambda)},
$
we obtain
\[
\|\{\hat{\mu}_{p,q,\alpha,\delta}(a_k)\}\|_{\ell_\Phi}
\lesssim
\|\hat{\mu}_{p,q,\alpha,\delta}\|_{L_\Phi(\mathbb B_n,d\lambda)}.
\]
This proves the statement regarding the Luxemburg norm.

For the statement regarding the \(L^\infty\) norm, we first note that for \(z\in D_k\), the previous argument gives
$
\hat{\mu}_{p,q,\alpha,\delta}(z)
\lesssim
\sup_j \hat{\mu}_{p,q,\alpha,\delta}(a_j).
$
Taking the essential supremum over \(z\in\mathbb B_n\), we obtain
\[
\|\hat{\mu}_{p,q,\alpha,\delta}\|_{L^\infty(\mathbb B_n,d\lambda)}
\lesssim
\|\{\hat{\mu}_{p,q,\alpha,\delta}(a_k)\}\|_{\ell^\infty}.
\]
Conversely, using the same finite covering of \(B(a_k,\delta)\) as above, we have
\begin{align}\label{combbbbi1}
\hat{\mu}_{p,q,\alpha,\delta}(a_k)
\lesssim
\max_{1\le m\le M}
\frac{\mu(B(\xi_{k,m},\delta/4))^{1/q}}
{v_\alpha(B(a_k,\delta))^{1/p}}.
\end{align}
Choose \(m_k\) for which the maximum is attained. Then, for
\(z\in B(\xi_{k,m_k},\delta/4)\), we have $
B(\xi_{k,m_k},\delta/4)\subset B(z,\delta)$ and $
v_\alpha(B(z,\delta))\simeq v_\alpha(B(a_k,\delta)).
$
Hence,
\begin{align}\label{combbbbi2}
\frac{\mu(B(\xi_{k,m_k},\delta/4))^{1/q}}
{v_\alpha(B(a_k,\delta))^{1/p}}
\lesssim
\hat{\mu}_{p,q,\alpha,\delta}(z).
\end{align}
Combining the inequalities \eqref{combbbbi1} and \eqref{combbbbi2}, we deduce that
\[
\hat{\mu}_{p,q,\alpha,\delta}(a_k)\lesssim \|\hat{\mu}_{p,q,\alpha,\delta}\|_{L^\infty(B(\xi_{k,m_k},\delta/4))}
\lesssim
\|\hat{\mu}_{p,q,\alpha,\delta}\|_{L^\infty(\mathbb B_n,d\lambda)}.
\]
Taking the supremum over \(k\), we obtain
\[
\|\{\hat{\mu}_{p,q,\alpha,\delta}(a_k)\}\|_{\ell^\infty}
\lesssim
\|\hat{\mu}_{p,q,\alpha,\delta}\|_{L^\infty(\mathbb B_n,d\lambda)}.
\]
The proof of Lemma \ref{lem:averaging-discretization} is complete.
\end{proof}

\begin{proof}[Proof of Theorem \ref{MainTheorem2}]
\eqref{carleson1} $\Leftrightarrow$ \eqref{carleson2}. It suffices to prove the equivalence for an arbitrary fixed \(\delta\in(0,1]\) and an arbitrary \(\delta\)-lattice \(\{a_k\}_{k\geq1}\) in \(\mathbb B_n\).

{\it (Sufficiency.)} For convenience, put
$
\boldsymbol{\Xi}
:=
\{\hat{\mu}_{p,q,\alpha,\delta}(a_k)\}_{k\geq1}.
$
Assume that
$
\boldsymbol{\Xi}\in\ell^\kappa.
$

By the finite overlap property \eqref{a1}, for \(p\geq1\), the
localization map
\[
M_1:
f\in A^p_\alpha
\longmapsto
\{f|_{B(a_k,3\delta)}\}_k
\in
\bigoplus\nolimits_{\ell^p}
A^p(B(a_k,3\delta),dv_\alpha)
\]
is bounded. Next, define a partition adapted to the lattice: set
\(E_1:=B(a_1,\delta)\) and, for \(j\geq2\),
\[
E_j
:=
B(a_j,\delta)
\setminus
\bigcup_{1\leq k\leq j-1}E_k.
\]
Then the sets \(\{E_j\}_{j\geq1}\) are pairwise disjoint, cover
\(\mathbb B_n\), and satisfy \(E_j\subset B(a_j,\delta)\) for all
\(j\). It follows that the gluing operator
\[
M_2:
\bigoplus\nolimits_{\ell^q}
L^q(B(a_k,\delta),d\mu)
\longrightarrow
L^q(d\mu),
\qquad
(M_2\{g_k\})(z):=g_j(z)
\quad\text{for }z\in E_j,
\]
is bounded, since
\[
\|M_2\{g_k\}\|_{L^q(d\mu)}^q
=
\sum_j\int_{E_j}|g_j(z)|^q\,d\mu(z)
\leq
\sum_j
\|g_j\|_{L^q(B(a_j,\delta),d\mu)}^q.
\]

By Lemma \ref{cc1p5}, for each \(k\), the restriction
\[
S_k:
A^p(B(a_k,3\delta),dv_\alpha)
\longrightarrow
L^q(B(a_k,\delta),d\mu),
\qquad
S_k(f):=f|_{B(a_k,\delta)},
\]
is \(1\)-summing and satisfies
$
\pi_1(S_k)
\lesssim
\hat{\mu}_{p,q,\alpha,\delta}(a_k).
$
Let \(S\) be the block-diagonal operator defined by
\[
S:
\{x_k\}
\in
\bigoplus\nolimits_{\ell^p}
A^p(B(a_k,3\delta),dv_\alpha)
\longmapsto
\{S_k(x_k)\}
\in
\bigoplus\nolimits_{\ell^q}
L^q(B(a_k,\delta),d\mu).
\]
Note that \(J_\mu\) admits the factorization
$
J_\mu=M_2\circ S\circ M_1.
$

Set
$
\mathcal{S}:=\{\pi_1(S_k)\}_{k\geq1},$ and
$
\mathcal{S}/\boldsymbol{\Xi}
:=
\{
\pi_1(S_k)/\hat{\mu}_{p,q,\alpha,\delta}(a_k)
\}_{k\geq1},
$
where the quotient is understood to be \(0\) whenever \(\hat{\mu}_{p,q,\alpha,\delta}(a_k)=0\).
Since
$
\pi_1(S_k)\lesssim\hat{\mu}_{p,q,\alpha,\delta}(a_k),
$
the diagonal operator $\mathscr M_{\mathcal{S}/\boldsymbol{\Xi}}$ is bounded on $\ell^q$. Moreover,
$
\mathscr M_{\mathcal S}
=
\mathscr M_{\mathcal{S}/\boldsymbol{\Xi}}
\circ
\mathscr M_{\boldsymbol{\Xi}}.
$

By Theorem \ref{GarlingTheorem},
$
\mathscr M_{\boldsymbol{\Xi}}
\in
\Pi_r(\ell^p,\ell^q)
$
and
$
\pi_r(\mathscr M_{\boldsymbol{\Xi}})
\lesssim
\|\boldsymbol{\Xi}\|_{\ell^\kappa}.
$
Hence, by the ideal property,
\[
\pi_r(\mathscr M_{\mathcal S})
\leq
\|\mathscr M_{\mathcal{S}/\boldsymbol{\Xi}}\|\,
\pi_r(\mathscr M_{\boldsymbol{\Xi}})
\lesssim
\|\boldsymbol{\Xi}\|_{\ell^\kappa}.
\]
Applying Proposition \ref{lem:block-diagonal-transfer2} to the
\(1\)-summing operators \(\{S_k\}_{k\geq1}\), we see that \(S\)
is \(r\)-summing and
\[
\pi_r(S)
\leq
\pi_r(\mathscr M_{\mathcal S})
\lesssim
\|\boldsymbol{\Xi}\|_{\ell^\kappa}.
\]
Therefore, by the ideal property of \(r\)-summing
operators,
\[
\pi_r(J_\mu)
\leq
\|M_2\|\pi_r(S)\|M_1\|
\lesssim
\|\boldsymbol{\Xi}\|_{\ell^\kappa}.
\]
This proves the sufficiency and the required upper estimate.

{\it (Necessity.)} For each $k$, by \eqref{a5s},
	\begin{align*}
		\int_{B(a_k, \delta)}|J_\mu(h_{p,a_k}^b)(z)|^q\,d\mu(z)=\int_{B(a_k, \delta)}|h_{p,a_k}^b(z)|^q\,d\mu(z)\simeq (1-|a_k|^2)^{-\frac{(n+1+\alpha)q}{p}}\ \mu\bigl(B(a_k,\delta)\bigr).
	\end{align*}
By taking $T=J_\mu$ in Proposition \ref{proq2pr}, we see that $\mathscr{M}_{\boldsymbol{\Xi}}$ is $r$-summing and
\begin{align}\label{JmuMbeta}
\pi_r(\mathscr{M}_{\boldsymbol{\Xi}})\lesssim \pi_r(J_\mu).
\end{align}
By Theorem \ref{GarlingTheorem},
\[
\pi_r(\mathscr{M}_{\boldsymbol{\Xi}})\gtrsim \|\{\hat{\mu}_{p,q,\alpha,\delta}(a_k)\}\|_{\ell^{\kappa}}.
\]
Combining these estimates, we obtain the necessity and the required
lower estimate.

\eqref{carleson2} $\Leftrightarrow$ \eqref{carleson3}. This follows directly from Lemma \ref{lem:averaging-discretization}.

This completes the proof of Theorem \ref{MainTheorem2}.
\end{proof}

%
%
%

\section{Absolutely summing Hankel operators}\label{Sect5}
\setcounter{equation}{0}

In this section, we characterize those symbols $f$ for which the associated Hankel operator $H^\beta_f:A^p_\alpha\to L^q(dv_\beta)$ is $r$-summing.
\subsection{The case of general symbols}
The following lemma provides a discretization of weighted IDA (semi-)norms: up to harmless enlargements of the Bergman-metric radius, the corresponding continuous $L^s(\mathbb{B}_n,d\lambda)$ norms are equivalent to their lattice $\ell^s$ counterparts.
\begin{lemma}\label{GG}
Let \(1\le q<\infty\), \(t\in\mathbb R\), and
\(\sigma,\tau,\upsilon>0\) with
$
\tau\ge \sigma+1$ and
$\upsilon\ge \tau+1$.
Let \(\Phi\) be a Young function and let \(f\in L^q_{\rm loc}(\mathbb B_n)\).
Then for every \(0<\delta\le1\) and every \(\delta\)-lattice
\(\{a_k\}\) in \(\mathbb B_n\), we have
\begin{align}\label{GG-general}
\big\|(1-|\cdot|^2)^t
G_{q,\sigma\delta}(f)\big\|_{L_\Phi(\mathbb B_n,d\lambda)}
\lesssim
\big\|\big\{
(1-|a_k|^2)^tG_{q,\tau\delta}(f)(a_k)
\big\}\big\|_{\ell_\Phi}\lesssim
\big\|(1-|\cdot|^2)^t
G_{q,\upsilon\delta}(f)\big\|_{L_\Phi(\mathbb B_n,d\lambda)}.
\end{align}
An analogous statement holds with
\(L_\Phi(\mathbb B_n,d\lambda)\) and \(\ell_\Phi\) replaced by
\(L^\infty(\mathbb B_n,d\lambda)\) and \(\ell^\infty\), respectively.
\end{lemma}
\begin{proof}

{\it(Lower inequality.)} Let $\{a_k\}$ be a $\delta$-lattice. Then the balls $\{B(a_k,\delta)\}$ cover $\mathbb{B}_n$ with finite overlap.
If $w\in B(a_k,\delta)$, then $B(w,\sigma\delta)\subset B(a_k,(\sigma+1)\delta)$ and $v(B(w,\sigma\delta))\simeq v(B(a_k,(\sigma+1)\delta))$. Then
\[
G_{q,\sigma\delta}(f)(w)\lesssim G_{q,(\sigma+1)\delta}(f)(a_k)\lesssim G_{q,\tau\delta}(f)(a_k)
\quad\text{and}\quad
(1-|w|^2)^{t}\simeq (1-|a_k|^2)^{t}.
\]
Hence, there exists a constant $\mathcal{C}_1>0$ such that
\begin{equation}\label{eq957}
\sup_{w\in B(a_k,\delta)} (1-|w|^2)^{t} G_{q,\sigma\delta}(f)(w)
\leq \mathcal{C}_1 (1-|a_k|^2)^{t} G_{q,\tau\delta}(f)(a_k).
\end{equation}

For the statement regarding the Luxemburg norm, let
$
\eta>\big\|\{(1-|a_k|^2)^{t}G_{q,\tau\delta}(f)(a_k)\}\big\|_{\ell_\Phi}
$ (if the right-hand side is $\infty$, then there is nothing to prove), then
by the definition of the Luxemburg norm,
\[
\sum_k
\Phi\left(
\frac{(1-|a_k|^2)^{t}G_{q,\tau\delta}(f)(a_k)}{\eta}
\right)
\le 1.
\]

Using the covering
\(\mathbb{B}_n=\bigcup_k B(a_k,\delta)\), the inequality \eqref{eq957} and
the monotonicity of \(\Phi\), we obtain
\begin{align*}
\int_{\mathbb{B}_n}
\Phi\left(
\frac{(1-|w|^2)^{t}G_{q,\sigma\delta}(f)(w)}
{\mathcal{C}_1\eta}
\right)\,d\lambda(w)
&\le
\sum_k\int_{B(a_k,\delta)}
\Phi\left(
\frac{(1-|w|^2)^{t}G_{q,\sigma\delta}(f)(w)}
{\mathcal{C}_1\eta}
\right)\,d\lambda(w)\\
&\le
\sum_k\int_{B(a_k,\delta)}
\Phi\left(
\frac{(1-|a_k|^2)^{t}G_{q,\tau\delta}(f)(a_k)}
{\eta}
\right)\,d\lambda(w)\\
&\lesssim
\sum_k
\Phi\left(
\frac{(1-|a_k|^2)^{t}G_{q,\tau\delta}(f)(a_k)}
{\eta}
\right)\leq \mathcal{C}_2,
\end{align*}
where \(\mathcal{C}_2\ge 1\) depends only on \(n\) and \(\delta\). Since \(\Phi\) is
convex and \(\Phi(0)=0\), we have
\[
\int_{\mathbb{B}_n}
\Phi\left(
\frac{(1-|w|^2)^{t}G_{q,\sigma\delta}(f)(w)}
{\mathcal{C}_1\mathcal{C}_2\eta}
\right)\,d\lambda(w)
\le 1.
\]
By the definition of the Luxemburg norm on
\(L_\Phi(\mathbb{B}_n,d\lambda)\), we get
\[
\|(1-|\cdot|^2)^{t}G_{q,\sigma\delta}(f)\|_{L_\Phi(\mathbb{B}_n,d\lambda)}
\le
\mathcal{C}_1\mathcal{C}_2\eta.
\]
Letting
$
\eta\downarrow
\|\{(1-|a_k|^2)^{t}G_{q,\tau\delta}(f)(a_k)\}\|_{\ell_\Phi},
$
we obtain
\[
\|(1-|\cdot|^2)^{t}G_{q,\sigma\delta}(f)\|_{L_\Phi(\mathbb{B}_n,d\lambda)}
\lesssim
\|\{(1-|a_k|^2)^{t}G_{q,\tau\delta}(f)(a_k)\}\|_{\ell_\Phi}.
\]

For the statement regarding the \(L^\infty\) norm, since the balls $\{B(a_k,\delta)\}$ cover $\mathbb{B}_n$, by the inequality \eqref{eq957}, we deduce that
\begin{align*}
\|(1-|\cdot|^2)^{t} G_{q,\sigma\delta}(f)\|_{L^\infty(\mathbb{B}_n,d\lambda)}
&\leq\sup_k\operatorname*{ess\,sup}_{w\in B(a_k,\delta)}
(1-|w|^2)^{t} G_{q,\sigma\delta}(f)(w)\\
&\lesssim
\|\{(1-|a_k|^2)^{t} G_{q,\tau\delta}(f)(a_k)\}\|_{\ell^\infty}.
\end{align*}

{\it(Upper inequality.)}
Let $\{a_k\}$ be a $\delta$-lattice. If $z\in B(a_k,\delta)$, then $B(a_k,\tau\delta)\subset B(z,\upsilon\delta)$ and $v(B(a_k,\tau\delta))\simeq v(B(z,\upsilon\delta))$. Therefore, for every $z\in B(a_k,\delta)$, we have $1-|a_k|^2\simeq 1-|z|^2$ and
\begin{equation}\label{eq957a}
G_{q,\tau\delta}(f)(a_k)\lesssim G_{q,\upsilon\delta}(f)(z).
\end{equation}

For the statement regarding the Luxemburg norm, let
$
\eta>\|(1-|\cdot|^2)^{t}G_{q,\upsilon\delta}(f)\|_{L_\Phi(\mathbb{B}_n,d\lambda)}
$ (if the right-hand side is $\infty$, then there is nothing to prove), then
by the definition of the Luxemburg norm,
\[
\int_{\mathbb{B}_n}
\Phi\left(
\frac{(1-|z|^2)^{t}G_{q,\upsilon\delta}(f)(z)}
{\eta}
\right)\,d\lambda(z)
\le 1.
\]

By the comparability \(1-|a_k|^2\simeq 1-|z|^2\) for every \(z\in B(a_k,\delta)\) and the inequality
\eqref{eq957a}, there exists a constant \(\mathcal{C}_3>0\) such that for every \(z\in B(a_k,\delta)\),
\[
(1-|a_k|^2)^{t}G_{q,\tau\delta}(f)(a_k)
\le
\mathcal{C}_3(1-|z|^2)^{t}G_{q,\upsilon\delta}(f)(z).
\]
Hence, by the monotonicity of \(\Phi\),
\[
\Phi\left(
\frac{(1-|a_k|^2)^{t}G_{q,\tau\delta}(f)(a_k)}
{\mathcal{C}_3\eta}
\right)
\le
\Phi\left(
\frac{(1-|z|^2)^{t}G_{q,\upsilon\delta}(f)(z)}
{\eta}
\right),
\qquad z\in B(a_k,\delta).
\]
Integrating over \(B(a_k,\delta)\) and using \(\lambda(B(a_k,\delta))\simeq 1\),
we obtain
\[
\Phi\left(
\frac{(1-|a_k|^2)^{t}G_{q,\tau\delta}(f)(a_k)}
{\mathcal{C}_3\eta}
\right)
\lesssim
\int_{B(a_k,\delta)}
\Phi\left(
\frac{(1-|z|^2)^{t}G_{q,\upsilon\delta}(f)(z)}
{\eta}
\right)\,d\lambda(z).
\]
Summing over \(k\), and using the finite overlap of the balls \(B(a_k,\delta)\),
we get
\begin{align*}
&\sum_k
\Phi\left(
\frac{(1-|a_k|^2)^{t}G_{q,\tau\delta}(f)(a_k)}
{\mathcal{C}_3\eta}
\right)\\
&\quad\lesssim
\sum_k\int_{B(a_k,\delta)}
\Phi\left(
\frac{(1-|z|^2)^{t}G_{q,\upsilon\delta}(f)(z)}
{\eta}
\right)\,d\lambda(z)\\
&\quad\lesssim
\int_{\mathbb{B}_n}
\Phi\left(
\frac{(1-|z|^2)^{t}G_{q,\upsilon\delta}(f)(z)}
{\eta}
\right)\,d\lambda(z)
\le \mathcal{C}_4,
\end{align*}
where \(\mathcal{C}_4\ge 1\) depends only on \(n\) and \(\delta\). Since \(\Phi\) is
convex and \(\Phi(0)=0\), we have
\[
\sum_k
\Phi\left(
\frac{(1-|a_k|^2)^{t}G_{q,\tau\delta}(f)(a_k)}
{\mathcal{C}_3\mathcal{C}_4\eta}
\right)
\le 1.
\]
By the definition of the Luxemburg norm on \(\ell_\Phi\), we get
\[
\|\{(1-|a_k|^2)^{t}G_{q,\tau\delta}(f)(a_k)\}\|_{\ell_\Phi}
\le
\mathcal{C}_3\mathcal{C}_4\eta.
\]
Letting
$
\eta\downarrow
\|(1-|\cdot|^2)^{t}G_{q,\upsilon\delta}(f)\|_{L_\Phi(\mathbb{B}_n,d\lambda)},
$
we obtain
\[
\|\{(1-|a_k|^2)^{t}G_{q,\tau\delta}(f)(a_k)\}\|_{\ell_\Phi}
\lesssim
\|(1-|\cdot|^2)^{t}G_{q,\upsilon\delta}(f)\|_{L_\Phi(\mathbb{B}_n,d\lambda)}.
\]

For the statement regarding the \(L^\infty\) norm, since the balls $\{B(a_k,\delta)\}$ cover $\mathbb{B}_n$, by  \eqref{eq957a},
\begin{align*}
\|\{(1-|a_k|^2)^{t} G_{q,\tau\delta}(f)(a_k)\}\|_{\ell^\infty}
&\lesssim
\sup_k\operatorname*{ess\,sup}_{z\in B(a_k,\delta)}
(1-|z|^2)^{t} G_{q,\upsilon\delta}(f)(z)\\
&=
\|(1-|\cdot|^2)^{t} G_{q,\upsilon\delta}(f)\|_{L^\infty(\mathbb{B}_n,d\lambda)}.
\end{align*}
Therefore, the proof of Lemma \ref{GG} is complete.
\end{proof}

\begin{definition}\label{localmean}
For any $\delta>0$, $q\geq 1$ and $f\in L^q_{\rm loc}(\mathbb{B}_n)$, we define the local \(L^q\)-mean
$$M_{q,\delta}(f)(z):=\left(\frac{1}{v(B(z,\delta))}\int_{B(z,\delta)}|f(w)|^q\,dv(w)\right)^{1/q}.$$
\end{definition}

\begin{nota}
	Let $1\leq p<\infty, 1<q<\infty, \alpha,\beta>-1$. Throughout the paper, we denote by
\begin{align}\label{defofgamma}
\gamma:=\frac{n+1+\beta}{q}-\frac{n+1+\alpha}{p}.
\end{align}
\end{nota}

The main result of this section is the following reformulation of
Theorem~\ref{Maintheorem}, which also provides a discrete
characterization in terms of \(\delta\)-lattices (recall that $\kappa$ is defined in \eqref{defofkappa} and \eqref{kappa-interpolation-relation}).
\begin{theorem}\label{mainhankel}
	Let $1\leq p<\infty, 1<q<\infty, 1\leq r<\infty, \alpha,\beta>-1$ and $f\in L^1(dv_\beta)$. Then the following statements are equivalent:
\begin{enumerate}
  \item\label{HankelAS1} $H^\beta_f\in\Pi_r(A^p_\alpha,L^q(dv_\beta));$
    \item\label{HankelAS3} $(1-|\cdot|^2)^\gamma G_{q,\delta}(f)\in L^\kappa(\mathbb{B}_n,d\lambda)$ \text{for some (equivalently, every)} $0<\delta\leq 1;$
  \item\label{HankelAS2} $\{(1-|a_k|^2)^{\gamma}G_{q,4\delta}(f)(a_k)\}\in\ell^\kappa$ \text{for some (equivalently, every)} $\delta$-lattice $\{a_k\}$ in $\mathbb{B}_n$, $0<\delta\leq 1.$
\end{enumerate}
Moreover, for every $0<\delta\leq 1$ and every $\delta$-lattice $\{a_k\}$ in $\mathbb{B}_n$,
\begin{align*}
\pi_r(H_f^\beta:A^p_\alpha\rightarrow L^q(dv_\beta))\simeq \|(1-|\cdot|^2)^\gamma G_{q,\delta}(f)\|_{L^\kappa(\mathbb{B}_n,d\lambda)}\simeq \|\{(1-|a_k|^2)^{\gamma}G_{q,4\delta}(f)(a_k)\}\|_{\ell^\kappa}.
\end{align*}
\end{theorem}

\begin{proof}
\eqref{HankelAS1} $\Leftrightarrow$ \eqref{HankelAS3}.

{\it (Sufficiency.)} We first note that the assumed integrability of
\((1-|\cdot|^2)^{\gamma}G_{q,\delta}(f)\) implies that
\(f\in L^q_{\rm loc}(\mathbb B_n)\).
Let \(\epsilon:=\delta/10\), and fix an \(\epsilon\)-lattice
\(\{a_k\}\) in \(\mathbb B_n\). For any $u\in C^1(\mathbb{B}_n)$, define
\[
D_\rho (u)(z):=|\rho(z) \bar{\partial}u(z)|
+ |\rho(z)|^{\frac{1}{2}}
|\bar{\partial}u(z) \wedge \bar{\partial} \rho(z)|,
\]
where $\rho(z):=1-|z|^2$. By \cite[(3.14) and (3.15)]{MR4229611},
there exists a decomposition $f=f_1+f_2$ such that
$f_1\in C^1(\mathbb{B}_n)$, and for all $z\in\mathbb{B}_n$,
\begin{align}\label{eq201}
	D_\rho (f_1)(z)
	\lesssim \sup_{w\in B(z,4\epsilon)} G_{q,2\epsilon}(f)(w),
\end{align}
and
\begin{align}\label{eq202}
	M_{q,\epsilon}(f_2)(z)
	\lesssim \sum_{j\in J_z} G_{q,2\epsilon}(f)(a_j),
\end{align}
where \(M_{q,\epsilon}\) is defined in Definition~\ref{localmean}, and
\[
J_z:=\{j:z\in B(a_j,2\epsilon)\}.
\]
Taking the $L^q$-average over \(B(a_k,\epsilon)\) on both sides of
\eqref{eq201} and using
\(B(z,4\epsilon)\subset B(a_k,5\epsilon)\) for
\(z\in B(a_k,\epsilon)\), we obtain
\[
M_{q,\epsilon}(D_\rho (f_1))(a_k)
\lesssim
\sup_{w\in B(a_k,5\epsilon)}G_{q,2\epsilon}(f)(w).
\]
Moreover, if \(w\in B(a_k,5\epsilon)\), then
\(B(w,2\epsilon)\subset B(a_k,8\epsilon)\), and
\(v(B(w,2\epsilon))\simeq v(B(a_k,8\epsilon))\) by \eqref{b1}
with \(\alpha=0\), which implies that
\[
\sup_{w\in B(a_k,5\epsilon)}G_{q,2\epsilon}(f)(w)
\lesssim
G_{q,8\epsilon}(f)(a_k).
\]
Therefore,
\begin{align}\label{brbb}
	M_{q,\epsilon}(D_\rho f_1)(a_k)
	\lesssim G_{q,8\epsilon}(f)(a_k).
\end{align}
Next, for \(j\in J_z\), we have
\(B(a_j,2\epsilon)\subset B(z,4\epsilon)\), and by \eqref{b1}
with \(\alpha=0\),
\(v(B(a_j,2\epsilon))\simeq v(B(z,4\epsilon))\). Hence,
\[
G_{q,2\epsilon}(f)(a_j)
\lesssim
G_{q,4\epsilon}(f)(z).
\]
Combining this with \eqref{eq202} and the finite overlap property
\eqref{a1} yields
\begin{align}\label{ebbbn1}
	M_{q,\epsilon}(f_2)(z)
	\lesssim
	G_{q,4\epsilon}(f)(z)
	\lesssim
	G_{q,8\epsilon}(f)(z).
\end{align}
Taking \(z=a_k\) in \eqref{ebbbn1} and using \eqref{brbb}, we see that
\begin{align}\label{eq3}
	M_{q,\epsilon}(D_\rho f_1)(a_k)
	+
	M_{q,\epsilon}( f_2)(a_k)
	\lesssim
	G_{q,8\epsilon}(f)(a_k).
\end{align}
Set \(d\mu:=|f_2|^qdv_\beta\) and
\(d\nu:=|D_\rho (f_1)|^qdv_\beta\). We apply \eqref{a2}, \eqref{a3}
and \eqref{eq3} to deduce that
\begin{align}\label{eq1124}
(1-|a_k|^2)^{\gamma}G_{q,8\epsilon}(f)(a_k)
&\gtrsim
(1-|a_k|^2)^{\gamma}
M_{q,\epsilon}( D_{\rho}(f_1))(a_k)\nonumber\\
&=
(1-|a_k|^2)^{\gamma}
\left(
\frac{1}{v(B(a_k,\epsilon))}
\int_{B(a_k,\epsilon)}|D_{\rho}(f_1)(w)|^q\,dv(w)
\right)^{1/q}\notag\\
&\simeq
(1-|a_k|^2)^{\gamma}
(1-|a_k|^2)^{-(n+1+\beta)/q}
\left(
\int_{B(a_k,\epsilon)}
|D_{\rho}(f_1)(w)|^q\,dv_\beta(w)
\right)^{1/q}\notag\\
&=
(1-|a_k|^2)^{-(n+1+\alpha)/p}
\nu(B(a_k,\epsilon))^{1/q}.
\end{align}
Similarly,
\begin{align}\label{eq1123}
(1-|a_k|^2)^{\gamma}G_{q,8\epsilon}(f)(a_k)
\gtrsim
(1-|a_k|^2)^{\gamma}M_{q,\epsilon}( f_2)(a_k)\simeq
(1-|a_k|^2)^{-(n+1+\alpha)/p}
\mu(B(a_k,\epsilon))^{1/q}.
\end{align}

We shall use the following estimates from \cite{MR4229611}. If \(1<p\le q\) and
$
(1-|\cdot|^2)^{\gamma}G_{q,2\epsilon}(f)
\in L^\infty(\mathbb{B}_n,d\lambda),
$
then \cite[(3.17) and (3.18)]{MR4229611} give
\begin{align}\label{claiminequality}
\|H^\beta_{f_1}(g)\|_{q,\beta}
\lesssim \|J_\nu(g)\|_{L^q(d\nu)},\quad
\|H^\beta_{f_2}(g)\|_{q,\beta}
\lesssim \|J_{\mu}(g)\|_{L^q(d\mu)}
\ {\rm for}\ {\rm every}\ g\in H^\infty.
\end{align}
If \(1<q<p\) and
$
(1-|\cdot|^2)^{\gamma}G_{q,2\epsilon}(f)
\in L^{\frac{pq}{p-q}}(\mathbb{B}_n,d\lambda),
$
then the same estimates follow from
\cite[(3.31) and (3.33)]{MR4229611}. Moreover, note that the proofs of
\cite[(3.17) and (3.18)]{MR4229611} also apply to the endpoint case
\(p=1\), \(q>1\). Hence, the above estimates remain valid in this
endpoint case.

In what follows, once the corresponding hypothesis, denoted by
\((\bf{H})\), has been verified and \(J_\mu,J_\nu\) are known to be bounded,
these estimates imply by density that \(H^\beta_{f_1}\) and \(H^\beta_{f_2}\),
initially defined on \(H^\infty\), extend uniquely to bounded operators on
\(A_\alpha^p\), and the same estimates remain valid for all
\(g\in A_\alpha^p\).

By the second inequality in Lemma~\ref{GG}, with
\(\tau=8\) and \(\upsilon=10\), we have
\begin{align}\label{epsilondelta}
\big\|\big\{
(1-|a_k|^2)^{\gamma}G_{q,8\epsilon}(f)(a_k)
\big\}\big\|_{\ell^{\kappa}}
\lesssim
\big\|(1-|\cdot|^2)^{\gamma}G_{q,10\epsilon}(f)
\big\|_{L^{\kappa}(\mathbb{B}_n,d\lambda)}
=
\big\|(1-|\cdot|^2)^{\gamma}G_{q,\delta}(f)
\big\|_{L^{\kappa}(\mathbb{B}_n,d\lambda)}.
\end{align}
It follows from \eqref{eq1124}, \eqref{eq1123} and
Theorem \ref{MainTheorem2} that both \(J_{\mu}\) and \(J_{\nu}\)
are \(r\)-summing with
\begin{align}\label{eq7}
&\pi_r(J_{\mu})+\pi_r(J_{\nu})\nonumber\\
&\lesssim
\big\|\big\{
(1-|a_k|^2)^{-(n+1+\alpha)/p}
\mu(B(a_k,\epsilon))^{1/q}
\big\}\big\|_{\ell^{\kappa}}
+
\big\|\big\{
(1-|a_k|^2)^{-(n+1+\alpha)/p}
\nu(B(a_k,\epsilon))^{1/q}
\big\}\big\|_{\ell^{\kappa}}\notag\\
&\lesssim
\big\|\big\{
(1-|a_k|^2)^{\gamma}G_{q,8\epsilon}(f)(a_k)
\big\}\big\|_{\ell^{\kappa}}.
\end{align}
We now verify the hypothesis \((\bf{H})\) required in the estimates
\eqref{claiminequality}.

(i) If \(p\leq q\), then
\begin{align}\label{case1embed}
\big\|\big\{
(1-|a_k|^2)^{\gamma}G_{q,8\epsilon}(f)(a_k)
\big\}\big\|_{\ell^\infty}
\leq
\big\|\big\{
(1-|a_k|^2)^{\gamma}G_{q,8\epsilon}(f)(a_k)
\big\}\big\|_{\ell^\kappa}
<+\infty.
\end{align}
Together with Lemma \ref{GG}, this gives
$
(1-|\cdot|^2)^\gamma G_{q,2\epsilon}(f)
\in L^\infty(\mathbb B_n,d\lambda),
$
and hence verifies the hypothesis \((\bf{H})\).

(ii) If \(q<p\), then
\[
\frac{1}{\kappa}=
\begin{cases}
\displaystyle
\frac{1}{p'}+\frac{1}{q}-\frac{1}{2}
=
\frac{1}{2}+\Big(\frac{1}{q}-\frac{1}{p}\Big)
\ \ge\ \frac{1}{q}-\frac{1}{p},
& {\rm if}\ 1\le p\le 2,\ 1<q\le 2,\\[10pt]
\displaystyle \frac{1}{q}
\ \ge\ \frac{1}{q}-\frac{1}{p},
& {\rm if}\ 2<p<\infty,\ 1<q<p'<2,\\[10pt]
\displaystyle \frac{1}{q^{-}}
\ge_{\textrm{formal}}
\frac{1}{q}
>
\frac{1}{q}-\frac{1}{p},
& {\rm if}\ 2<p<\infty,\ 1<q=p'<2,\ 1\le r<q,\\[10pt]
\displaystyle \frac{1}{q}
\ \ge\ \frac{1}{q}-\frac{1}{p},
& {\rm if}\ 2<p<\infty,\ 1<q=p'<2,\ q\le r<\infty,\\[10pt]
\displaystyle \frac{1}{\max\{p',\min\{r,q\}\}}
\ \ge\ \frac{1}{q}
\ \ge\ \frac{1}{q}-\frac{1}{p},
& {\rm if}\ 2\le p\le\infty,\ 1\le p'<q<\infty.
\end{cases}
\]
where the expression
\(\frac{1}{q^{-}}\ge_{\mathrm{formal}}\frac{1}{q}\) should not be interpreted
as a strict inequality. Instead, it should be interpreted as the inclusion
\(\ell^{q^-}\subset\ell^{q}\).

In each subcase, we have
\begin{align*}
\big\|\big\{
(1-|a_k|^2)^{\gamma}G_{q,8\epsilon}(f)(a_k)
\big\}\big\|_{\ell^{\frac{pq}{p-q}}}
\lesssim
\big\|\big\{
(1-|a_k|^2)^{\gamma}G_{q,8\epsilon}(f)(a_k)
\big\}\big\|_{\ell^\kappa}
<+\infty.
\end{align*}
Together with Lemma \ref{GG}, this gives
$
(1-|\cdot|^2)^\gamma G_{q,2\epsilon}(f)
\in L^{\frac{pq}{p-q}}(\mathbb B_n,d\lambda),
$
and hence verifies the hypothesis \((\bf{H})\).

Combining \eqref{claiminequality} and \eqref{eq7} with the definition of
\(r\)-summing operators, we conclude that \(H^\beta_f\) is \(r\)-summing with
\begin{align}\label{eq8}
\pi_r(H^\beta_{f})
&\leq
\pi_r(H^\beta_{f_1})+\pi_r(H^\beta_{f_2})\notag\\
&\lesssim
\pi_r(J_{\mu})+\pi_r(J_{\nu})\notag\\
&\lesssim
\big\|\big\{
(1-|a_k|^2)^{\gamma}G_{q,8\epsilon}(f)(a_k)
\big\}\big\|_{\ell^{\kappa}}.
\end{align}
Together with \eqref{epsilondelta}, this gives
\begin{align}\label{similarargumentdelta}
\pi_r(H^\beta_f)
\lesssim
\big\|(1-|\cdot|^2)^{\gamma}G_{q,\delta}(f)
\big\|_{L^{\kappa}(\mathbb{B}_n,d\lambda)}.
\end{align}
This proves the sufficiency and the required upper estimate.

	{\it(Necessity.)}
  We first make a standard reduction. If \(H_f^\beta\) is
bounded from \(A_\alpha^p\) to \(L^q(dv_\beta)\), then
$
f-P_\beta f=H_f^\beta(1)\in L^q(dv_\beta).
$
Moreover, \(f-P_\beta f\in L^1(dv_\beta)\), and hence \(P_\beta f\in L^1(dv_\beta)\).
Since \(P_\beta f\) is holomorphic, we have
$
H_f^\beta=H_{f-P_\beta f}^\beta
$
and \(G_{q,\delta}(f)=G_{q,\delta}(f-P_\beta f)\).
Thus, we may assume without loss of generality that
\(f\in L^q(dv_\beta)\).

Let \(\{a_k\}\) be a \(\delta\)-lattice in \(\mathbb B_n\).
	By \eqref{a5s}, we have $|h_{p,a_k}^b(z)|\simeq (1-|a_k|^2)^{-(n+1+\alpha)/p}$ for $z\in B(a_k,4\delta)$. Moreover, since \(h_{p,a_k}^b\) has no zeros on \(\mathbb B_n\), the function
\[
z\mapsto \frac{P_\beta(fh_{p,a_k}^b)(z)}{h_{p,a_k}^b(z)}
\]
is holomorphic on \(\mathbb B_n\). It follows that for $p\geq 1$ and $q>1$,
	\begin{align}\label{similarlower}
		\int_{B(a_k, 4\delta)}|H_f^\beta(h_{p,a_k}^b)(z)|^q\,dv_\beta(z)
		&=\int_{B(a_k, 4\delta)}|f(z)h_{p,a_k}^b(z)-P_\beta(fh_{p,a_k}^b)(z)|^q\,dv_\beta(z)\notag\\
		&=\int_{B(a_k, 4\delta)}\left|f(z)-\frac{P_\beta(fh_{p,a_k}^b)(z)}{h_{p,a_k}^b(z)}\right|^q |h_{p,a_k}^b(z)|^{q}\,dv_\beta(z)\notag\\
		&\simeq \frac{(1-|a_k|^2)^\beta}{(1-|a_k|^2)^{(n+1+\alpha)q/p}}\int_{B(a_k, 4\delta)}\left|f(z)-\frac{P_\beta(fh_{p,a_k}^b)(z)}{h_{p,a_k}^b(z)}\right|^q \,dv(z)\notag\\
		&\gtrsim (1-|a_k|^2)^{q\gamma}G_{q,4\delta}^q(f)(a_k).
	\end{align}

By taking $T=H_f^\beta$ and $d\mu=dv_\beta$ in Proposition \ref{proq2pr} and applying Theorem \ref{GarlingTheorem}, we conclude that if  $H^\beta_f:A^p_\alpha\to L^q(dv_\beta)$ is \(r\)-summing, then
	$\big\{(1-|a_k|^2)^{\gamma}G_{q,4\delta}(f)(a_k)\big\}\in \ell^{\kappa}$ and
	\begin{equation}\label{newadd1}
		\big\|\big\{(1-|a_k|^2)^{\gamma}G_{q,4\delta}(f)(a_k)\big\}\big\|_{\ell^{\kappa}}\lesssim \pi_r(H^\beta_f).
	\end{equation}
This, together with the first inequality in Lemma \ref{GG}, implies the necessity of the statement.

\eqref{HankelAS1} $\Leftrightarrow$ \eqref{HankelAS2}. This follows directly from the first inequality in Lemma \ref{GG} with
\eqref{similarargumentdelta} and \eqref{newadd1}.
	
The proof of Theorem \ref{mainhankel} is complete.
\end{proof}
%
%
%
%
%
%
\subsection{The case of simultaneous \(r\)-summability}
\begin{definition}
	For any $\delta>0$, $1\leq q<\infty$, and $f\in L^q_{\rm loc}(\mathbb{B}_n)$,  define the local \(L^q\)-mean oscillation $$MO_{q,\delta}(f)(z):=\left(\frac{1}{v(B(z,\delta))}\int_{B(z,\delta)}\left|f(\xi)-\frac{1}{v(B(z,\delta))}\int_{B(z, \delta)}f(w)\,dv(w)\right|^q\,dv(\xi)\right)^{1/q}.$$
\end{definition}
The following corollary characterizes those symbols $f$ for which the associated Hankel operators $H^\beta_f:A^p_\alpha\to L^q(dv_\beta)$ and $H^\beta_{\bar{f}}:A^p_\alpha\to L^q(dv_\beta)$ are simultaneously $r$-summing.
\begin{coro}\label{mainhankelb}
	Let $1\leq p<\infty, 1<q<\infty, 1\leq r<\infty, \alpha,\beta>-1$, and $f\in L^1(dv_\beta)$. Then the following statements are equivalent:
\begin{enumerate}
  \item $H^\beta_f,H^\beta_{\bar{f}}\in\Pi_r(A^p_\alpha,L^q(dv_\beta));$
  \item $(1-|\cdot|^2)^{\gamma}MO_{q,\delta}(f)\in L^{\kappa }(\mathbb{B}_n,d\lambda)$ \text{for some (equivalently, every)} $0<\delta\leq 1;$
    \item $\{(1-|a_k|^2)^{\gamma}MO_{q,4\delta}(f)(a_k)\}\in \ell^{\kappa }$ \text{for some (equivalently, every)} $\delta$-lattice $\{a_k\}$ in $\mathbb{B}_n$, $0<\delta\leq 1.$
\end{enumerate}
Moreover, for every $0<\delta\leq 1$ and every $\delta$-lattice $\{a_k\}$ in $\mathbb{B}_n$,
\begin{align*}
\pi_r(H^\beta_f)+\pi_r(H^\beta_{\bar{f}})\simeq \|(1-|\cdot|^2)^{\gamma}MO_{q,\delta}(f)\|_{L^{\kappa }(\mathbb{B}_n,d\lambda)} \simeq \|\{(1-|a_k|^2)^{\gamma}MO_{q,4\delta}(f)(a_k)\}\|_{\ell^\kappa}.
\end{align*}
\end{coro}
\begin{proof}
	It follows from \cite[(4.3) and (4.4)]{MR4229611} that for all $z\in \mathbb{B}_n$,
	\begin{align}\label{MOandG}
MO_{q,\delta}(f)(z)\simeq G_{q,\delta}(f)(z)+G_{q,\delta}(\bar{f})(z).
\end{align}
	Combining this equivalence with Theorem \ref{mainhankel}, we complete the proof of Corollary \ref{mainhankelb}.
\end{proof}
%
\subsection{The case of anti-analytic symbols}
\begin{definition}
For $f\in H(\mathbb{B}_n)$ and $z\in\mathbb{B}_n$, the complex gradient of $f$ at $z$ is defined by
$$\nabla f(z):=\left(\frac{\partial f}{\partial  {z}_1}(z),\ldots,\frac{\partial f}{\partial  {z}_n}(z)\right),$$
and the invariant gradient of $f$ at $z$ is defined by
$$\widetilde{\nabla}f(z):=\nabla (f\circ \varphi_z)(0),$$
where $\varphi_z$ denotes the involutive automorphism of $\mathbb{B}_n $ interchanging $z$ and $0$ {\rm \cite[(1.2)]{MR2115155}}.
\end{definition}
\begin{definition}
Let $\displaystyle s\in[1,\infty]\cup \Big(\bigcup_{1<q<2}\{q^{-}\}\Big)$ and $t\in \mathbb{R}$. The weighted Besov space $B_s^t$ consists of all $f\in H(\mathbb{B}_n)$ such that $(1-|\cdot|^2)^{t}\widetilde{\nabla}f\in L^{s}(\mathbb{B}_n,d\lambda) $. It is equipped with the seminorm
$$\|f\|_{B^t_s}:=\|(1-|\cdot|^2)^{t}|\widetilde{\nabla}f|\|_{L^{s}(\mathbb{B}_n,d\lambda)}.$$
When \(s=\infty\), \(B_\infty^t\) is the corresponding weighted Bloch
space. In particular, \(B_\infty^0\) is the usual Bloch space (see {\rm \cite[Chapter 3]{MR2115155}}).
\end{definition}
The following lemma establishes a mean-oscillation characterization of the weighted Besov spaces.
\begin{lemma}\label{lem:MO-Besov-qminus}
Let $t\in\mathbb{R}$ and $f\in H(\mathbb B_n)$. Let $\Phi$ be a Young function. Then
$
(1-|\cdot|^2)^{t}\widetilde{\nabla}f\in L_\Phi(\mathbb B_n,d\lambda)
$
if and only if
$
(1-|\cdot|^2)^t MO_{q,\delta}(f)\in
L_\Phi(\mathbb B_n,d\lambda)
$  \text{for some (equivalently, every)} $0<\delta\leq 1$ and \text{for some (equivalently, every)} $1\leq q<\infty.$
Moreover, for every $0<\delta\leq 1$, and $1\leq q<\infty$,
\[
\|(1-|\cdot|^2)^{t}\widetilde{\nabla}f\|_{L_\Phi(\mathbb B_n,d\lambda)}
\simeq \big\|(1-|\cdot|^2)^t MO_{q,\delta}(f)\big\|_{L_\Phi(\mathbb B_n,d\lambda)}.
\]
An analogous statement holds with
\(L_\Phi(\mathbb B_n,d\lambda)\) replaced by
\(L^\infty(\mathbb B_n,d\lambda)\).
\end{lemma}

\begin{proof}
Fix \(0<\delta\leq 1\) and \(1\leq q<\infty\).
By the estimate in \cite[p.~182, line~6]{MR2115155}, we have
\begin{align}\label{pointwisecon}
|\widetilde{\nabla}f(z)|
\lesssim
MO_{q,\delta}(f)(z),
\qquad z\in\mathbb B_n.
\end{align}
Hence, by the monotonicity of the Luxemburg norm,
\begin{align}\label{onesideinequ}
\|(1-|\cdot|^2)^t |\widetilde{\nabla}f|\|_{L_\Phi(\mathbb B_n,d\lambda)}
\lesssim
\|(1-|\cdot|^2)^t MO_{q,\delta}(f)\|_{L_\Phi(\mathbb B_n,d\lambda)}.
\end{align}

For the converse estimate, we first use the argument in the proof of
\cite[Theorem~3.6]{MR2115155}, applied to the Bergman geodesic segment
joining \(z\) and \(\xi\), to see that for $\xi\in B(z,\delta)$,
\[
|f(\xi)-f(z)|
\leq
d(z,\xi)\sup_{w\in B(z,\delta)}Q_f(w)
\leq
\delta\sup_{w\in B(z,\delta)}Q_f(w),
\]
where $Q_f(w)$ is defined in \cite[(3.1)]{MR2115155}. By
\cite[Theorem 3.1]{MR2115155},
$Q_f(w)=|\widetilde{\nabla}f(w)|$ for all $w\in \mathbb{B}_n$. Thus,
\[
|f(\xi)-f(z)|
\leq
\delta\sup_{w\in B(z,\delta)}
|\widetilde{\nabla}f(w)|.
\]
This, in combination with the elementary inequality
\begin{align*}
MO_{q,\delta}(f)(z)
&\le
2\left(
\frac1{v(B(z,\delta))}
\int_{B(z,\delta)}
|f(\xi)-f(z)|^q\,dv(\xi)
\right)^{1/q}  \le
2\sup_{\xi\in B(z,\delta)}
|f(\xi)-f(z)|,
\end{align*}
and \eqref{a3}, implies that there exists a constant $K_0>0$ such that for all $z\in\mathbb{B}_n$,
\begin{equation}\label{eq:MO-by-local-sup-grad}
(1-|z|^2)^t MO_{q,\delta}(f)(z)
\leq K_0
\sup_{w\in B(z,\delta)}
(1-|w|^2)^t |\widetilde{\nabla}f(w)|.
\end{equation}

Next, let \(\{a_k\}\) be a
\(\delta\)-lattice. For each \(k\),
we choose \(\xi_k\in B(a_k,2\delta)\) such that
\[
\sup_{w\in B(a_k,2\delta)}
(1-|w|^2)^t |\widetilde{\nabla}f(w)|
\le
2(1-|\xi_k|^2)^t |\widetilde{\nabla}f(\xi_k)|.
\]

We claim that there exist constants \(C_0>0\) and \(0<\varepsilon<\delta\), independent of \(k\), \(f\), and \(\eta\),  such that
\begin{equation}\label{eq:orlicz-submean-grad}
\Phi\Bigg(
\frac{(1-|\xi_k|^2)^t
|\widetilde{\nabla}f(\xi_k)|}
{C_0\eta}
\Bigg)
\lesssim
\int_{B(\xi_k,\varepsilon)}
\Phi\Bigg(
\frac{(1-|w|^2)^t
|\widetilde{\nabla}f(w)|}
{\eta}
\Bigg)\,d\lambda(w)
\quad {\rm for}\ {\rm every}\ \eta>0.
\end{equation}

To prove the inequality \eqref{eq:orlicz-submean-grad}, we note that
for each $k$, the function
\[
u\longmapsto
\Phi\left(
\frac{(1-|\xi_k|^2)^t
|\nabla(f\circ\varphi_{\xi_k})(u)|}
{C_0\eta}
\right)
\]
is plurisubharmonic on \(B(0,\varepsilon)\)
(see e.g., \cite[Theorem~5.6, p.~39, and Example~5.12, p.~41]{DemaillyCADG}).
By the sub-mean value inequality
(see \cite[(5.3)]{DemaillyCADG} and
\cite[Theorem and definition 4.12]{DemaillyCADG})
on the fixed ball \(B(0,\varepsilon)\),
\[
\Phi\left(
\frac{(1-|\xi_k|^2)^t
|\nabla(f\circ\varphi_{\xi_k})(0)|}
{C_0\eta}
\right)
\lesssim
\int_{B(0,\varepsilon)}
\Phi\left(
\frac{(1-|\xi_k|^2)^t
|\nabla(f\circ\varphi_{\xi_k})(u)|}
{C_0\eta}
\right)
\,dv(u).
\]
For \(u\in B(0,\varepsilon)\), put
\(w=\varphi_{\xi_k}(u)\). By \cite[Lemma 2.14]{MR2115155}, since
\(u\) stays in a fixed compact subset of \(\mathbb B_n\), we have
\[
|\nabla(f\circ\varphi_{\xi_k})(u)|
\leq
(1-|u|^2)^{-1}
|\widetilde{\nabla}(f\circ\varphi_{\xi_k})(u)|
\lesssim
|\widetilde{\nabla}(f\circ\varphi_{\xi_k})(u)|
=
|\widetilde{\nabla}f(w)|.
\]
Moreover, \(w\in B(\xi_k,\varepsilon)\), and hence
$
(1-|\xi_k|^2)^t
\simeq
(1-|w|^2)^t.
$
Choosing \(C_0\) large enough and using the invariance of \(d\lambda\),
together with
\(dv(u)\lesssim d\lambda(u)=d\lambda(w)\) on the fixed ball
\(B(0,\varepsilon)\), we obtain
\eqref{eq:orlicz-submean-grad}.

Let
$\eta>
\|(1-|\cdot|^2)^t|\widetilde{\nabla}f|
\|_{L_\Phi(\mathbb B_n,d\lambda)}$
(if the right-hand side is $\infty$, then there is nothing to prove).
Then
\[
\int_{\mathbb B_n}
\Phi\Bigg(
\frac{(1-|w|^2)^t
|\widetilde{\nabla}f(w)|}{\eta}
\Bigg)\,d\lambda(w)
\leq1.
\]

Since \(\xi_k\in B(a_k,2\delta)\), the balls
\(B(\xi_k,\varepsilon)\) have finite overlap. Thus, summing
\eqref{eq:orlicz-submean-grad} over \(k\), we obtain
\begin{align*}
\sum_k
\Phi\Bigg(
\frac{\sup_{w\in B(a_k,2\delta)}
(1-|w|^2)^t|\widetilde{\nabla}f(w)|}
{2C_0\eta}
\Bigg)
&\leq
\sum_k
\Phi\Bigg(
\frac{(1-|\xi_k|^2)^t
|\widetilde{\nabla}f(\xi_k)|}
{C_0\eta}
\Bigg)
\\
&\lesssim
\sum_k
\int_{B(\xi_k,\varepsilon)}
\Phi\Bigg(
\frac{(1-|w|^2)^t
|\widetilde{\nabla}f(w)|}{\eta}
\Bigg)
\,d\lambda(w)
\\
&\lesssim
\int_{\mathbb B_n}
\Phi\Bigg(
\frac{(1-|w|^2)^t
|\widetilde{\nabla}f(w)|}{\eta}
\Bigg)
\,d\lambda(w)
\lesssim1.
\end{align*}

Using \eqref{eq:MO-by-local-sup-grad} and the inclusion
$B(z,\delta)\subset B(a_k,2\delta)$ for all
$z\in B(a_k,\delta)$,
\begin{align*}
\int_{\mathbb B_n}
\Phi\left(
\frac{(1-|z|^2)^tMO_{q,\delta}(f)(z)}
{2K_0C_0\eta}
\right)\,d\lambda(z)
&\le
\sum_k\int_{B(a_k,\delta)}
\Phi\left(
\frac{
\sup_{w\in B(a_k,2\delta)}
(1-|w|^2)^t|\widetilde{\nabla}f(w)|
}
{2C_0\eta}
\right)\,d\lambda(z)\\
&\lesssim
\sum_k
\Phi\left(
\frac{
\sup_{w\in B(a_k,2\delta)}
(1-|w|^2)^t|\widetilde{\nabla}f(w)|
}
{2C_0\eta}
\right)
\le C_1.
\end{align*}
for some constant \(C_1\geq1\) independent of \(f\) and \(\eta\).
Since \(\Phi\) is convex and \(\Phi(0)=0\),
\[
\int_{\mathbb B_n}
\Phi\left(
\frac{(1-|z|^2)^t MO_{q,\delta}(f)(z)}
{2K_0C_0C_1\eta}
\right)
\,d\lambda(z)
\leq1.
\]

By the definition of the Luxemburg norm,
\[
\|(1-|\cdot|^2)^tMO_{q,\delta}(f)
\|_{L_\Phi(\mathbb B_n,d\lambda)}
\lesssim\eta.
\]
Letting
$
\eta\downarrow
\|(1-|\cdot|^2)^t|\widetilde{\nabla}f|
\|_{L_\Phi(\mathbb B_n,d\lambda)},
$
we obtain the converse inequality of \eqref{onesideinequ}.

The statement regarding the $L^\infty$ norm follows directly from the inequalities \eqref{pointwisecon} and \eqref{eq:MO-by-local-sup-grad}. This completes the proof of Lemma \ref{lem:MO-Besov-qminus}.
\end{proof}

The following corollary characterizes those holomorphic functions \(f\) for which the Hankel operator with anti-analytic symbol \(\bar f\),
\(H^\beta_{\bar f}:A^p_\alpha\to L^q(dv_\beta)\), is \(r\)-summing.
\begin{coro}\label{mainhankelc}
	Let $1\leq p<\infty, 1<q<\infty, 1\leq r<\infty, \alpha,\beta>-1$, and $f\in A^1_\beta$. Then the following statements are equivalent:
\begin{enumerate}
  \item $H^\beta_{\bar{f}}\in \Pi_r(A^p_\alpha, L^q(dv_\beta));$
  \item $f\in B_{\kappa}^{\gamma};$
  \item $(1-|\cdot|^2)^\gamma MO_{q,\delta}(f)\in L^\kappa(\mathbb B_n,d\lambda)$ \text{for some (equivalently, every)} $0<\delta\leq 1;$
  \item $\{(1-|a_k|^2)^\gamma MO_{q,4\delta}(f)(a_k)\}\in \ell^\kappa$ \text{for some (equivalently, every)} $\delta$-lattice $\{a_k\}$ in $\mathbb{B}_n$, $0<\delta\leq 1$.
\end{enumerate}
Moreover, for every $0<\delta\leq 1$ and every $\delta$-lattice $\{a_k\}$ in $\mathbb{B}_n$,
\begin{align*}
\pi_r(H^\beta_{\bar{f}}:A^p_\alpha\rightarrow L^q(dv_\beta))\simeq \|f\|_{B_{\kappa}^{\gamma}}\simeq \|(1-|\cdot|^2)^\gamma MO_{q,\delta}(f)\|_{L^\kappa(\mathbb B_n,d\lambda)}\simeq \|\{(1-|a_k|^2)^\gamma MO_{q,4\delta}(f)(a_k)\}\|_{\ell^\kappa}.
\end{align*}
\end{coro}
\begin{proof}
Fix $0<\delta\leq 1$. Since \(f\in H(\mathbb B_n)\), we have \(G_{q,\delta}(f)\equiv 0\) and \(H_f^\beta=0\) on
\(A_\alpha^p\). This, together with Lemma \ref{lem:MO-Besov-qminus} and Corollary~\ref{mainhankelb}, completes the proof of Corollary \ref{mainhankelc}.
\end{proof}

Hahn and Youssfi \cite{MR1123807} proved that, when \(n>1\), the unweighted M\"{o}bius invariant Besov space \(B_s^0\) is nontrivial if and only if \(s>2n\). Equivalently, \(B_s^0=\mathbb C\) for \(1<s\le 2n\). See also \cite{MR1013987,MR2115155,MR1113389} for related discussion. The following lemma provides the corresponding characterization of when a weighted Besov space consists only of constants.
\begin{lemma}\label{constantcharBesov}
Let \(t\in\mathbb R\). For
$
s\in[1,\infty)\cup\{s_0^{-}:1<s_0<2\},
$
define
\[
s_\natural:=
\begin{cases}
s, & s\in[1,\infty),\\
s_0, & s=s_0^{-}\ \text{ for some }1<s_0<2.
\end{cases}
\]
Then
\begin{enumerate}
\item If \(s\in[1,\infty)\cup \{s_0^{-}:1<s_0<2\}\), then
\[
B_s^t=\mathbb C
\Longleftrightarrow
\begin{cases}
s_\natural(t+1)\le 1, & n=1,\\[1mm]
s_\natural\left(t+\frac12\right)\le n, & n\ge2.
\end{cases}
\]
\item If \(s=\infty\), then
\[
B_\infty^t=\mathbb C
\Longleftrightarrow
\begin{cases}
t<-1, & n=1,\\[1mm]
t<-\dfrac12, & n\ge2.
\end{cases}
\]
\end{enumerate}
\end{lemma}

\begin{proof}
It is clear that all constant functions belong to \(B_s^t\). It remains to identify the parameter ranges in which every element of \(B_s^t\)
is constant, and to show that outside these ranges \(B_s^t\) contains a
non-constant function.

{\bf Case 1.} If \(n=1\), then we recall from \cite[Definition~2.11 and Lemma~2.13]{MR2115155} that
\begin{align}\label{nablaandprime}
|\widetilde{\nabla}f(z)|
=(1-|z|^2)|f'(z)|.
\end{align}
Let \(f\in H(\mathbb D)\) be non-constant. Write
$
f(z)=a_0+\sum_{k=m}^\infty a_k z^k,
$
where \(m\ge1\) and \(a_m\ne0\). Then
$
f'(z)=\sum_{k=m}^\infty k a_k z^{k-1}.
$
Substituting \(z=re^{i\theta}\), multiplying by \(e^{-i(m-1)\theta}\), and
integrating over \([0,2\pi]\), we see that
\[
m a_m r^{m-1}
=
\frac1{2\pi}\int_0^{2\pi}
e^{-i(m-1)\theta}f'(re^{i\theta})\,d\theta .
\]
By H\"{o}lder's inequality, for every \(1\le p\le\infty\) and \(r\ge1/2\),
\begin{align}\label{fprimelower}
\|f'(re^{i\cdot})\|_{L^p(\partial\mathbb D)}
\ge m|a_m|r^{m-1}\geq  m|a_m|2^{1-m}>0.
\end{align}

{\bf Case 1.1.}  If \(1\le s<\infty\), then by \eqref{nablaandprime} and \eqref{fprimelower},
\[
\begin{aligned}
\|(1-|\cdot|^2)^t\widetilde{\nabla}f\|_{L^s(\mathbb D,d\lambda)}^s
&=
\int_{\mathbb D}
(1-|z|^2)^{s(t+1)}|f'(z)|^s\,d\lambda(z)\ge C
\int_{1/2}^1
(1-r^2)^{s(t+1)-2}\,dr.
\end{aligned}
\]
The last integral diverges if and only if
$
s(t+1)\le1.
$
Thus, if \(s(t+1)\le1\), every function in \(B_s^t\) is constant.
Conversely, if \(s(t+1)>1\), then we take \(f(z)=z\). We have
$
|\widetilde{\nabla}z|=1-|z|^2.
$
Then
\[
\|(1-|\cdot|^2)^t\widetilde{\nabla}z\|_{L^s(\mathbb D,d\lambda)}^s
=
\int_{\mathbb D}(1-|z|^2)^{s(t+1)}\,d\lambda(z)<\infty.
\]
Thus, \(z\in B_s^t\). This implies that \(B_s^t\ne\mathbb C\).

{\bf Case 1.2.} If \(s=q^{-}\) for some \(1<q<2\), then we recall the
argument in Case 1.1 that every
non-constant \(f\in H(\mathbb D)\) satisfies
$
(1-|\cdot|^2)^t\widetilde{\nabla}f\notin L^q(\mathbb D,d\lambda)$ if $q(t+1)\le1.
$
Since $
L^{q^-}(\mathbb D,d\lambda)\subset L^q(\mathbb D,d\lambda)
$, such an \(f\) cannot belong to \(B_{q^-}^t\). Hence,
\(B_{q^-}^t=\mathbb C\).
Conversely, if \(q(t+1)>1\), then for \(f(z)=z\),
\[
(1-|z|^2)^t|\widetilde{\nabla}z|
=
(1-|z|^2)^{t+1}\in L^{q^-}(\mathbb D,d\lambda).
\]
Thus, \(z\in B_{q^-}^t\). This implies that \(B_{q^-}^t\ne\mathbb C\).

{\bf Case 1.3.} If $s=\infty$, then we apply \eqref{nablaandprime} and  \eqref{fprimelower} to see that for every non-constant \(f\),
\[
\sup_{|z|=r}(1-|z|^2)^t|\widetilde{\nabla}f(z)|
\ge C(1-r^2)^{t+1},
\qquad r\ge1/2.
\]
If \(t<-1\), the right-hand side tends to \(+\infty\) as \(r\to1^{-}\).
Thus, every function in \(B_\infty^t\) is constant. If \(t\ge-1\), then \(f(z)=z\) satisfies
\[
(1-|z|^2)^t|\widetilde{\nabla}z|
=(1-|z|^2)^{t+1}\le C.
\]
Thus, \(z\in B_\infty^t\). This implies that \(B_\infty^t\ne\mathbb C\).

{\bf Case 2.} If \(n\ge2\), then we let \(d\sigma\) denote the normalized surface measure on
\(\mathbb S_n=\partial\mathbb B_n\). Let \(|\nabla_T f(z)|\) be the tangential gradient of \(f\) at \(z\) (see \cite[Section~7.6]{MR2115155}).
For \(z=r\zeta\), \(0<r<1\) and \(\zeta\in\mathbb S_n\), we have
\[
|\nabla_T f(z)|
:=
\max_{\substack{|u|=1\\ u\perp z}}
\left|\frac{\partial f}{\partial u}(z)\right|
=
\max_{\substack{|u|=1\\ u\perp z}}
\left|
\left\langle
u,\overline{\nabla f(z)}
\right\rangle
\right|
=\max_{\substack{|u|=1\\ u\perp \zeta}}
\left|
\left\langle
u,\overline{\nabla f(r\zeta)}
\right\rangle
\right|.
\]
Since this maximum equals the norm of the orthogonal projection of
\(\overline{\nabla f(r\zeta)}\) onto \(\zeta^\perp\),
\begin{align}\label{orthproj}
|\nabla_T f(r\zeta)|^2
=
\left|
\overline{\nabla f(r\zeta)}
-
\left\langle
\overline{\nabla f(r\zeta)},\zeta
\right\rangle\zeta
\right|^2
=
|\nabla f(r\zeta)|^2
-
\left|
\sum_{j=1}^n
\zeta_j\frac{\partial f}{\partial z_j}(r\zeta)
\right|^2 .
\end{align}
By \cite[Theorem~7.22]{MR2115155} (see also the fourth line in the proof of \cite[Theorem~7.23]{MR2115155}), for $z\in\mathbb{B}_n\backslash \{0\}$,
\begin{align}\label{invarianttangentlower}
|\widetilde{\nabla}f(z)|
\ge
(1-|z|^2)^{1/2}|\nabla_T f(z)|.
\end{align}

Let \(f\in H(\mathbb B_n)\) be non-constant. Write its homogeneous expansion as
$
f(z)=\sum_{k=0}^\infty P_k(z),
$
where each \(P_k\) is a homogeneous holomorphic polynomial of degree \(k\).
Let \(m\ge1\) be the smallest integer such that \(P_m\not\equiv0\).
We {\it claim} that
$
|\nabla_T P_m|\not\equiv0$ on $\mathbb{S}_n.$
Indeed, differentiating the identity \(P_m(\lambda z)=\lambda^mP_m(z)\) with respect to \(\lambda\) and then setting \(\lambda=1\) and $z=\zeta$ (see also \cite[(1.27)--(1.28)]{MR2115155}), we
obtain
\[
\sum_{j=1}^n \zeta_j\frac{\partial P_m}{\partial z_j}(\zeta)=mP_m(\zeta),
\qquad \zeta\in\mathbb S_n.
\]
Moreover, by the orthogonality of monomials on \(\mathbb S_n\) (see \cite[(1.21) and Lemma~1.11]{MR2115155}),
\[
\int_{\mathbb S_n}|\nabla P_m(\zeta)|^2\,d\sigma(\zeta)
=
m(n+m-1)
\int_{\mathbb S_n}|P_m(\zeta)|^2\,d\sigma(\zeta).
\]
This, together with \eqref{orthproj}, yields
\[
\begin{aligned}
\int_{\mathbb S_n}|\nabla_T P_m(\zeta)|^2\,d\sigma(\zeta)
&=
\int_{\mathbb S_n}|\nabla P_m(\zeta)|^2\,d\sigma(\zeta)
-
\int_{\mathbb S_n}
\left|
\sum_{j=1}^n \zeta_j\frac{\partial P_m}{\partial z_j}(\zeta)
\right|^2
d\sigma(\zeta)\\
&=
m(n-1)\int_{\mathbb S_n}|P_m(\zeta)|^2\,d\sigma(\zeta).
\end{aligned}
\]
Since \(n\ge2\), \(m\ge1\), and \(P_m\not\equiv0\), this shows that $
|\nabla_T P_m|\not\equiv0$ on $\mathbb{S}_n.$

Next we prove a lower estimate for the tangential gradient on spheres. Fix
\(0<r<1\), \(\zeta\in\mathbb S_n\), and a unit vector \(u\perp\zeta\). Since
\(u\perp e^{i\theta}\zeta\), the homogeneity of \(P_k\) gives
\[
\frac{\partial f}{\partial u}(re^{i\theta}\zeta)
=
\sum_{k=1}^\infty
r^{k-1}e^{i(k-1)\theta}
\frac{\partial P_k}{\partial u}(\zeta).
\]
Multiplying by \(e^{-i(m-1)\theta}\), and
integrating over \([0,2\pi]\), we see that
\[
r^{m-1}\frac{\partial P_m}{\partial u}(\zeta)
=
\frac1{2\pi}\int_0^{2\pi}
e^{-i(m-1)\theta}
\frac{\partial f}{\partial u}(re^{i\theta}\zeta)\,d\theta .
\]
Taking absolute values and then the supremum over all unit vectors
\(u\perp\zeta\), we obtain
\[
r^{m-1}|\nabla_T P_m(\zeta)|
\le
\frac1{2\pi}\int_0^{2\pi}
|\nabla_T f(re^{i\theta}\zeta)|\,d\theta .
\]
Taking the \(L^p(\mathbb S_n)\)-norm in \(\zeta\), applying Minkowski's inequality and the rotation invariance of \(d\sigma\), we see that for every
\(1\le p\le\infty\),
\begin{align}\label{tangentialLpLower}
r^{m-1}\|\nabla_T P_m\|_{L^p(\mathbb S_n)}
\le
\|\nabla_T f(r\cdot)\|_{L^p(\mathbb S_n)} .
\end{align}
Since \(|\nabla_T P_m|\not\equiv0\) on \(\mathbb S_n\), it follows that, for
\(r\ge1/2\),
\begin{align}\label{tangentialPositiveLower}
\|\nabla_T f(r\cdot)\|_{L^p(\mathbb S_n)}
\ge C_{f,p}>0 .
\end{align}

{\bf Case 2.1.} If \(1\le s<\infty\), then by \eqref{invarianttangentlower} and \eqref{tangentialPositiveLower},
\[
\begin{aligned}
\|(1-|\cdot|^2)^t\widetilde{\nabla}f\|_{L^s(\mathbb B_n,d\lambda)}^s&=2n
\int_0^1
\int_{\mathbb S_n}
(1-r^2)^{st}|\widetilde{\nabla}f(r\zeta)|^s
\,d\sigma(\zeta)
\frac{r^{2n-1}\,dr}{(1-r^2)^{n+1}}\\
&\ge 2n
\int_{1/2}^1
(1-r^2)^{s(t+1/2)}
\int_{\mathbb S_n}|\nabla_T f(r\zeta)|^s\,d\sigma(\zeta)
\frac{r^{2n-1}\,dr}{(1-r^2)^{n+1}}\\
&\ge
C\int_{1/2}^1
(1-r^2)^{s(t+1/2)-n-1}\,dr.
\end{aligned}
\]
The last integral diverges if and only if
$
s(t+1/2)\le n.
$
Hence, if \(s(t+1/2)\le n\), every function in \(B_s^t\) is constant.
Conversely, if \(s(t+1/2)>n\), then we take \(f(z)=z_1\). By \cite[Lemma~2.13]{MR2115155},
\begin{align}\label{zoneinvariant}
|\widetilde{\nabla}z_1(z)|^2
=
(1-|z|^2)(1-|z_1|^2)
\le
1-|z|^2.
\end{align}
Hence, $(1-|\cdot|^2)^t\widetilde{\nabla}z_1\in L^s(\mathbb{B}_n,d\lambda)$.
Thus, \(z_1\in B_s^t\). This implies that \(B_s^t\ne\mathbb C\).

{\bf Case 2.2.} If \(s=q^{-}\) for some \(1<q<2\), then we recall the
argument in Case 2.1 that every
non-constant \(f\in H(\mathbb B_n)\) satisfies
$
(1-|\cdot|^2)^t\widetilde{\nabla}f
\notin L^q(\mathbb B_n,d\lambda)
$ if \(q(t+1/2)\le n\).
Since $
L^{q^-}(\mathbb B_n,d\lambda)\subset L^q(\mathbb B_n,d\lambda)
$,  such an \(f\) cannot belong to \(B_{q^-}^t\). Hence,
\(B_{q^-}^t=\mathbb C\). Conversely, if \(q(t+1/2)>n\), then it follows from \eqref{zoneinvariant} that
$
(1-|z|^2)^t|\widetilde{\nabla}z_1(z)|
\in L^{q^-}(\mathbb B_n,d\lambda).
$
Thus, \(z_1\in B_{q^-}^t\). This implies that \(B_{q^-}^t\ne\mathbb C\).

{\bf Case 2.3.} If $s=\infty$, then we apply \eqref{invarianttangentlower} and \eqref{tangentialPositiveLower} to see that every non-constant
\(f\in H(\mathbb B_n)\) satisfies
\[
\sup_{|z|=r}
(1-|z|^2)^t|\widetilde{\nabla}f(z)|
\ge
C(1-r^2)^{t+1/2},
\qquad r\ge1/2.
\]
If \(t<-1/2\), the right-hand side tends to \(+\infty\) as \(r\to1^{-}\).
Thus, every function in \(B_\infty^t\) is constant. Conversely, if
\(t\ge-1/2\), then by \eqref{zoneinvariant},
\[
(1-|z|^2)^t|\widetilde{\nabla}z_1(z)|
\le
(1-|z|^2)^{t+1/2}\le C.
\]
Thus, \(z_1\in B_\infty^t\). This implies that \(B_\infty^t\ne\mathbb C\).

The proof of Lemma \ref{constantcharBesov} is
complete.
\end{proof}

The following statement says that for certain parameter ranges, the $r$-summing condition for $H^\beta_{\bar f}:A^p_\alpha\to L^q(dv_\beta)$ forces the holomorphic symbol $f$ to be constant.
\begin{coro}\label{mainhankelcconstant}
Let \(1\leq p<\infty\), \(1<q<\infty\), \(1\leq r<\infty\),
\(\alpha,\beta>-1\), and $f\in A^1_\beta$.
Assume that one of the following conditions holds:
\begin{enumerate}
\item \(\kappa\neq\infty\) and
\[
\begin{cases}
\kappa_\natural(\gamma+1)\leq 1,
& n=1,\\[1mm]
\kappa_\natural\left(\gamma+\dfrac12\right)\leq n,
& n\geq2,
\end{cases}
\]
where \(\kappa_\natural=\kappa\) if \(\kappa\in[1,\infty)\), and
\(\kappa_\natural=q\) if \(\kappa=q^{-}\).

\item \(\kappa=\infty\) and
\[
\begin{cases}
\gamma<-1,
& n=1,\\[1mm]
\gamma<-\dfrac12,
& n\geq2.
\end{cases}
\]
\end{enumerate}
Then
$H^\beta_{\bar{f}}\in \Pi_r(A^p_\alpha, L^q(dv_\beta))$ if and only if
$f$ is constant.
Moreover, the above parameter conditions are sharp.
\end{coro}
\begin{proof}
The required equivalence follows directly by combining Lemma \ref{constantcharBesov} with Corollary \ref{mainhankelc}.
It remains to verify sharpness. Suppose that the parameter condition stated above for the given values of \(n\) and \(\kappa\) is not satisfied. The proof of Lemma~\ref{constantcharBesov} shows that
\(z\in B^\gamma_\kappa\) when \(n=1\), and that
\(z_1\in B^\gamma_\kappa\) when \(n\geq2\).
These coordinate functions belong to \(A^1_\beta\) and are
non-constant. Corollary~\ref{mainhankelc} therefore implies that the
corresponding Hankel operator is \(r\)-summing. Hence the conclusion
fails outside the stated parameter range.
\end{proof}

\section{Absolutely summing little Hankel operators}\label{Sect6}
\setcounter{equation}{0}

In this section, we characterize those symbols $f$ for which the associated little Hankel operator $h^\beta_{\bar{f}}:A^p_\alpha\to L^q(dv_\beta)$ is $r$-summing.

The following proposition gives a slight variant of \cite[Proposition 8.37]{MR2311536}: the a priori assumption on the symbol is weakened, and the reduction \(h_{\bar f}^\beta=h_{\overline{P_\beta f}}^\beta\) is replaced by $ h^\beta_{\bar f}=h^\beta_{\overline{V_{b,\beta}(f)}}. $ This variant is what allows Theorem \ref{Maintheoremlittle} to be proved with only the natural \(L^1(dv_\beta)\) a priori hypothesis on the symbol \(f\).
\begin{proposition}\label{p61}
Let \(\beta>-1\), \(b>0\), and \(f\in L^1(dv_\beta)\). Suppose that
\(V_{b,\beta}(f)\in L^1(dv_\beta)\). Then
\begin{align}\label{littleequa}
h^\beta_{\bar f}(g)=h^\beta_{\overline{V_{b,\beta}(f)}}(g),
\qquad g\in H^\infty.
\end{align}
\end{proposition}
\begin{proof}
Set \(F:=V_{b,\beta}(f)\) and
$
u(z):=(1-|z|^2)^{-b}f(z).
$
Then \(u\in L^1(dv_{\beta+b})\) and
$
F=(1-|\cdot|^2)^b P_{\beta+b}u.
$
Since \(F\in L^1(dv_\beta)\), we also have
\(P_{\beta+b}u\in L^1(dv_{\beta+b})\). Hence, for every
\(\varphi\in H^\infty\),
\begin{align}\label{varphiequa}
\int_{\mathbb B_n} F(z)\overline{\varphi(z)}\,dv_\beta(z)
&=
\frac{c_\beta}{c_{\beta+b}}
\int_{\mathbb B_n} P_{\beta+b}u(z)\overline{\varphi(z)}\,dv_{\beta+b}(z)\nonumber\\
&=
\frac{c_\beta}{c_{\beta+b}}
\int_{\mathbb B_n} u(z)\overline{\varphi(z)}\,dv_{\beta+b}(z)\nonumber\\
&=
\int_{\mathbb B_n} f(z)\overline{\varphi(z)}\,dv_\beta(z).
\end{align}

Now fix \(g\in H^\infty\) and \(z\in\mathbb B_n\). The function
\[
\varphi_z(w):=\frac{g(w)}{(1-\langle w,z\rangle)^{n+1+\beta}}
\]
belongs to \(H^\infty\). Applying the identity \eqref{varphiequa} to \(\varphi_z\) and
then conjugating gives
\[
\int_{\mathbb B_n}\frac{\overline{f(w)}g(w)}
{(1-\langle w,z\rangle)^{n+1+\beta}}\,dv_\beta(w)
=
\int_{\mathbb B_n}\frac{\overline{F(w)}g(w)}
{(1-\langle w,z\rangle)^{n+1+\beta}}\,dv_\beta(w).
\]
That is,
$
h^\beta_{\bar f}(g)(z)
=
h^\beta_{\overline{V_{b,\beta}(f)}}(g)(z).
$
The proof of Proposition \ref{p61} is complete.
\end{proof}

%

%
%

Recall that \(\mathfrak L_\delta\) denotes the collection of all
\(\delta\)-lattices in \(\mathbb B_n\). The following lemma establishes a discretization principle.
\begin{lemma}\label{kl}
Let \(0<\delta\leq1\), \(\beta>-1\), \(b>0\), \(t\in\mathbb R\), and
\(f\in L^1(dv_\beta)\). Let \(\Phi\) be a Young function. Then
$
(1-|\cdot|^2)^t |V_{b,\beta}(f)|
\in L_\Phi(\mathbb B_n,d\lambda)
$
if and only if
$
\left\{
(1-|a_k|^2)^t |V_{b,\beta}(f)(a_k)|
\right\}_{k\geq1}
\in\ell_\Phi
$
for every \(\{a_k\}\in\mathfrak L_\delta\). Moreover,
\[
\left\|
(1-|\cdot|^2)^t |V_{b,\beta}(f)|
\right\|_{L_\Phi(\mathbb B_n,d\lambda)}
\simeq
\sup_{\{a_k\}\in\mathfrak L_\delta}
\left\|
\left\{
(1-|a_k|^2)^t |V_{b,\beta}(f)(a_k)|
\right\}_{k\geq1}
\right\|_{\ell_\Phi}.
\]
An analogous statement holds with
\(L_\Phi(\mathbb B_n,d\lambda)\) and \(\ell_\Phi\) replaced by
\(L^\infty(\mathbb B_n,d\lambda)\) and \(\ell^\infty\), respectively.
\end{lemma}
\begin{proof}
For $z\in\mathbb B_n$, we let
\[
F_{b,\beta}(z):=\int_{\mathbb{B}_n}
\frac{f(w)}{(1-\langle z,w\rangle)^{n+1+\beta+b}}\,dv_\beta(w),\qquad
G_{b,\beta,t}(z):=(1-|z|^2)^t |V_{b,\beta}(f)(z)|.
\]
Since \(f\in L^1(dv_\beta)\), we have
\(F_{b,\beta}\in H(\mathbb B_n)\), and
\begin{align}\label{Vbbeta}
V_{b,\beta}(f)(z)=\frac{c_{b+\beta}}{c_\beta}
(1-|z|^2)^b F_{b,\beta}(z).
\end{align}

{\it (Sufficiency.)} Assume that $
\{
(1-|a_k|^2)^t |V_{b,\beta}(f)(a_k)|
\}_{k\geq1}
\in\ell_\Phi
$
for every \(\{a_k\}\in\mathfrak L_\delta\). Fix a \(\delta\)-lattice \(\{c_k\}_{k\geq1}\). Since \(F_{b,\beta}\in H(\mathbb B_n)\), we apply \eqref{Vbbeta} to see that $G_{b,\beta,t}(z)$ is continuous on \(\mathbb B_n\). Therefore, for each \(k\), there exists
\(\xi_k\in\overline{B(c_k,\delta)}\) such that
\[
G_{b,\beta,t}(\xi_k)
=
\max_{z\in\overline{B(c_k,\delta)}}
G_{b,\beta,t}(z).
\]

Since \(\xi_k\in\overline{B(c_k,\delta)}\), the finite multiplicity property of \(\{c_k\}\) implies that the sequence
\(\{\xi_k\}\) can be partitioned into
at most \(m=m(n,\delta)\) \(\delta\)-separated subsequences.
Thus, there are finitely many index sets \(I_1,\ldots,I_m\) such that
$
\mathbb N=I_1\cup\cdots\cup I_m
$
and each sequence \(\{\xi_k:k\in I_i\}\) is $\delta$-separated. For each $i$, extend
\(\{\xi_k:k\in I_i\}\) to a maximal \(\delta\)-separated sequence
\(\mathcal L_i\). By maximality,
$
\{B(\xi,\delta):\xi\in\mathcal L_i\}
$
covers \(\mathbb B_n\), while the balls
$
\{B(\xi,\frac{\delta}{4}):\xi\in\mathcal L_i\}
$
are pairwise disjoint. Thus \(\mathcal L_i\) is a
\(\delta\)-lattice. By the assumption, we have
\[
\|
\{G_{b,\beta,t}(\xi_k)\}_{k\geq1}
\|_{\ell_\Phi}
\leq
\sum_{i=1}^{m}
\|
\{G_{b,\beta,t}(\xi):\xi\in\mathcal L_i\}
\|_{\ell_\Phi}
<\infty.
\]

For the statement regarding the Luxemburg norm, we let
$
\eta>
\|
\{
G_{b,\beta,t}(\xi_k)
\}_{k\geq1}
\|_{\ell_\Phi}.
$
By the definition of the Luxemburg norm,
\[
\sum_k
\Phi\left(
\frac{G_{b,\beta,t}(\xi_k)}{\eta}
\right)
\leq1.
\]
Since the balls \(\{B(c_k,\delta)\}\) cover \(\mathbb B_n\),
 and \(d\lambda\) is M\"obius-invariant, we obtain
\begin{align*}
\int_{\mathbb B_n}
\Phi\left(
\frac{G_{b,\beta,t}(z)}{C\eta}
\right)
\,d\lambda(z)
\leq
\sum_k
\int_{B(c_k,\delta)}
\Phi\left(
\frac{G_{b,\beta,t}(z)}{C\eta}
\right)
\,d\lambda(z)
\lesssim
\sum_k
\Phi\left(
\frac{G_{b,\beta,t}(\xi_k)}{C\eta}
\right)\leq 1
\end{align*}
for a sufficiently large constant \(C\geq1\).

It follows from the definition of the Luxemburg norm that
$
\|
G_{b,\beta,t}
\|_{L_\Phi(\mathbb B_n,d\lambda)}
\leq C\eta.
$
Letting
$
\eta\downarrow
\|
\{
G_{b,\beta,t}(\xi_k)
\}_{k\geq1}
\|_{\ell_\Phi},
$
we obtain
\begin{equation}\label{eq:kl-reverse-sampling}
\|
G_{b,\beta,t}
\|_{L_\Phi(\mathbb B_n,d\lambda)}
\lesssim \|
\{
G_{b,\beta,t}(\xi_k)
\}_{k\geq1}
\|_{\ell_\Phi}\lesssim \sum_{i=1}^{m}
\|
\{
G_{b,\beta,t}(\xi):\xi\in\mathcal L_i
\}
\|_{\ell_\Phi}\lesssim\sup_{\{a_k\}\in\mathfrak L_\delta}
\|
\{
G_{b,\beta,t}(a_k)
\}_{k\geq1}
\|_{\ell_\Phi}
.
\end{equation}
This proves both the qualitative implication and the asserted reverse
norm estimate.

For the statement regarding the $L^\infty$ norm, assume that
$
\{
(1-|a_k|^2)^t |V_{b,\beta}(f)(a_k)|
\}_{k\geq1}
\in\ell^\infty
$
for every \(\{a_k\}\in\mathfrak L_\delta\).
The covering property yields
\[
\left\|
G_{b,\beta,t}
\right\|_{L^\infty(\mathbb B_n,d\lambda)}
\leq
\|
\{
G_{b,\beta,t}(\xi_k)
\}_{k\geq1}
\|_{\ell^\infty}\leq
\max_{1\leq i\leq m}
\|
\{
G_{b,\beta,t}(\xi):\xi\in\mathcal L_i
\}
\|_{\ell^\infty}\leq \sup_{\{a_k\}\in\mathfrak L_\delta}
\|
\{
G_{b,\beta,t}(a_k)
\}_{k\geq1}
\|_{\ell^\infty}.
\]
This proves both the qualitative implication and the asserted reverse
norm estimate.

{\it (Necessity.)} Assume that $
(1-|\cdot|^2)^t |V_{b,\beta}(f)|
\in L_\Phi(\mathbb B_n,d\lambda)
$. Fix a \(\delta\)-lattice
\(\{a_k\}_{k\geq1}\).

For the statement regarding the Luxemburg norm,
we let $
\eta>\|G_{b,\beta,t}\|_{L_\Phi(\mathbb B_n,d\lambda)}
$. By the definition of the Luxemburg norm,
\[
\int_{\mathbb B_n}
\Phi\left(\frac{G_{b,\beta,t}(w)}{\eta}\right)\,d\lambda(w)\le 1.
\]
We note that for each $k$, the function
\[
z\mapsto
\Phi\left(
\frac{(1-|a_k|^2)^{t+b}|F_{b,\beta}(z)|}{\eta}
\right)
\]
is plurisubharmonic on \(\mathbb B_n\) (see, e.g., \cite[Theorem~5.6, p.~39, and Example~5.12, p.~41]{DemaillyCADG}). Hence, by the sub-mean property for plurisubharmonic functions (see \cite[(5.3)]{DemaillyCADG} and \cite[Theorem and Definition~4.12]{DemaillyCADG}) and \eqref{a3}, there exists a constant \(C_1>0\), independent of \(k\) and \(\eta\), such that
\[
\Phi\left(\frac{G_{b,\beta,t}(a_k)}{C_1\eta}\right)
\lesssim
\int_{B(a_k,\delta/2)}
\Phi\left(\frac{G_{b,\beta,t}(w)}{\eta}\right)\,d\lambda(w).
\]
By the finite multiplicity property
\eqref{a1}, the family \(\{B(a_k,\delta/2)\}\) has finite overlap. Consequently,
\begin{align*}
\sum_k
\Phi\left(\frac{G_{b,\beta,t}(a_k)}{C_1\eta}\right)
&\lesssim
\sum_k
\int_{B(a_k,\delta/2)}
\Phi\left(\frac{G_{b,\beta,t}(w)}{\eta}\right)\,d\lambda(w) \lesssim
\int_{\mathbb B_n}
\Phi\left(\frac{G_{b,\beta,t}(w)}{\eta}\right)\,d\lambda(w)
\le C_2
\end{align*}
for some constant \(C_2\ge 1\). Since \(\Phi\) is convex and \(\Phi(0)=0\),
\[
\sum_k
\Phi\left(\frac{G_{b,\beta,t}(a_k)}{C_1C_2\eta}\right)
\le
\frac{1}{C_2}
\sum_k
\Phi\left(\frac{G_{b,\beta,t}(a_k)}{C_1\eta}\right)
\le 1.
\]
By the definition of the Luxemburg norm,
$
\|
\{
G_{b,\beta,t}(a_k)
\}
\|_{\ell_\Phi}
\le C_1C_2\eta .
$
Letting
$
\eta\downarrow \|G_{b,\beta,t}\|_{L_\Phi(\mathbb B_n,d\lambda)}$,
we deduce that
\begin{align}
\sup_{\{a_k\}\in\mathfrak L_\delta}
\|
\{
G_{b,\beta,t}(a_k)
\}_{k\geq1}
\|_{\ell_\Phi}
\lesssim
\|
G_{b,\beta,t}
\|_{L_\Phi(\mathbb B_n,d\lambda)}.
\label{eq:kl-sampling-upper}
\end{align}

For the statement regarding the $L^\infty$ norm, assume that
$
(1-|\cdot|^2)^t |V_{b,\beta}(f)|
\in L^\infty(\mathbb B_n,d\lambda).
$
We apply the continuity of \(G_{b,\beta,t}\) to see that
\[
\sup_{\{a_k\}\in\mathfrak L_\delta}
\|\{G_{b,\beta,t}(a_k)\}_{k\ge1}\|_{\ell^\infty}
\le
\sup_{z\in\mathbb B_n}G_{b,\beta,t}(z)
=
\|G_{b,\beta,t}\|_{L^\infty(\mathbb B_n,d\lambda)}.
\]

The proof of Lemma \ref{kl} is complete.
\end{proof}

\begin{proof}[Proof of Theorem \ref{Maintheoremlittle}]
Since the first statement is independent of the choice of $b$ and $\delta$, it suffices to show the equivalence for fixed $b>n+(\alpha+1)/p$ and \(0<\delta\leq1\).

\eqref{littlehankel3} $\Rightarrow$ \eqref{littlehankel1}.
	Set
	$
	d\mu(z)=|V_{b,\beta}(f)(z)|^q\,dv_\beta(z).
	$
Let \(\{a_k\}\) be a \(\delta\)-lattice. By \eqref{a2} and \eqref{a3}, we have
	\begin{align}\label{muestimates}
	\hat{\mu}_{p,q,\alpha,\delta}(a_k)
	&\simeq
	(1-|a_k|^2)^{-\frac{n+1+\alpha}{p}}
	\left(
	\int_{B(a_k,\delta)}
	|V_{b,\beta}(f)(z)|^q\,dv_\beta(z)
	\right)^{1/q}\nonumber\\
    &\simeq
\Bigg(\int_{B(a_k,\delta)}
(1-|z|^2)^{q\gamma}
|V_{b,\beta}(f)(z)|^{q}
d\lambda(z)\Bigg)^{1/q}\nonumber\\
&\lesssim
\max_{z\in \overline{B(a_k,\delta)}}
(1-|z|^2)^\gamma
|V_{b,\beta}(f)(z)|.
	\end{align}
For each \(k\), choose \(\xi_k\in\overline{B(a_k,\delta)}\) such that
\[
(1-|\xi_k|^2)^\gamma
|V_{b,\beta}(f)(\xi_k)|
=
\max_{z\in\overline{B(a_k,\delta)}}
(1-|z|^2)^\gamma
|V_{b,\beta}(f)(z)|.
\]
	By the finite multiplicity property \eqref{a1}, the sequence \(\{\xi_k\}\) is a finite union of separated sequences. Thus, there are finitely many index sets \(I_1,\ldots,I_m\) such that
$
	\mathbb N=I_1\cup\cdots\cup I_m
$
	and each sequence \(\{\xi_k:k\in I_i\}\) is $\delta$-separated. For each $i$, extend
\(\{\xi_k:k\in I_i\}\) to a maximal \(\delta\)-separated sequence
\(\mathcal L_i\). By maximality,
$
\{B(\xi,\delta):\xi\in\mathcal L_i\}
$
covers \(\mathbb B_n\), while the balls
$
\{B(\xi,\frac{\delta}{4}):\xi\in\mathcal L_i\}
$
are pairwise disjoint. Thus \(\mathcal L_i\) is a
\(\delta\)-lattice.

Applying Lemma \ref{kl}, with lattice radius \(\delta\), gives
\begin{align}\label{eq:lemma-kl-applied-perturbed}
\left\|
\left\{
(1-|\xi|^2)^\gamma
|V_{b,\beta}(f)(\xi)|
\right\}_{\xi\in\mathcal L_i}
\right\|_{\ell^\kappa}\lesssim
\left\|
(1-|\cdot|^2)^\gamma
|V_{b,\beta}(f)|
\right\|_{L^\kappa(\mathbb B_n,d\lambda)}.
\end{align}
Since \(\{\xi_k:k\in I_i\}\) is a subsequence of
\(\mathcal L_i\) and \(\mathbb N=I_1\cup\cdots\cup I_m\), we obtain
\begin{align}\label{finallylittle}
\left\|
\left\{
\max_{z\in \overline{B(a_k,\delta)}}
(1-|z|^2)^\gamma
|V_{b,\beta}(f)(z)|
\right\}_{k\geq1}
\right\|_{\ell^\kappa}
&=
\left\|
\left\{
(1-|\xi_k|^2)^\gamma
|V_{b,\beta}(f)(\xi_k)|
\right\}_{k\geq1}
\right\|_{\ell^\kappa}
\nonumber\\
&\leq
\sum_{i=1}^{m}
\left\|
\left\{
(1-|\xi_k|^2)^\gamma
|V_{b,\beta}(f)(\xi_k)|
\right\}_{k\in I_i}
\right\|_{\ell^\kappa}
\nonumber\\
&\lesssim
\left\|
(1-|\cdot|^2)^\gamma
|V_{b,\beta}(f)|
\right\|_{L^\kappa(\mathbb B_n,d\lambda)}.
\end{align}
Combining the estimates \eqref{muestimates} and
\eqref{finallylittle}, we deduce that
\begin{align*}
\|\{\hat{\mu}_{p,q,\alpha,\delta}(a_k)\}\|_{\ell^\kappa}
\lesssim
\left\|
\left\{
\max_{z\in \overline{B(a_k,\delta)}}
(1-|z|^2)^\gamma
|V_{b,\beta}(f)(z)|
\right\}_{k\geq1}
\right\|_{\ell^\kappa}\lesssim
\left\|
(1-|\cdot|^2)^\gamma
|V_{b,\beta}(f)|
\right\|_{L^\kappa(\mathbb B_n,d\lambda)}.
\end{align*}
Combining with Theorem \ref{MainTheorem2}, we see that \(J_\mu:A^p_\alpha\to L^q(d\mu)\) is \(r\)-summing and
\begin{align}\label{hatinequality}
\pi_r(J_\mu)
	\lesssim
	\|\{\hat{\mu}_{p,q,\alpha,\delta}(a_k)\}\|_{\ell^{\kappa}}
	\lesssim
	\|(1-|\cdot|^2)^\gamma |V_{b,\beta}(f)|\|_{L^{\kappa}(\mathbb B_n,d\lambda)}.
\end{align}
	

Since \(1\in A_\alpha^p\), the boundedness of \(J_\mu\) implies that $J_\mu(1)\in L^q(d\mu)$. Thus,
$
V_{b,\beta}(f)\in L^q(dv_\beta)\subset L^1(dv_\beta).
$
By Proposition \ref{p61},
\begin{equation}\label{eq52}
h^\beta_{\bar f}(g)
=
h^\beta_{\overline{V_{b,\beta}(f)}}(g),
\qquad g\in H^\infty.
\end{equation}

	By the boundedness of $\overline{P_\beta}:L^q(dv_\beta) \to \overline{A^q_\beta}$ and the formula \eqref{eq52}, we have
	\[
	\|h^\beta_{\bar{f}}(g)\|_{q,\beta}
	=
	\|h^\beta_{\overline{V_{b,\beta}(f)}}(g)\|_{q,\beta}
	=
	\|\overline{P_\beta}(\overline{V_{b,\beta}(f)}g)\|_{q,\beta}
	\lesssim
	\|\overline{V_{b,\beta}(f)}g\|_{q,\beta}
	=
	\|J_\mu (g)\|_{L^q(d\mu)},
	\]
	for $g\in H^\infty$. Since \(H^\infty\) is dense in \(A^p_\alpha\) and
\(J_\mu:A^p_\alpha\to L^q(d\mu)\) is bounded, \(h^\beta_{\bar f}\)
extends uniquely to a bounded operator on \(A^p_\alpha\), and the above
estimate remains valid for every \(g\in A^p_\alpha\). Consequently, by the definition of \(r\)-summing operators and \eqref{hatinequality}, we see that \(h^\beta_{\bar f}\) is \(r\)-summing and
	\[
	\pi_r(h^\beta_{\bar f})\lesssim \pi_r(J_\mu)\lesssim
	\|(1-|\cdot|^2)^\gamma |V_{b,\beta}(f)|\|_{L^{\kappa}(\mathbb B_n,d\lambda)}.
	\]
	
\eqref{littlehankel1} $\Rightarrow$ \eqref{littlehankel3}. By \eqref{a4s}, for every $z\in\mathbb{B}_n$,
	\begin{equation*}
		\begin{aligned}
			|h^\beta_{\bar{f}}(h_{p,z}^b)(z)|&=\left|\int_{\mathbb{B}_n} \frac{\overline{f(w)} h_{p,z}^b(w)}{(1-\langle w,z\rangle )^{n+1+\beta}} \, dv_\beta(w)\right|\\
			&\simeq \left|\int_{\mathbb{B}_n} \frac{(1-|z|^2)^{b-(n+1+\alpha)/p}\overline{f(w)}}{(1-\langle w,z\rangle )^{n+1+\beta+b}} \, dv_\beta(w)\right|\\
			&\simeq (1-|z|^2)^{-(n+1+\alpha)/p}|V_{b,\beta} (f)(z)|.
		\end{aligned}
	\end{equation*}
	On the other hand, by Lemma \ref{Bergmaninequality}, we have
	\begin{align*}
		|h^\beta_{\bar{f}}(h_{p,z}^b)(z)|^q\lesssim \frac{1}{v_\beta(B(z,\delta))} \int_{B(z,\delta)}|h^\beta_{\bar{f}}(h_{p,z}^b)(w)|^q\,dv_\beta(w).
	\end{align*}
	Therefore, by \eqref{a2} we have
	\begin{align}\label{doubleuse}
		(1-|z|^2)^{\gamma}|V_{b,\beta} (f)(z)|\lesssim \left( \int_{B(z,\delta)}|h^\beta_{\bar{f}}(h_{p,z}^b)(w)|^q\,dv_\beta(w)\right)^{\frac{1}{q}}.
	\end{align}

Let
\(\{a_k\}\) be an arbitrary \(\delta\)-lattice in \(\mathbb B_n\). Set \(\mathbf d=\{d_k\}\) with
\[
d_k
:=
\left(
\int_{B(a_k,\delta)}
\left|
h^\beta_{\bar f}(h_{p,a_k}^b)(w)
\right|^q
\,dv_\beta(w)
\right)^{1/q}.
\]

Applying Proposition~\ref{proq2pr} with
$T=h^\beta_{\bar f}$ and $d\mu=dv_\beta$, we see that $\mathscr M_{\mathbf d}\in\Pi_r(\ell^p,\ell^q)$ and satisfies
$
\pi_r(\mathscr M_{\mathbf d})
\lesssim
\pi_r(h^\beta_{\bar f}).
$
Combining \eqref{doubleuse} with Theorem~\ref{GarlingTheorem}, we deduce that
\begin{equation}\label{eq:little-hankel-d-sequence}
\left\|
\left\{
(1-|a_k|^2)^\gamma
|V_{b,\beta}(f)(a_k)|
\right\}_{k\geq1}
\right\|_{\ell^\kappa}\lesssim\|\mathbf d\|_{\ell^\kappa}
\lesssim
\pi_r(\mathscr M_{\mathbf d})
\lesssim
\pi_r(h^\beta_{\bar f}).
\end{equation}
This, in combination with Lemma~\ref{kl}, yields
\begin{equation}\label{eq:little-hankel-necessity-final}
\left\|
(1-|\cdot|^2)^\gamma
|V_{b,\beta}(f)|
\right\|_{L^\kappa(\mathbb B_n,d\lambda)}
\lesssim \sup_{\{a_k\}\in\mathfrak L_\delta}
\left\|
\left\{
(1-|a_k|^2)^\gamma
|V_{b,\beta}(f)(a_k)|
\right\}_{k\geq1}
\right\|_{\ell^\kappa}
\lesssim
\pi_r(h^\beta_{\bar f}).
\end{equation}

\eqref{littlehankel3} $\Leftrightarrow$ \eqref{littlehankel2}. This follows directly from Lemma \ref{kl}.

The proof of Theorem \ref{Maintheoremlittle} is complete.
\end{proof}

	\bigskip

\noindent
 {\bf Acknowledgements:}
Z. Fan is supported by the National Natural Science Foundation of China  (No. 12601166), and by the Guangdong Basic and Applied Basic Research Foundation (No. 2023A1515110879). B. He is supported by the National Natural Science Foundation of China  (No. 12501155). X. Wang is supported by  the National Natural Science Foundation of China (Grant No. 12471119 and No. 11971125).

\end{document}